\documentclass{amsart}
\usepackage{geometry}
\usepackage{doi}
\usepackage{graphicx}
\usepackage[numbers]{natbib}
\usepackage{amssymb}
\usepackage{amsthm}
\usepackage{amsmath}
\usepackage{mathtools}

\newtheorem{definition}{Definition}
 
\usepackage{multicol}
\usepackage{array}
\usepackage{boldline}
\usepackage{hyperref}
\usepackage[nosort]{cleveref}
\usepackage{tikz}
\usepackage{subfigure}
\usetikzlibrary{shapes.geometric, arrows}
\usetikzlibrary{er,positioning}
\usetikzlibrary{shapes}
\usepackage[font=small,skip=5pt]{caption}

\newcommand*\circled[1]{\tikz[baseline=(char.base)]{%
		\node[shape=circle,fill=white!80!black,draw,minimum size=20] (char) {#1};}}

\usepackage{algorithm}
\usepackage{algpseudocode}

\algnewcommand\True{\textbf{True}\space}
\algnewcommand\False{\textbf{False}\space}
\algloopdefx{NFor}[1]{\textbf{for all} #1 \textbf{do}}
\algloopdefx{NIf}[1]{\textbf{if} #1 \textbf{then}}

\usepackage{multirow}
\usepackage{lscape}
\usepackage{adjustbox}

\DeclareMathOperator{\dom}{dom}
\DeclareMathOperator{\opt}{opt}
\DeclareMathOperator*{\argmax}{arg\,max}

\newcommand{\linprogref}[1]{\textup{(#1)}}

\title[Improving algorithms for $\Gamma$-maximin, $\Gamma$-maximax and interval dominance]{Improving and benchmarking of algorithms for $\Gamma$-maximin, $\Gamma$-maximax and interval dominance}

\author{Nawapon Nakharutai}
\address{Chiang Mai University, 239 Huaykaew Rd, Suthep, Mueang District, Chiang Mai 50200, Thailand}
\email{nawapon.nakharutai@cmu.ac.th}

\author{Matthias C. M. Troffaes}
\address{Durham University, Stockton Rd, Durham DH1 3LE, United Kingdom}
\email{matthias.troffaes@durham.ac.uk}

\author{Camila C. S. Caiado}
\address{Durham University, Stockton Rd, Durham DH1 3LE, United Kingdom}
\email{c.c.s.caiado@durham.ac.uk}

\keywords{$\Gamma$-maximin; $\Gamma$-maximax; interval dominance; primal-dual; algorithm; benchmarking}

\begin{document}

\begin{abstract}
 $\Gamma$-maximin, $\Gamma$-maximax and inteval dominance are familiar decision criteria for making decisions under severe uncertainty, when probability distributions can only be partially identified. One can apply these three criteria by solving sequences of linear programs. In this study, we present new algorithms for these criteria and compare their performance to existing standard algorithms. 
Specifically, we use efficient ways, based on previous work, to find common initial feasible points for these algorithms. Exploiting these initial feasible points, we develop early stopping criteria to determine whether gambles are either $\Gamma$-maximin, $\Gamma$-maximax or interval dominant. We observe that the primal-dual interior point method benefits considerably from these improvements.
In our simulation, we find that our proposed algorithms outperform the standard algorithms when the size of the domain of lower previsions is less or equal to the sizes of decisions and outcomes. However,  our proposed algorithms do not outperform the standard algorithms in the case that the size of the domain of lower previsions is much larger than the sizes of decisions and outcomes.
\end{abstract}

\maketitle

\section{Introduction}\label{intro}

Consider a situation where a subject must make a decision by selecting any options from a set of possible decisions. Each option leads to an uncertain reward depending on the option and on some uncertain outcome. The reward can be anything, for example, money or casino chips.  However, we suppose that rewards can be related to a utility scale, and hence any uncertain reward can be considered as a bounded real-valued function on a set of outcomes. We call such function a \emph{gamble}. In other words, we now consider a situation where a subject would like to choose gambles from a set of possible gambles.

We suppose that, for any set of possible gambles, the subject is able to state gambles that she does not want to select. Once we delete all gambles that she does not choose, the remaining gambles in the set are called \emph{optimal}. In a classical decision theory, the subject should simply select a gamble that yields the maximum expected utility if probabilities for all relevant events are available for her \citep{anscombe:1963}. Unfortunately, she may be unable to specify all probabilities exactly. To handle this situation, \citet{1991:walley} suggested that the subject can still model her belief using \emph{lower previsions} which are equivalent to expectation bounds.

There are several well-known decision criteria associated with lower previsions \cite{2007:troffaes:decision:intro}. For example, \textit{E-admissibility}, \textit{maximality} and \textit{interval dominance} will return a set of decisions \cite{2007:troffaes:decision:intro}. In this study, we are interested in interval dominance as well as two other popular criteria which give us a set of equivalent decisions (usually a singleton for generic decision problems), that is, $\Gamma$-maximin which maximizes the lower prevision of gambles, and $\Gamma$-maximax which maximizes the upper prevision of gambles \citep[p.~193]{2014:troffaes:itip:decision}.

Algorithms for decision criteria associated with lower previsions have been studied in the literature for quite some time. For example, see \citep{Utkin2005PowerfulAF} for E-admissibility and \citep{2011:Kikuti:Cozman:Filho} \citep{MATT201763}, \citep{2019:Nakharutai:Troffase:Caiado:maximal} and \citep[p.~336]{2014:troffaes:itip:computation} for maximality. From these, only \citep{2019:Nakharutai:Troffase:Caiado:maximal} provides in-depth benchmarking of their algorithms.

Algorithms for $\Gamma$-maximin and $\Gamma$-maximax have been studied for instance in  
\citet{2011:Kikuti:Cozman:Filho} and \citet[p.~335]{2014:troffaes:itip:computation}.
In those algorithms, to determine whether a gamble is $\Gamma$-maximin (or $\Gamma$-maximax), 
one has to evaluate the value of the lower (or the upper) prevision of all gambles in a set of gambles. This can be done by solving linear programs \citep{2011:Kikuti:Cozman:Filho,2014:troffaes:itip:computation}. Therefore, the number of linear programs that we have to solve is the cardinality of the set of gambles.
We emphasize that the scope of our study only concerns lower previsions defined on a finite set of gambles over a finite possibility space, and only concerns finite sets of decisions. More general situations where, say, non-linear relationships are imposed to construct the credal set (which can happen under certain independence conditions), and thus where non-linear programming applies, are beyond the scope of this paper.

In previous work, we studied efficient ways to solve linear programs inside algorithms for decision making with lower prevision. Specifically, in \citep{2017:Nakharutai:Troffaes:Caiado,2018:Nakharutai:Troffaes:Caiado} we investigated efficient algorithms to solve linear programs for checking avoiding sure loss. We presented an efficient way to find initial feasible points \citep[\S 4.2]{2018:Nakharutai:Troffaes:Caiado}, and extra stopping criteria to check avoiding sure loss. The results in that work showed that the primal-dual method can exploit these improvements and performs much more better compared to the simplex and the affine scaling methods when working with lower previsions.
Since we knew that the primal-dual method is suitable for handling lower previsions, in \citep{2019:Nakharutai:Troffase:Caiado:maximal} we improved and benchmarked algorithms for maximality, also using the primal-dual method. Some improvements, such as efficient ways to access common feasible starting points or early stopping criteria, also have clear potential to improve algorithms for $\Gamma$-maximin, $\Gamma$-maximax and interval dominance. This is what we aim to do in this paper.

The contributions of this paper are as follows. Based on our earlier work \citep{2017:Nakharutai:Troffaes:Caiado,2018:Nakharutai:Troffaes:Caiado,2019:Nakharutai:Troffase:Caiado:maximal}, we propose new algorithms for $\Gamma$-maximin, $\Gamma$-maximax and interval dominance. Specifically, we implement a quick way to obtain feasible starting points for the linear programs and early stopping criteria to determine whether a gamble is $\Gamma$-maximin (or $\Gamma$-maximax). As the primal-dual method solves both the primal and dual at the same time, we show how it can exploit both feasible starting points and early stopping criteria to further improve the efficiency of the algorithms for $\Gamma$-maximin and $\Gamma$-maximax. We also show that these improvements can be applied to algorithms for interval dominance. In addition, we benchmark these improvements by doing a comparative study of different improved algorithms on random generated sets of gambles for $\Gamma$-maximin, $\Gamma$-maximax and interval dominance.

We organise the paper as follows. In \cref{sec:decision}, we review the theory of lower previsions and natural extension, and we review the definitions of $\Gamma$-maximin, $\Gamma$-maximax and interval dominance. In \cref{sec:algorithm}, we present algorithms for  $\Gamma$-maximin, $\Gamma$-maximax and interval dominance and propose improvements for these algorithms. For benchmarking these improvements, we use an algorithm from \citep{2020:IUKM:Nakharutai} and \citep{2017:Nakharutai:Troffaes:Caiado} to generate sets of gambles, and compare the performance of algorithms with different improvements in \cref{sec:benchmark}. Finally, a conclusion is presented in \cref{sec:conclusion}.

\section{Decision making with lower previsions}\label{sec:decision}

In this section, we briefly state the necessary preliminary definitions and notation used in this work: lower previsions, natural extension (see \cite{1991:walley,2008:miranda::survey:lowprevs,2014:miranda:itip:lowerprevision,2014:troffaes:decooman::lower:previsions}), and the $\Gamma$-maximin, $\Gamma$-maximax and interval dominance decision criteria (see \cite{2007:troffaes:decision:intro}). 

\subsection{Lower previsions}

The \emph{possibility space} is denoted by $\Omega$ and represents the set of states of nature.
A \emph{gamble} is a function $f\colon\Omega \to \mathbb{R}$, viewed as an uncertain reward: after the true state of nature $\omega$ is revealed, the gamble $f$ will yield the payoff $f(\omega)$. Let $\mathcal{L}$ be the set of all gambles on $\Omega$.

A \emph{lower prevision} $\underline{P}$ is a real-valued function defined on
$\dom\underline{P} \subseteq\mathcal{L}$, where $\dom\underline{P}$ denotes the domain of $\underline{P}$. Let $f$ be a gamble in $\dom\underline{P}$. We interpret $\underline{P}(f)$ as the subject's supremum buying price for $f$. This means that for all $\alpha<\underline{P}(f)$, the gamble $f-\alpha$ is desirable to the subject.
Lower previsions are suitable to handle the subject's uncertainty when little information is available, because they do not require a full probability specification
\cite{1991:walley,2008:miranda::survey:lowprevs,2014:miranda:itip:lowerprevision,2014:troffaes:decooman::lower:previsions}.

Let us now state a basic consistency requirement for lower previsions. 
\begin{definition} 
A lower prevision $\underline{P}$ \emph{avoids sure loss} if 
for all $ n \in \mathbb{N}$, all  $\lambda_{1}, \dots,\lambda_{n} \geq 0$, and all $f_{1}, \dots,f_{n} \in \dom\underline{P}$, it holds that \cite[p.~42]{2014:troffaes:decooman::lower:previsions}:
\begin{equation}\label{eq:asl}
\max_{\omega\in \Omega} \left( \sum_{i=1}^{n} \lambda_{i}\left[f_{i}(\omega)-\underline{P}(f_{i})\right] \right) \geq 0. 
\end{equation} 
\end{definition}

When \cref{eq:asl} does not hold, then there are $f_{1}, \dots,f_{n} \in \dom\underline{P}$ and $\lambda_{1}, \dots,\lambda_{n} \geq 0$ such that: 
\begin{equation}
\max_{\omega\in \Omega} \left(\sum_{i=1}^{n} \lambda_{i}f_{i}(\omega) \right)  < \sum_{i=1}^{n} \lambda_{i}\underline{P}(f_{i}),
\end{equation}  
and we say that $\underline{P}$ does not avoid sure loss \cite[p.~44]{2014:troffaes:decooman::lower:previsions}.
This implies that the subject agrees to pay more than the maximum possible gain, which is not sensible.
All lower previsions further are assumed to avoid sure loss. Note that we could require lower previsions to be coherent which is a stronger rationality requirement, namely, that the highest price a subject is willing to pay for a gamble cannot be increased by considering combinations of other acceptable gambles. However, coherence is not necessary here as all formulas that we will use rely on natural extension, which only requires avoiding sure loss (see \cite[\S 3.1]{1991:walley} for further details).

Given any lower prevision $\underline{P}$, we 
can also consider its conjugate defined by $\overline{P}(f)\coloneqq -\underline{P}(-f)$ on $\dom \overline{P}\coloneqq\{-f: f\in \dom \underline{P}\}$. We call $\overline{P}$ an \emph{upper prevision}
which can be viewed as the infimum acceptable selling price for $f$ to the subject \cite[p.~41]{2014:troffaes:decooman::lower:previsions}.

A lower prevision $\underline{P}$ can be extended to the set of all gambles $\mathcal{L}$ via the following method:
\begin{definition} \cite[p.~47]{2014:troffaes:decooman::lower:previsions}
Let $\underline{P}$ be a lower prevision. The natural extension of $\underline{P}$ is given, for all $g \in  \mathcal{L}$, by:
\begin{equation}\label{eq:lowerE}
\underline{E}(g)\coloneqq\sup \left\lbrace\alpha \in \mathbb{R}\colon g -\alpha \geq \sum_{i=1}^{n} \lambda_{i}(f_{i}-\underline{P}(f_{i})) , n\in\mathbb{N},\,f_{i} \in  \dom\underline{P},\,\lambda_{i} \geq 0\right\rbrace.
\end{equation} 
\end{definition}
This can be interpreted as follows: based on the prices $\underline{P}(f_i)$ for all $f_i \in \dom \underline{P}$,
$\underline{E}(g)$ is the subject's supremum price that subject is willing to pay for $g\in\mathcal{L}$ \cite[p.~47]{2014:troffaes:decooman::lower:previsions}.
Note that if $\underline{P}$ avoids sure loss, then $\underline{E}$ is finite, and therefore, is a lower prevision \cite[p.~68]{2014:troffaes:decooman::lower:previsions}. 
Moreover, if both $\Omega$ and $\dom \underline{P}$ are finite, then $\underline{E}(g)$ can be computed directly by solving a linear program \cite[p.~331]{2014:troffaes:itip:computation}.
Therefore, throughout this study, both $\Omega$ and $\dom \underline{P}$ are assumed finite.

Similarly, we consider the conjugate of $\underline{E}$ which is denoted by $\overline{E}$ and given by
\begin{align}
  \overline{E}(g)& \coloneqq -\underline{E}(-g)\\ \label{eq:underE}
 & = \inf \left\lbrace\beta \in \mathbb{R}\colon \beta - g \geq \sum_{i=1}^{n} \lambda_{i}(f_{i}-\underline{P}(f_{i})), n \in \mathbb{N}, f_{i} \in \dom\underline{P}, \lambda_{i} \geq 0\right\rbrace.
\end{align}
An upper prevision $\overline{E}$ can be interpreted as the subject's infimum selling price for a gamble $g \in \mathcal{L}$ \cite[p.~47]{2014:troffaes:decooman::lower:previsions}.

\subsection{Decision criteria}

In this section, we review three decision criteria that are used in the context of lower previsions: $\Gamma$-maximin, $\Gamma$-maximax and interval dominance.

\begin{definition}\cite{2007:troffaes:decision:intro}
	Let $\mathcal{K}\subseteq\mathcal{L}$. A gamble is called $\Gamma$-maximin in $\mathcal{K}$ if it maximizes the lower natural extension in $\mathcal{K}$, that is, if it belongs to the set
	\begin{align}\label{eq:maximin}
		\opt_{\underline{E}}(\mathcal{K}) & \coloneqq \{f \in \mathcal{K}\colon (\forall g \in \mathcal{K})(\underline{E}(f) \geq \underline{E}(g) )\} \\
& = \argmax_{f\in \mathcal{K}}\underline{E}(f).
	\end{align}
\end{definition}

In the above, $\argmax_{f\in \mathcal{K}}\underline{E}(f)$ denotes the set of gambles whose elements attain the lower natural extension's largest value. So, basically, a $\Gamma$-maximin gamble is a gamble whose lower prevision (after natural extension) is maximal across all gambles in $\mathcal{K}$. Similarly, we have:

\begin{definition}\cite{2007:troffaes:decision:intro}
  Let $\mathcal{K}\subseteq\mathcal{L}$. A gamble is called $\Gamma$-maximax in $\mathcal{K}$ if it maximizes the upper natural extension in $\mathcal{K}$, that is, if it belongs to the set
\begin{align}\label{eq:maximax}
\opt_{\overline{E}}(\mathcal{K}) & \coloneqq \{f \in \mathcal{K}\colon (\forall g \in \mathcal{K})(\overline{E}(f) \geq \overline{E}(g) )\} \\
& = \argmax_{f\in \mathcal{K}}\overline{E}(f),
\end{align}
\end{definition}
Here, $\argmax_{f\in \mathcal{K}}\overline{E}(f)$ denotes the set of gambles whose elements attain the upper natural extension's largest value.

Next,  we consider \textit{interval dominance} which is based on a strict partial preference order. 
\begin{definition}
	The set of interval dominant gambles of $\mathcal{K}$
	is
	\begin{align}
		\opt_{\sqsupset}(\mathcal{K}) & \coloneqq \{f \in \mathcal{K}\colon (\forall g \in \mathcal{K})(g \not\sqsupset f )\} \\
		&=\{f \in \mathcal{K}: \overline{E}(f) \geq \max_{g \in \mathcal{K}}\underline{E}(g)\}.
	\end{align}
        where we used the following strict partial preference order:
	\begin{equation}
		f \sqsupset g \text{ if } \underline{E}(f) > \overline{E}(g)
	\end{equation}
\end{definition}
Note that \citep{2014:troffaes:itip:decision}:
\begin{equation}
	\opt_{\underline{E}}(\mathcal{K}) \cup \opt_{\overline{E}}(\mathcal{K})  \subseteq \opt_{\sqsupset}(\mathcal{K}).
\end{equation}

\section{Algorithms}\label{sec:algorithm}

In this section, we will discuss algorithms for  $\Gamma$-maximin, $\Gamma$-maximax and interval dominance and we will propose new algorithms based on improvements in  \citep{2019:Nakharutai:Troffase:Caiado:maximal}. 

\subsection{Base Algorithms}

Recall that we only consider the situation where the domain of the lower prevision $\underline{P}$ and the possibility space $\Omega$ are finite.

A first straightforward algorithm is simply to calculate $\underline{E}(f)$ for each $f\in\mathcal{K}$, by solving linear programs. Let $k = |\mathcal{K}|$. To find $\underline{E}(f)$, we can either solve \linprogref{P1} or \linprogref{D1}: 
\begin{align}
\label{P1:1}\tag{P1a}
\linprogref{P1} &&
\min \quad & \sum_{\omega \in \Omega} f(\omega)p(\omega)  \\
\label{P1:2}\tag{P1b}
&& \text{subject to}\quad & \forall g_i\in\dom\underline{P} \colon \sum_{\omega \in \Omega} (g_i(\omega)-\underline{P}(g_i))p(\omega) \geq 0\\
\label{P1:3}\tag{P1c}
&& & \sum_{\omega \in \Omega}p(\omega) = 1\\
\label{P1:4}\tag{P1d}
&&\text{where} \quad  & \forall \omega\colon p(\omega) \geq 0,
\end{align}

\begin{align}
\tag{D1a}
\linprogref{D1}  &&
\max\quad & \alpha \\
\label{thm1:5}\tag{D1b}
&& \text{subject to}\quad & \forall \omega \in\Omega\colon  \sum_{i=1}^{k} (g_i(\omega)-\underline{P}(g_i))\lambda_{i} + \alpha \leq f(\omega) \\
\label{thm1:3}\tag{D1c}
&& \text{where} \quad & \forall i\colon \lambda_{i} \geq 0 \quad (\alpha \text{ free}).
\end{align}
Note that $\underline{E}(f)$ is exactly the optimal value of these linear programs and  \linprogref{D1} corresponds to the unconditional case of the linear program in \cite[p.~331]{2014:troffaes:itip:computation}. Also note that \cref{P1:2} is equivalent to requiring that the expectation $E_{p}(g_i)$ is greater than $\underline{P}(g_i)$, \cref{P1:3,P1:4} guarantee that $p$ is a probability mass function, and \cref{P1:2,P1:3,P1:4} are used to compute the feasible region of \linprogref{P1} which coincides with the credal set of $\underline{P}$, where the credal set of $\underline{P}$ is a closed convex set of probability mass functions over $\Omega$  \citep[p.~135]{1991:walley}. Moreover, \linprogref{D1} is equivalent to \cref{eq:lowerE}.

To find $\argmax_{f\in \mathcal{K}}\underline{E}(f)$, we need to solve $k$ linear programs to get their corresponding optimal values; see \cref{alg:maximin1}.

\begin{algorithm}\caption{$\Gamma$-maximin (original)}\label{alg:maximin1}
	\begin{algorithmic}[1]
		\Require a set of $k$ gambles $\mathcal{K} = \{f_1,\dots, f_k \}$; a small number $\epsilon$; initial states $x_i^P$ and $x_i^D$ for \linprogref{P1} and \linprogref{D1} corresponding to $f_i$ respectively;
		\Ensure a single $\Gamma$-maximin gamble
		\Procedure{$\Gamma$-maximin1}{$\mathcal{K}, \epsilon$} 
		\State $\underline{e} \gets -\infty$
		\For {$i = 1:k$}
		\Repeat 
		\State{$(x_i^P,x_i^D,r^P_i,r^D_i) \gets \phi(f_i,x_i^P,x_i^D)$} \Comment{next iteration to compute $\underline{E}(f_i)$} \\
                \Comment{$r^P_i$ is the primal residual} \\
                \Comment{$r^D_i$ is the dual residual}
		\State{$\ell_i \gets  \underline{e}_*(f_i,x_i^P,x_i^D)$}  \Comment{the current lower bound for $\underline{E}(f_i)$}
		\State{$u_i \gets \underline{e}^*(f_i, x_i^P,x_i^D)$} \Comment{the current upper bound for $\underline{E}(f_i)$}
		\State{$\gamma_i =  u_i-\ell_i$}		
		\Until{$\max\{\gamma_i,r^P_i,r^D_i\}< \epsilon$}
		\If{$\ell_i > \underline{e}$} 
		\State $i^* \gets i $; $\underline{e} \gets \ell_i$ 
		\EndIf
		\EndFor 
		\State \Return $i^*$  \Comment{$i^*$ is $\Gamma$-maximin index}
	\EndProcedure
	\end{algorithmic}
\end{algorithm}
Even though \cref{alg:maximin1} is straightforward, we present the full algorithm here to introduce the reader to the syntax, and also for ease of comparison with algorithms that will be introduced later.

There are several common linear programming methods, such as the simplex method, and interior-point methods such as the affine scaling and the primal-dual method. These methods are iterative, i.e., given a starting point, the methods will find a next feasible solution that comes closer to, and eventually converges to, an optimal solution (provided that an optimal solution exists).

From earlier work \citep{2018:Nakharutai:Troffaes:Caiado} where we performed in-depth analysis of algorithms for working with lower previsions, we know that the primal-dual method, which simultaneously solves both primal and dual problems, performs very well for solving such linear programs. So in the simulation study in \cref{sec:benchmark}, we will solve all linear programs only by the primal-dual method and we do not consider other methods.

What is interesting about the primal-dual method is that it can start from an arbitrary (not necessarily feasible) point. The method will then iterate over points that simultaneously move the solution towards feasibility and optimality. Once feasibility is achieved, all next iterations will retain feasibility, and the solution will eventually converge to an optimal solution \citep[\S 7.3]{1993:Fang:Puthenpura}.

We now explain how $\underline{E}(f_i)$ is calculated in \cref{alg:maximin1} through the primal-dual method.
Let $(x_i^P,x_i^D)$ denote a current pair of (potentially infeasible) points, one for \linprogref{P1}, and one for \linprogref{D1}, corresponding to the linear program for $f_i$. The corresponding \linprogref{P1} value is denoted by $\underline{e}^*(f_i,x_i^P,x_i^D)$, and the \linprogref{D1} value by $\underline{e}_*(f_i,x_i^P,x_i^D)$. Let $\phi$ denote a function that updates this pair of points, representing one iterative step of the algorithm. At each iteration, the method will return an updated point $\phi(f_i,x_i^P,x_i^D)$.
The process will be repeated until $(x_i^P,x_i^D)$ is feasible (i.e., $r^P_i,r^D_i < \epsilon$) and the duality gap $\underline{e}^*(f_i,x_i^P,x_i^D) - \underline{e}_*(f_i,x_i^P,x_i^D)$ is small enough. Once the method terminates, $\underline{E}(f_i)$ can be read from $\underline{e}_*(f_i,x_i^P,x_i^D)$ (or $\underline{e}^*(f_i,x_i^P,x_i^D)$).

Similarly, for $\Gamma$-maximax, a straightforward algorithm is to compute $\overline{E}(f)$ for all $f \in \mathcal{K}$. Note that $\overline{E}(f)$ can be obtained either solving \linprogref{P2} or \linprogref{D2}: 
\begin{align}
\label{P2:1}\tag{P2a}
\linprogref{P2} &&
\min \quad & \beta  \\
\label{P2:2}\tag{P2b}
&& \text{subject to}\quad & \forall \omega \in\Omega\colon \beta - \sum_{i=1}^{k} (g_i(\omega)-\underline{P}(g_i))\lambda_{i} \geq f(\omega)\\
\label{P2:3}\tag{P2c}
&&  \text{where} \quad & \forall i\colon \lambda_{i} \geq 0 \quad (\beta \text{ free}),
\end{align}

\begin{align}
\tag{D2a}
\linprogref{D2}  &&
\max\quad & \sum_{\omega \in \Omega} f(\omega)p(\omega)  \\
\label{D2:2}\tag{D2b}
&& \text{subject to}\quad &  \forall g_i\in\dom\underline{P} \colon \sum_{\omega \in \Omega} (g_i(\omega)-\underline{P}(g_i))p(\omega) \geq 0 \\
\label{D2:3}\tag{D2c}
&& & \sum_{\omega \in \Omega}p(\omega) = 1\\
\label{D2:4}\tag{P2d}
&&\text{where} \quad  & \forall \omega\colon p(\omega) \geq 0.
\end{align}
To find $\argmax_{f\in \mathcal{K}}\overline{E}(f)$, again we have to solve $k = |\mathcal{K}|$ linear programs to get their corresponding optimal values; see \cref{alg:maximax1}. In this case, \cref{D2:2} implies that $E_{p}(g_i) \geq \underline{P}(g_i)$ and \cref{D2:3,D2:4} ensures that $p$ is a probability mass function. In addition, \cref{D2:2,D2:3,D2:4} are used to find the feasible region of \linprogref{D2} which is exactly the credal set of $\underline{P}$ while \linprogref{P2} is equivalent to \cref{eq:underE}. 

Note that we set up the linear programs for $\Gamma$-maximax differently from the linear programs for $\Gamma$-maximin, because we want to match the standard form for the primal-dual method, namely the primal problem is to minimise the objective function and the dual problem is to maximise the objective function. As the pair of problems are dual to each other, we can choose either one of the pair as the primal problem and the other one then becomes its dual problem \citep[p.~174]{2009:Griva:Nash:Sofer}.

\begin{algorithm}
\caption{$\Gamma$-maximax (original)}\label{alg:maximax1}
	\begin{algorithmic}[1]
		\Require a set of $k$ gambles $\mathcal{K} = \{f_1,\dots, f_k \}$; a small number $\epsilon$; initial states $y_i^P$ and $y_i^D$ for \linprogref{P2} and \linprogref{D2} corresponding to $f_i$ respectively;
		\Ensure a single $\Gamma$-maximax gamble
		\Procedure{$\Gamma$-maximax1}{$\mathcal{K}, \epsilon$} 
		\State $\overline{e} \gets -\infty$
		\For {$i = 1:k$}
		\Repeat 
		\State{$(y_i^P,y_i^D,r^P_i,r^D_i) \gets \psi(f_i,y_i^P,y_i^D)$} \Comment{next iteration to compute $\overline{E}(f_i)$} \\
                \Comment{$r^P_i$ is the primal residual} \\
                \Comment{$r^D_i$ is the dual residual}
		\State{$\ell_i \gets  \overline{e}_*(f_i,y_i^P,y_i^D)$}  \Comment{the current lower bound for $\overline{E}(f_i)$}
		\State{$u_i \gets \overline{e}^*(f_i,y_i^P,y_i^D)$} \Comment{the current upper bound for $\overline{E}(f_i)$}
		\State{$\gamma_i =  u_i-\ell_i$}		
		\Until{$\max\{\gamma_i,r^P_i,r^D_i\} < \epsilon$}
		\If{$u_i > \overline{e}$} 
		\State $i^* \gets i $; $\overline{e} \gets u_i$ 
		\EndIf
		\EndFor 
		\State \Return $i^*$  \Comment{$i^*$ is a $\Gamma$-maximax index}
	\EndProcedure
	\end{algorithmic}
\end{algorithm}

The explanation of \cref{alg:maximax1} is similar as for \cref{alg:maximin1}. Here, $(y_i^P,y_i^D)$ denotes a current pair of (potentially infeasible) points for \linprogref{P2} and \linprogref{D2} respectively, corresponding to $f_i$.
The corresponding \linprogref{P2} value is denoted by $\overline{e}^*(f_i,y_i^P,y_i^D)$, and the \linprogref{D2} value by $\overline{e}_*(f_i,y_i^P,y_i^D)$.
At each iteration, the method will return an updated point $\psi(f_i,y_i^P,y_i^D)$.
This process is repeated until $(y_i^P,y_i^D)$ is feasible (i.e., $r^P_i,r^D_i < \epsilon$) and the duality gap $\overline{e}^*(f_i,y_i^P,y_i^D) - \overline{e}_*(f_i,y_i^P,y_i^D)$ is small enough.
When the method terminates, $\overline{E}(f_i)$ can be read from $\overline{e}_*(f_i,y_i^P,y_i^D)$ (or $\overline{e}^*(f_i,y_i^P,y_i^D)$). 

Next, we consider an algorithm for interval dominance. Let $\mathcal{K}$ be a set of $k$ gambles. We are going to find a set of interval dominant gambles in $\mathcal{K}$. To do so, for each $f_i$, we are going to check whether $f_i$ is interval dominant or not.  Note that $f_i$ is interval dominant in the set $\mathcal{K}$ if 
\begin{equation}\label{eq:check_id}
	\overline{E}(f_i) \geq \max_{f_j\in \mathcal{K}}\underline{E}(f_j),
\end{equation}
where $\argmax_{f_j \in \mathcal{K}}\underline{E}(f_j)$ is the $\Gamma$-maximin gamble. 
In other words,  for each $f$,  we compare it against the $\Gamma$-maximin gamble.  We can calculate $\underline{E}(f_i)$ by either solving \linprogref{P1} or \linprogref{D1} and calculate $\overline{E}(f_i)$ by either solving \linprogref{P2} or \linprogref{D2}. Overall,  to find out a set of interval dominant gambles from a set of $k$ gambles,  we have to solve $2k-1$ linear programs (as $f_\ell \coloneqq \argmax_{f_j\in \mathcal{K}}\underline{E}(f_j)$ is immediately interval dominant) \citep[p.~337]{2014:troffaes:itip:computation}.  An algorithm for finding interval dominant gambles is given in \cref{alg:ID1}.

\begin{algorithm}\caption{Interval dominance (original)}\label{alg:ID1}
	\begin{algorithmic}[1]
		\Require a set of $k$ gambles $\mathcal{K} = \{f_1,\dots, f_k \}$; a small number $\epsilon$; initial states $x_i^P$ and $x_i^D$  for \linprogref{P1} and \linprogref{D1} corresponding to $f_i$ respectively; initial states  $y_i^P$  and $y_i^D$ for \linprogref{P2} and \linprogref{D2} corresponding to $f_i$ respectively;
		\Ensure a single $\Gamma$-maximin gamble; an index set of $\opt_{\sqsupset}(\mathcal{K})$
		\Procedure{IntervalDominant1}{$\mathcal{K}, \epsilon$} 
		\State{\textbf{Stage I:} Apply \textsc{$\Gamma$-maximin1} to obtain $i^*$ and $ \underline{e} $} \Comment{$i^*$ is an index of the $\Gamma$-maximin and $ \underline{e} =\max_{f_j\in \mathcal{K}}\underline{E}(f_j)$}
		
		~
		
		\State{\noindent \textbf{Stage II: Interval dominance}}
		\State $I\gets \{i^*\}$  
		\For{$i \in \{1,2,\dots,k\} \setminus\{i^*\} $}
		\Repeat 
		\State{$(y_i^P,y_i^D,r^P_i,r^D_i)\gets \psi(f_i,y_i^P,y_i^D)$} \Comment{next iteration to compute $\overline{E}(f_i)$}\\
		\Comment{$r^P_i$ is the primal residual} \\
		\Comment{$r^D_i$ is the dual residual}
		\State{$\ell_i \gets  \overline{e}_*(f_i,y_i^P,y_i^D)$}  \Comment{a lower bound for $\overline{E}(f_i)$}
		\State{$u_i \gets \overline{e}^*(f_i,y_i^P,y_i^D)$} \Comment{an upper bound for $\overline{E}(f_i)$}
		\State{$\gamma_i =  u_i-\ell_i$}
		
		\Until{$\max\{\gamma_i,r^P_i,r^D_i\}< \epsilon$}
		\If{$u_i \geq \underline{e} $} 
		\State  $I\gets I \cup \{i\}$ \Comment{$f_i$ is interval dominant}
		\EndIf
		\EndFor
		\State \Return $(i^*, I)$  \Comment{$I$ is an index set of $\opt_{\sqsupset}(\mathcal{K})$ } 
		\EndProcedure
	\end{algorithmic}
\end{algorithm}

\subsection{Improvements}\label{subsec:improve_LP}

We now discuss how to improve \cref{alg:maximin1,alg:maximax1,alg:ID1}, by exploiting several features of the primal-dual method itself.

\begin{figure}
	\centering
	\begin{tikzpicture}
	\draw (4,0.7) node[align=center] {$\underline{E}(f_i)$}; 
	\draw [-] (4,0.3)--(4,-0.01);
	\draw (0,0)--(12,0);
	\draw [thick,red,->] (0,0)--(3.5,0);
	\draw (2,-0.5) node[anchor=center] {$\underline{e}_*(f_i,x_i^P,x_i^D)$}; 
	\draw [thick,blue,<-] (4.5,0)--(8,0);
	\draw (6,-0.5) node[anchor=center]{$ \underline{e}^*(f_i,x_i^P,x_i^D)$}; 
	\draw (9,-1.5) node[align=center] {Stop if $ \underline{e}^*(f_i,x_i^P,x_i^D) < \underline{e}$};
	\draw [thick,->] (8.5,-1)--(8.5,-0.5)--(7.8,-0.5); 
	\draw [dotted](9.7,0.3)--(9.5,0.3)--(9.5,-0.3)--(9.7,-0.3);
	\draw [dotted](10.3,0.3)--(10.5,0.3)--(10.5,-0.3)--(10.3,-0.3);
	\draw [-](10,0.3)--(10,-0.01);
	\draw (10,0.7) node[align=center] {$\underline{E}(f_{i^*})$}; 
	\draw (9.5,-0.7) node[anchor=center]{$ \underline{e}$}; 
	\end{tikzpicture}
	\caption{Early stopping criterion for $\Gamma$-maximin}\label{fig:gamma-min}
\end{figure}
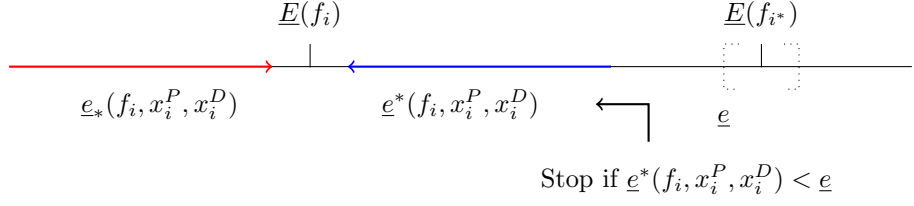

Recall that we only need to verify whether $\underline{E}(f_i) > \underline{e}$ or not.
Note that for any $(x_i^P,x_i^D)$,
\begin{equation}
\underline{e}_*(f_i,x_i^P,x_i^D) \leq \underline{E}(f_i) \leq \underline{e}^*(f_i,x_i^P,x_i^D).
\end{equation} 
Therefore, provided that the current solution is feasible, we can terminate as soon as we find upper bound $\underline{e}^*(f_i,x_i^P,x_i^D)$ that is less than $\underline{e}$ since we immediately know that the optimal value of \linprogref{P1} (or $\underline{E}(f_i)$) is also less than $\underline{e}$ (see \cref{fig:gamma-min}). Hence, it is not necessary to find the optimal value. This gives us an early stopping criterion for $\Gamma$-maximin. 

For this to work best, ideally, we start the primal-dual method with feasible points. To obtain such feasible points, we can reuse the fast techniques proposed in \citep{2017:Nakharutai:Troffaes:Caiado,2018:Nakharutai:Troffaes:Caiado,2019:Nakharutai:Troffase:Caiado:maximal}. 
Specifically, to find a starting feasible solution $p(\omega)$ for \linprogref{P1} and \linprogref{D2}, we could employ the first phase of the two-phase method, similar to earlier work \citep[\S 4.2]{2018:Nakharutai:Troffaes:Caiado}.
Note that the constraints of each \linprogref{P1} and \linprogref{D2} do not depend on the gamble $f\in\mathcal{K}$. If we find their feasible starting points, then these feasible solutions can be reused for \linprogref{P1} and \linprogref{D2} across all $f\in\mathcal{K}$. In other words, for any given lower prevision, we only have to find their corresponding starting points once at the beginning. 
On the other hand, to obtain initial feasible solutions for \linprogref{D1} and \linprogref{P2}, there is no need to solve any linear program as we can use a result from \citep[Theorem~7]{2018:Nakharutai:Troffaes:Caiado} to directly calculate feasible starting points.

Also note that if we can identify a $\Gamma$-maximin gamble early, then fewer iterations are required. Fortunately, if we sort all gambles in $\mathcal{K}$ as $f_1$, \dots, $f_k$ such that for some probability mass function $p$ (which we have to find anyway as we need an initial feasible solution), for all $i < j$:
\begin{equation}\label{eq:E_p(i-f_j)}
E_p(f_i)\geq E_p(f_j) ,
\end{equation}
then it is more likely that $f_i$ is $\Gamma$-maximin for a lower index $i$, and we have to do fewer computations as a result.
Note that this technique helps in most cases, but not always. In \citep[Alg 2]{2019:Nakharutai:Troffase:Caiado:maximal}, we showed that this technique did help in the most cases even though it is just heuristic and may not always help. An algorithm that implements this ordering trick, along with the early stopping criterion, is presented in \cref{alg:maximin2}. 
Note that we do not need to check feasibility in this algorithm as we start with feasible points \citep[\S 10.7]{2009:Griva:Nash:Sofer}.

\begin{algorithm}
	\caption{$\Gamma$-maximin}\label{alg:maximin2}
	\begin{algorithmic}[1]
	\Require  a set of $k$ gambles $\mathcal{K} = \{f_1,\dots, f_k \}$ such that for some probability mass function $p$, $E_p(f_1)\ge E_p(f_2)\ge \dots \ge E_p(f_k)$; a small number $\epsilon$; initial feasible states $x_i^P$ and $x_i^D$ for \linprogref{P1} and \linprogref{D1} corresponding to $f_i$ respectively;
	\Ensure a single $\Gamma$-maximin gamble
	\Procedure{$\Gamma$-maximin2}{$\mathcal{K}, \epsilon$}
	\State $\underline{e} \gets -\infty$
	\For {$i = 1:k$}
	\Repeat 
	\State{$(x_i^P,x_i^D)  \gets \phi(f_i,x_i^P,x_i^D)$} \Comment{next iteration to compute $\underline{E}(f_i)$}
	\State{$\ell_i \gets  \underline{e}_*(f_i,x_i^P,x_i^D)$}  \Comment{the current lower bound for $\underline{E}(f_i)$}
	\State{$u_i \gets \underline{e}^*(f_i,x_i^P,x_i^D)$} \Comment{the current upper bound for $\underline{E}(f_i)$}
	\State{$\gamma_i =  u_i-\ell_i$}
	\Until{$u_i < \underline{e}$ \textbf{or} $\gamma_i < \epsilon$}  \Comment{an early stopping criterion}
	\If{ $\ell_i > \underline{e}$}
	\State {$\underline{e} \gets \ell_i$; $i^* \gets i$} 
	\EndIf
	\EndFor 
	\State \Return $i^*$   \Comment{$i^*$ is a $\Gamma$-maximin index}
	\EndProcedure
	\end{algorithmic}
\end{algorithm}

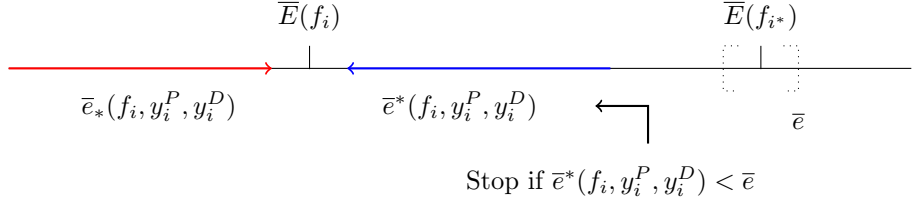
\begin{figure}
	\centering
	\begin{tikzpicture}
	\draw (4,0.7) node[align=center] {$\overline{E}(f_i)$}; 
	\draw [-] (4,0.3)--(4,-0.01);
	\draw (0,0)--(12,0);
	\draw [thick,red,->] (0,0)--(3.5,0);
	\draw (2,-0.5) node[anchor=center] {$\overline{e}_*(f_i,y_i^P,y_i^D)$}; 
	\draw [thick,blue,<-] (4.5,0)--(8,0);
	\draw (6,-0.5) node[anchor=center]{$ \overline{e}^*(f_i,y_i^P,y_i^D)$}; 
	\draw (8,-1.5) node[align=center] {Stop if $ \overline{e}^*(f_i,y_i^P,y_i^D) < \overline{e}$};
	\draw [thick,->] (8.5,-1)--(8.5,-0.5)--(7.8,-0.5); 
	\draw [dotted](9.7,0.3)--(9.5,0.3)--(9.5,-0.3)--(9.7,-0.3);
	\draw [dotted](10.3,0.3)--(10.5,0.3)--(10.5,-0.3)--(10.3,-0.3);
	\draw [-](10,0.3)--(10,-0.01);
	\draw (10,0.7) node[align=center] {$\overline{E}(f_{i^*})$}; 
	\draw (10.5,-0.7) node[anchor=center]{$ \overline{e}$}; 
	\end{tikzpicture}
	\caption{Early stopping criterion for $\Gamma$-maximax}\label{fig:gamma-max}
\end{figure}

These improvements can be applied to $\Gamma$-maximax as well. Specifically, we exploit the fact that we only have to verify whether $\overline{E}(f_i) > \overline{e}$ or not. Since for any $(y_i^P,y_i^D)$,
\begin{equation}
\overline{e}_*(f_i,y_i^P,y_i^D) \leq \overline{E}(f_i) \leq \overline{e}^*(f_i,y_i^P,y_i^D),
\end{equation} 
we can terminate as soon as we find $\overline{e}^*(f_i,y_i^P,y_i^D)$ that is less than $\overline{e}$ because then we know that the optimal value of \linprogref{P2}, i.e. $\overline{E}(f_i)$, is less than $\overline{e}$ (see \cref{fig:gamma-max}).

Again, if we sort the gambles in $\mathcal{K} = \{f_1,\dots, f_k\}$ such that \cref{eq:E_p(i-f_j)} is satisfied, then it is more likely that $f_i$ is $\Gamma$-maximax for a lower index $i$, and we can save computations.
See \cref{alg:maximax2} for a description of an improved algorithm for $\Gamma$-maximax. 

\begin{algorithm}
	\caption{$\Gamma$-maximax}\label{alg:maximax2}
	\begin{algorithmic}[1]
	\Require  a set of $k$ gambles $\mathcal{K} = \{f_1,\dots, f_k \}$ such that for some probability mass function $p$, $E_p(f_1)\ge E_p(f_2)\ge \dots \ge E_p(f_k)$; a small number $\epsilon$; initial feasible states $y_i^P$ and $y_i^D$ for \linprogref{P2} and \linprogref{D2} corresponding to $f_i$ respectively;
	\Ensure a single $\Gamma$-maximax gamble
	\Procedure{$\Gamma$-maximax2}{$\mathcal{K}, \epsilon$}
	\State $\overline{e} \gets -\infty$
	\For {$i = 1:k$}
	\Repeat 
	\State{$(y_i^P,y_i^D) \gets \psi(f_i,y_i^P,y_i^D)$} \Comment{next iteration to compute $\overline{E}(f_i)$}
	\State{$\ell_i \gets  \overline{e}_*(f_i,y_i^P,y_i^D)$}  \Comment{the current lower bound for $\overline{E}(f_i)$}
	\State{$u_i \gets \overline{e}^*(f_i,y_i^P,y_i^D)$} \Comment{the current upper bound for $\overline{E}(f_i)$}
	\State{$\gamma_i =  u_i-\ell_i$}
	\Until{$u_i < \overline{e}$ \textbf{or} $\gamma_i < \epsilon$}  \Comment{an early stopping criterion}
	\If{$u_i > \overline{e}$}
	\State{$\overline{e} \gets u_i$; $i^* \gets i$} 
	\EndIf
	\EndFor 
	\State \Return $i^*$   \Comment{$i^*$ is a $\Gamma$-maximax index}
	\EndProcedure
	\end{algorithmic}
\end{algorithm}

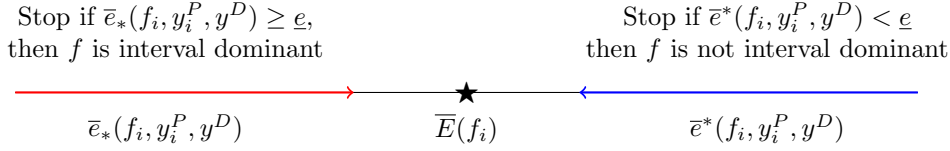
\begin{figure}
	\centering
	\begin{tikzpicture}
		\draw (2,0.8) node[align=center] {Stop if $\overline{e}_*(f_i,y_i^P,y^D) \geq \underline{e}$,\\ then $f$ is interval dominant};
		\draw (10,0.8) node[align=center] {Stop if $\overline{e}^*(f_i,y_i^P,y^D) <  \underline{e}$\\ then $f$ is not interval dominant}; 
		\draw (6,-0.8) node[anchor=south]{$\overline{E}(f_i)$};  
		\draw (6,-0.01) node{$\bigstar$};
		\draw (0,0)--(12,0);
		\draw [thick,red,->] (0,0)--(4.5,0);
		\draw (2,-0.8) node[anchor=south]{$\overline{e}_*(f_i,y_i^P,y^D)$}; 
		\draw [thick,blue,<-] (7.5,0)--(12,0);
		\draw (10,-0.8) node[anchor=south]{ $\overline{e}^*(f_i,y_i^P,y^D)$}; 
	\end{tikzpicture}
	\caption{Early stopping criterion for interval dominance}\label{fig:early-stop-id}
\end{figure}

For interval dominance, even though we have to solve $2k-1$ linear programs, we do not need to find all these optimal values since we can apply these improvements for $\Gamma$-maximin and $\Gamma$-maximax. Specifically, in the stage I of \cref{alg:ID1}, we can apply \cref{alg:maximin2} to quickly obtain the $\Gamma$-maximin gamble, say $f_{i^*}$, and $ \underline{e} =\max_{f_j\in \mathcal{K}}\underline{E}(f_j)$. 
In the stage II of of \cref{alg:ID1}, we can speed up the process of evaluating $\overline{E}(f)$ through solving both \linprogref{P2} and \linprogref{D2} similarly to what we have done in  \cref{alg:maximax2}, exploiting the fact that we only need to verify whether or not $\overline{E}(f_i) \geq \underline{e} $ (in this stage, $\underline{e}$ is $\max_{f_j\in \mathcal{K}}\underline{E}(f_j)$). 
Note that for any $(y_i^P,y_i^D)$,
\begin{equation}
	\overline{e}_*(f_i,y_i^P,y_i^D) \leq \overline{E}(f_i) \leq \overline{e}^*(f_i,y_i^P,y_i^D).
\end{equation} 
So, we can stop as soon as we find $\overline{e}^*(f_i,y_i^P,y_i^D)$ that is less than $\underline{e}$ because then we know that the optimal value of \linprogref{P2}, i.e. $\overline{E}(f_i)$, is immediately less than $\underline{e}$. In this case, $f$ is not interval dominant in $\mathcal{K}$. Meanwhile, we can stop as soon as we find $\overline{e}_*(f_i,y_i^P,y_i^D)$ that is larger or equal to $\underline{e}$ since we know that the optimal value of \linprogref{D2}, i.e. $\overline{E}(f_i)$, will be larger or equal to  $\underline{e}$ as well. In this case, $f$ is interval dominant in $\mathcal{K}$. \Cref{fig:early-stop-id} summarises these early stopping criteria.
An algorithm for interval dominance implemented the early stopping criteria is presented in \cref{alg:ID2}.

\begin{algorithm}
	\caption{Interval dominance}\label{alg:ID2}
	\begin{algorithmic}[1]
		\Require  a set of $k$ gambles $\mathcal{K} = \{f_1,\dots, f_k \}$ such that for some probability mass function $p$, $E_p(f_1)\ge E_p(f_2)\ge \dots \ge E_p(f_k)$; a small number $\epsilon$; initial feasible states $x_i^P$ and $x_i^D$ for \linprogref{P1} and \linprogref{D1} corresponding to $f_i$ respectively; initial feasible states $y_i^P$ and $y_i^D$ for \linprogref{P2} and \linprogref{D2} corresponding to $f_i$ respectively;
		\Ensure a single $\Gamma$-maximin gamble; an index set of $\opt_{\sqsupset}(\mathcal{K})$
		\Procedure{IntervalDominant2}{$\mathcal{K}, \epsilon$}
		\State{\textbf{Stage I:} Apply \textsc{$\Gamma$-maximin2} to obtain $i^*$ and $\underline{e}$} \Comment{$i^*$ is an index of the $\Gamma$-maximin gamble and $\underline{e}=\max_{f_j\in \mathcal{K}}\underline{E}(f_j)$}
		
		~
		
		\State{\noindent \textbf{Stage II: Interval dominance}}
		\State $I \gets \{i^*\}$ \Comment{$i^*$ is an index of the $\Gamma$-maximin gamble}
		\For{$i \in \{1,2,\dots,k\} \setminus\{i^*\} $}
		\Repeat 
		\State{$(y_i^P,y_i^D) \gets \psi(f_i,y_i^P,y_i^D)$} \Comment{next iteration to compute $\overline{E}(f_i)$} 
		\State{$\ell_i \gets  \overline{e}_*(f_i,y_i^P,y_i^D)$}  \Comment{a lower bound for $\overline{E}(f_i)$}
		\State{$u_i \gets \overline{e}^*(f_i,y_i^P,y_i^D)$} \Comment{an upper bound for $\overline{E}(f_i)$}
		\State{$\gamma_i = u_i-\ell_i$}
		\State{\textbf{if} $u_i < \underline{e}$ \textbf{then break}} \Comment{$f_i$ is not interval dominant}
		\State{\textbf{else if} $\ell_i \geq \underline{e}$  \textbf{or} $\gamma_i < \epsilon$  \textbf{then} }\Comment{$f_i$ is interval dominant}
		\State{\ \ \ \ \ \ $I \gets I 	\cup \{i\}$}; \textbf{break} 
		\State{\textbf{end if}}
		\Until{\False}
		\EndFor
		\State \Return $(i^*,I)$   \Comment{$I$ is an index set of $\opt_{\sqsupset}(\mathcal{K})$} 
		\EndProcedure
	\end{algorithmic}
\end{algorithm}

\subsection{Further Improvements}
\label{sec:furtherimprovements}

Note that in order to find a $\Gamma$-maximin gamble
by \cref{alg:maximin2}, in the first stage, we have to iterate until we obtain the optimal value for $f_1$. Only then can we start applying the early stopping criteria.
In addition, the output of \cref{alg:maximin1,alg:maximin2} is just a single gamble even though there can be more than one gamble that is $\Gamma$-maximin. We will consider another algorithm that can find $\Gamma$-maximin gambles without having to obtain the optimal value even of the first gamble, and that can find all $\Gamma$-maximin gambles.

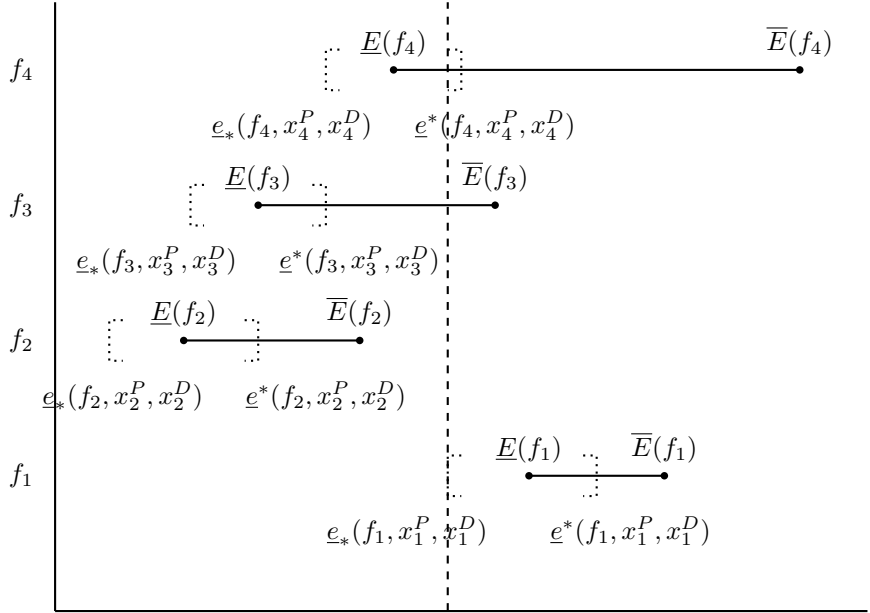
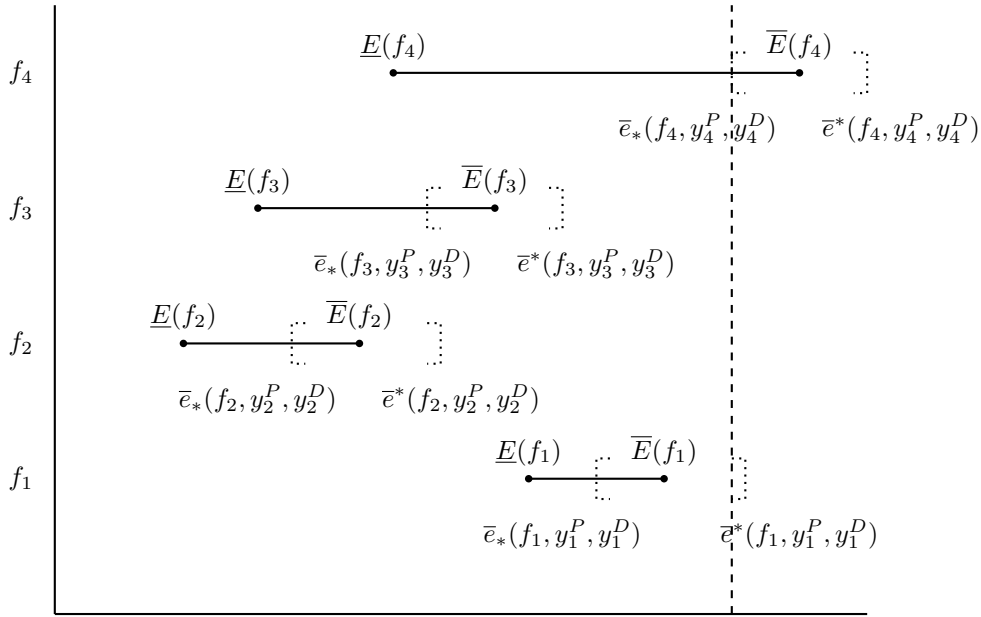
\begin{figure}
	\centering
\subfigure[Elimination of non $\Gamma$-maximin gambles.]{	
	\begin{tikzpicture}[ scale=.9]  
	\draw[thick] (0,0)--(12,0);
	\draw[thick] (0,0)--(0,9);
	\draw (-0.5,2) node[align=left] {$f_1$}; 
	\draw (-0.5,4) node[align=left] {$f_2$}; 
	\draw (-0.5,6) node[align=left] {$f_3$}; 
	\draw (-0.5,8) node[align=left] {$f_4$}; 
	
	\draw[thick] (7,2)--(9,2);
	\draw (7,2) node[circle,fill,inner sep=1pt,label=above:$\underline{E}(f_1)$]{}; 
	\draw (9,2) node[circle,fill,inner sep=1pt,label=above:$\overline{E}(f_1)$]{}; 
	\draw [thick,dotted](6,2.3)--(5.8,2.3)--(5.8,1.7)--(6,1.7);
	\draw (5.2,1.7) node[below=0.15cm]{$\underline{e}_*(f_1,x_1^P,x_1^D)$}; 
	\draw [thick,dotted](7.8,2.3)--(8,2.3)--(8,1.7)--(7.8,1.7);
	\draw (8.5,1.7) node[below=0.15cm]{$\underline{e}^*(f_1,x_1^P,x_1^D)$}; 
	\draw [thick, dashed] (5.8,0)--(5.8,9);
	
	\draw[thick] (1.9,4)--(4.5,4);
	\draw (1.9,4) node[circle,fill,inner sep=1pt,label=above:$\underline{E}(f_2)$]{};
	\draw [thick,dotted](1.0,4.3)--(0.8,4.3)--(0.8,3.7)--(1.0,3.7);	
	\draw (1,3.7) node[below=0.15cm]{$\underline{e}_*(f_2,x_2^P,x_2^D)$}; 
	\draw [thick,dotted](2.8,4.3)--(3,4.3)--(3,3.7)--(2.8,3.7);	
	\draw (4,3.7) node[below=0.15cm]{$\underline{e}^*(f_2,x_2^P,x_2^D)$}; 
	\draw (4.5,4) node[circle,fill,inner sep=1pt,label=above:$\overline{E}(f_2)$]{};

	\draw[thick] (3,6)--(6.5,6);
	\draw (3,6) node[circle,fill,inner sep=1pt,label=above:$\underline{E}(f_3)$]{};
	\draw [thick,dotted](2.2,6.3)--(2,6.3)--(2,5.7)--(2.2,5.7);
	\draw (1.5,5.7) node[below=0.15cm]{$\underline{e}_*(f_3,x_3^P,x_3^D)$}; 
	\draw [thick,dotted](3.8,6.3)--(4,6.3)--(4,5.7)--(3.8,5.7);
	\draw (4.5,5.7) node[below=0.15cm]{$\underline{e}^*(f_3,x_3^P,x_3^D)$}; 	
	\draw (6.5,6) node[circle,fill,inner sep=1pt,label=above:$\overline{E}(f_3)$]{}; 
	
	\draw[thick] (5,8)--(11,8);
	\draw (5,8) node[circle,fill,inner sep=1pt,label=above:$\underline{E}(f_4)$]{};
	\draw [thick,dotted](4.2,8.3)--(4,8.3)--(4,7.7)--(4.2,7.7);	
	\draw (3.5,7.7) node[below=0.15cm]{$\underline{e}_*(f_4,x_4^P,x_4^D)$}; 	
	\draw [thick,dotted](5.8,8.3)--(6,8.3)--(6,7.7)--(5.8,7.7);	
	\draw (6.5,7.7) node[below=0.15cm]{$\underline{e}^*(f_4,x_4^P,x_4^D)$}; 
	\draw (11,8) node[circle,fill,inner sep=1pt,label=above:$\overline{E}(f_4)$]{}; 	

	\end{tikzpicture}
} 
   
\subfigure[Elimination of non $\Gamma$-maximax gambles.]{
	\begin{tikzpicture}[ scale=.9]  
\draw[thick] (0,0)--(12,0);
\draw[thick] (0,0)--(0,9);
\draw (-0.5,2) node[align=left] {$f_1$}; 
\draw (-0.5,4) node[align=left] {$f_2$}; 
\draw (-0.5,6) node[align=left] {$f_3$}; 
\draw (-0.5,8) node[align=left] {$f_4$}; 

\draw[thick] (7,2)--(9,2);
\draw (7,2) node[circle,fill,inner sep=1pt,label=above:$\underline{E}(f_1)$]{}; 
\draw (9,2) node[circle,fill,inner sep=1pt,label=above:$\overline{E}(f_1)$]{}; 
\draw [thick,dotted](8.2,2.3)--(8,2.3)--(8,1.7)--(8.2,1.7);
\draw (7.5,1.7) node[below=0.15cm]{$\overline{e}_*(f_1,y_1^P,y_1^D)$}; 
\draw [thick,dotted](10,2.3)--(10.2,2.3)--(10.2,1.7)--(10,1.7);	
\draw (11.0,1.7) node[below=0.15cm]{$\overline{e}^*(f_1,y_1^P,y_1^D)$}; 
\draw [thick, dashed] (10,0)--(10,9);	

\draw[thick] (1.9,4)--(4.5,4);
\draw (1.9,4) node[circle,fill,inner sep=1pt,label=above:$\underline{E}(f_2)$]{};
\draw (4.5,4) node[circle,fill,inner sep=1pt,label=above:$\overline{E}(f_2)$]{}; 
\draw [thick,dotted](3.7,4.3)--(3.5,4.3)--(3.5,3.7)--(3.7,3.7);	
\draw (3,3.7) node[below=0.15cm]{$\overline{e}_*(f_2,y_2^P,y_2^D)$}; 
\draw [thick,dotted](5.5,4.3)--(5.7,4.3)--(5.7,3.7)--(5.5,3.7);	
\draw (6,3.7) node[below=0.15cm]{$\overline{e}^*(f_2,y_2^P,y_2^D)$};

\draw[thick] (3,6)--(6.5,6);
\draw (3,6) node[circle,fill,inner sep=1pt,label=above:$\underline{E}(f_3)$]{};
\draw (6.5,6) node[circle,fill,inner sep=1pt,label=above:$\overline{E}(f_3)$]{}; 
\draw [thick,dotted](5.7,6.3)--(5.5,6.3)--(5.5,5.7)--(5.7,5.7);
\draw (5,5.7) node[below=0.15cm]{$\overline{e}_*(f_3,y_3^P,y_3^D)$}; 
\draw [thick,dotted](7.3,6.3)--(7.5,6.3)--(7.5,5.7)--(7.3,5.7);
\draw (8,5.7) node[below=0.15cm]{$\overline{e}^*(f_3,y_3^P,y_3^D)$}; 

\draw[thick] (5,8)--(11,8);
\draw (5,8) node[circle,fill,inner sep=1pt,label=above:$\underline{E}(f_4)$]{};
\draw (11,8) node[circle,fill,inner sep=1pt,label=above:$\overline{E}(f_4)$]{}; 
\draw [thick,dotted](10.2,8.3)--(10,8.3)--(10,7.7)--(10.2,7.7);	
\draw (9.5,7.7) node[below=0.15cm]{$\overline{e}_*(f_4,y_4^P,y_4^D)$};
\draw [thick,dotted](11.8,8.3)--(12,8.3)--(12,7.7)--(11.8,7.7);	
\draw (12.5,7.7) node[below=0.15cm]{$\overline{e}^*(f_4,y_4^P,y_4^D)$};  	

\end{tikzpicture}	
}
	\caption{Elimination of many gambles.}\label{fig:eli:gamma-min}	
\end{figure}

As an example, consider a set of gambles $\mathcal{K} = \{f_1,f_2,f_3,f_4\}$. Suppose that we would like to find $\Gamma$-maximin gambles in $\mathcal{K}$. To do so, we might start by calculating $\underline{e}_*(f_i,x_i^P,x_i^D)$ and $\underline{e}^*(f_i,x_i^P,x_i^D)$ for each $f_i$ in $\mathcal{K}$ where $(x_i^P,x_i^D)$ is initial states for the linear programs associated with $f_i$. Assume that the result is illustrated as in \cref{fig:eli:gamma-min} (a). Note that $\underline{e}_*(f_1,x_1^P,x_1^D) = \max_{i}\underline{e}_*(f_i,x_i^P,x_i^D)$ and both $\underline{e}^*(f_2,x_2^P,x_2^D)$ and $\underline{e}^*(f_3,x_3^P,x_3^D)$ are smaller than $\underline{e}_*(f_1,x_3^P,x_3^D)$.  At this point, we already know that $f_2$ and $f_3$ cannot be $\Gamma$-maximin since $\underline{E}(f_2)$ and $\underline{E}(f_3)$ must be less than $\underline{e}_*(f_1,x_1^P,x_1^D)$. From this, we can see that we can quickly eliminate gambles that are not $\Gamma$-maximin, before obtaining the optimal value of the first linear program.

We can translate this argument into an algorithm that can sequentially narrow the set of potentially $\Gamma$-maximin gambles. Specifically, for each gamble $f_i$ in $\mathcal{K}$, we evaluate $\underline{e}_*(f_i,x_i^P,x_i^D)$ and $\underline{e}^*(f_i,x_i^P,x_i^D)$ which are a lower bound and an upper bound for $\underline{E}(f_i)$. Next, we calculate $M_*\coloneqq\max_{i}\underline{e}_*(f_i,x_i^P,x_i^D)$ and $M^*\coloneqq\max_{i}\underline{e}^*(f_i,x_i^P,x_i^D)$.
Then, any gamble $f_i$ for which $\underline{e}^*(f_i,x_i^P,x_i^D) < M_*$ can not be $\Gamma$-maximin and therefore can be eliminated. After that, we update the states $(x_i^P,x_i^D)$ for all gambles $f_i$ that are still potentially $\Gamma$-maximin. We repeat this process until either only one gamble is left in a set of potentially $\Gamma$-maximin gambles, or until the difference between $M_*$ and $M^*$ is less than a given tolerance, in which case we have found multiple $\Gamma$-maximin gambles. An algorithm that implements these arguments is presented in \cref{alg:maximin3}. 

\begin{algorithm}
	\caption{$\Gamma$-maximin}\label{alg:maximin3}
	\begin{algorithmic}[1]
		\Require  a set of $k$ gambles $\mathcal{K} = \{f_1,\dots, f_k \}$; a small number $\epsilon$; initial feasible states $x_i^P$ and $x_i^D$ for \linprogref{P1} and \linprogref{D1} corresponding to $f_i$ respectively;
		\Ensure a set of $\Gamma$-maximin gambles		
		\Procedure{$\Gamma$-maximin3}{$\mathcal{K},\epsilon$}
		\State $R \gets \{1,2,\dots, k\}$ \Comment{ an index set of potentially $\Gamma$-maximin}
		\Repeat 
		\State{$\forall i\in R: (x_i^P,x_i^D) \gets \phi(f_i,x_i^P,x_i^D)$} \Comment{the next iteration to compute $\underline{E}(f_i)$}
		\State{$\forall i\in R: \ell_i \gets  \underline{e}_*(f_i,x_i^P,x_i^D)$}  \Comment{the current lower bound for $\underline{E}(f_i)$}
		\State{$\forall i\in R: u_i \gets \underline{e}^*(f_i,x_i^P,x_i^D)$} \Comment{the current upper bound for $\underline{E}(f_i)$}
		\State{$M_{*} = \max_{i\in R}\ell_i$} \Comment{a lower bound for $\Gamma$-maximin value}
		\State{$M^{*} = \max_{i\in R}u_i$} \Comment{an upper bound for $\Gamma$-maximin value}
		\State{$R \gets \{i\in R\colon (u_i \geq M_{*})\}$}
		\Until{$(|R| = 1)$ \textbf{ or } $(M^{*} - M_{*} < \epsilon)$} 
	\State \Return $R$  
	\EndProcedure
	\end{algorithmic}
\end{algorithm}
A key difference between \cref{alg:maximin3} and the previous two algorithms for $\Gamma$-maximin (\cref{alg:maximin1,alg:maximin2}) lies in the output of the algorithms.
When there is more than one $\Gamma$-maximin gamble, \cref{alg:maximin3} will return all $\Gamma$-maximin gambles, while the other two algorithms will return only one of the $\Gamma$-maximin gambles. Also note that \cref{alg:maximin1,alg:maximin3} can be parallelized but not \cref{alg:maximin2}. In particular, \cref{alg:maximin3} is basically a parallel version of \cref{alg:maximin1} with an extra step to avoid unnecessary calculations.

The ideas in \cref{alg:maximin3} can also be applied to eliminate non $\Gamma$-maximax gambles; see \cref{fig:eli:gamma-min} (b) and \cref{alg:maximax3}.

\begin{algorithm}[]
	\caption{$\Gamma$-maximax}\label{alg:maximax3}
	\begin{algorithmic}[1]
		\Require  a set of $k$ gambles $\mathcal{K} = \{f_1,\dots, f_k \}$; a small number $\epsilon$; initial feasible state $y_i^P$ and $y_i^D$ for \linprogref{P2} and \linprogref{D2} corresponding to $f_i$;
		\Ensure a set of $\Gamma$-maximax gambles
		\Procedure{$\Gamma$-maximax3}{$\mathcal{K},\epsilon$}
		\State $R \gets \{1,2,\dots, k\}$ \Comment{ an index set of potentially $\Gamma$-maximax}
			\Repeat 
			\State{$\forall i\in R: (y_i^P,y_i^D)\gets \psi(f_i,y_i^P,y_i^D)$} \Comment{the next iteration to compute $\overline{E}(f_i)$}
		\State{$\forall i\in R: \ell_i \gets  \overline{e}_*(f_i,y_i^P,y_i^D)$}  \Comment{the current lower bound for $\overline{E}(f_i)$}
		\State{$\forall i\in R: u_i \gets \overline{e}^*(f_i,y_i^P,y_i^D)$} \Comment{the current upper bound for $\overline{E}(f_i)$}
		\State{$M_{*} = \max_{i\in R}\ell_i$} \Comment{a lower bound for $\Gamma$-maximax value}
		\State{$M^{*} = \max_{i\in R}u_i$} \Comment{an upper bound for $\Gamma$-maximax value}
		\State{$R \gets \{i\in R\colon (u_i \geq M_{*})\}$}
		\Until{$(|R| = 1)$ \textbf{ or } $(M^{*} - M_{*} < \epsilon)$} 
	\State \Return $R$  
	\EndProcedure
	\end{algorithmic}
\end{algorithm}

Similarly, the ideas in \cref{alg:maximin3,alg:maximax3} can be applied to interval dominance by simultaneously computing both upper and lower bounds for the lower and upper natural extensions.
We have found however that, perhaps counterintuitively at first, this
does not necessarily lead to improved algorithms for interval dominance.
Understanding why this happens requires a somewhat detailed analysis,
to which we turn next.

First, to bound the $\Gamma$-maximin value, for each $f_i$, we compute
$\underline{e}_*(f_i,x_i^P,x_i^D)$ and
$\underline{e}^*(f_i,x_i^P,x_i^D)$, which bound $\underline{E}(f_i)$.
We then compute $M_*\coloneqq\max_{i}\underline{e}_*(f_i,x_i^P,x_i^D)$
and $M^*\coloneqq\max_{i}\underline{e}^*(f_i,x_i^P,x_i^D)$ to bound
the $\Gamma$-maximin value $\max_{i}\underline{E}(f_i)$. We know that
any $f_i$ such that $\underline{e}^*(f_i,x_i^P,x_i^D) < M^*$ is
definitely not $\Gamma$-maximin, and thus can be removed from future
iterations for improving $M^*$ and $M_*$. Now, we do not have to wait
for $M_*$ and $M^*$ to have converged.  Indeed, we can already
immediately calculate $ \overline{e}_*(f_i,y_i^P,y_i^D)$ and
$\overline{e}^*(f_i,y_i^P,y_i^D)$, which bound $\overline{E}(f_i)$,
and also compute the duality gap
$\gamma_i=\overline{e}^*(f_i,y_i^P,y_i^D)-\overline{e}_*(f_i,y_i^P,y_i^D)$. We
then know that any $f_i$ such that
$\overline{e}_*(f_j,y_j^P,y_j^D) > M^*$ will be interval dominant
while any $f_i$ such that $\overline{e}^*(f_i,y_j^P,y_j^D) < M_*$ will
not be interval dominant. The entire process is repeated until there
are no more undetermined gambles or until all $\gamma_i$ and $M^*-M_*$
are less than a given tolerance. Additionally, we can also stop
computations in the first part of the algorithm as soon as $M^*-M_*$
is small enough, i.e. as soon as we know the $\Gamma$-maximin
value to sufficient precision. An algorithm that implements all of
these ideas is shown in \cref{alg:ID}.

\begin{algorithm}
	\caption{Interval dominance}\label{alg:ID}
	\begin{algorithmic}[1]
		\Require  a set of $k$ gambles $\mathcal{K} = \{f_1,\dots, f_k \}$;  a small number $\epsilon$; initial feasible states $x_i^P$ and $x_i^D$ for \linprogref{P1} and \linprogref{D1} corresponding to $f_i$ respectively; initial feasible states $y_i^P$ and $y_i^D$ for \linprogref{P2} and \linprogref{D2} corresponding to $f_i$ respectively;
		\Ensure an index set of $\opt_{\sqsupset}(\mathcal{K})$
		\Procedure{IntervalDominant3}{$\mathcal{K},\epsilon$}
		\State{$R \gets \{1,2\dots,k\}$} \Comment{index set of potentially $\Gamma$-maximin gambles}
		\State{$J \gets \{1,2\dots,k\}$} \Comment{index set of potentially interval dominant gambles}
		\State{$I \gets \emptyset$} \Comment{index set of known interval dominant gambles}
		\State $N \gets \emptyset$ \Comment{index set of known non interval dominant gambles}  
		\State $(M_*,M^*) \gets (-\infty,+\infty)$ \Comment{bounds on $\Gamma$-maximin value}
		\Repeat
                \State{\textbf{if $(M^{*}-M_{*}\ge\epsilon)$ then}}
		\State{\ \ \ $\forall i\in R$:}
		\State{\ \ \ \ \ \ $(x_i^P,x_i^D) \gets \phi(f_i,x_i^P,x_i^D)$} \Comment{next iteration to compute $\underline{E}(f_i)$}
		\State{\ \ \ \ \ \ $(\underline{\ell}_i,\underline{u}_i) \gets  (\underline{e}_*(f_i,x_i^P,x_i^D),\underline{e}^*(f_i,x_i^P,x_i^D))$}  \Comment{bounds on $\underline{E}(f_i)$}
		\State{\ \ \ $M_{*} = \max_{i\in R}\underline{\ell}_i$} \Comment{lower bound on $\Gamma$-maximin value}
		\State{\ \ \ $M^{*} = \max_{i\in R}\underline{u}_i$} \Comment{upper bound on $\Gamma$-maximin value}
		\State{\ \ \ $R \gets \{i\in R\colon \underline{u}_i \geq M_{*}\}$}		
                \State{\textbf{end if}}
		~
		\State{$\forall j \in J$:}
                \State{\ \ \ $(y_j^P,y_j^D) \gets \psi(f_j,y_j^P,y_j^D)$} \Comment{next iteration to compute $\overline{E}(f_j)$} 
		\State{\ \ \ $(\overline{\ell}_j,\overline{u}_j) \gets  (\overline{e}_*(f_j,y_j^P,y_j^D),\overline{e}^*(f_i,y_j^P,y_j^D))$}  \Comment{bounds on $\overline{E}(f_j)$}
		\State{\ \ \ $\gamma_j = \overline{u}_j-\overline{\ell}_j$} \Comment{duality gap}

		\State{$I \gets I	\cup \{j \in J\colon\overline{\ell}_j \geq M^* \}$} \Comment{$f_j$ is interval dominant}
		\State{$N \gets N \cup \{j \in J \colon \overline{u}_j < M_*\}$}
		\Comment{$f_j$ is not interval dominant}
		\State{$J = (I \cup N)^{\mathsf{c}}$}
		\Until{$(J = \emptyset)$ \textbf{ or } $\big( (\forall j \in J: \gamma_j < \epsilon)$\textbf{ and }$M^*-M_*<\epsilon \big)$} 
		\State $I \gets I	\cup J$
		\State \Return $I$  
		\EndProcedure
	\end{algorithmic}
\end{algorithm}

At first sight, this algorithm has potential to be faster because it
provides an opportunity to eliminate interval dominated gambles before
the $\Gamma$-maximin value is known to full precision. However, the
algorithm may need to do more iterations to calculate the upper
previsions to sufficient precision as long as the bounds on the
$\Gamma$-maximin value remain too wide. Additionally, prior sorting of
gambles does not help here, because all gambles are treated
simultaneously, and not in sequence, as in \cref{alg:ID2}. Our
numerical experiments have found that \cref{alg:ID2} beats
\cref{alg:ID}. Basically, by postponing our evaluation of the
$\Gamma$-maximin value, the algorithm cannot eliminate gambles early
on, and needs to do more overall iterations in the end, for most of
the benchmarking problems considered.

In our numerical experiments, we observed that in the second part of
\cref{alg:ID}, when we compare the lower bound of $\overline{E}(f)$
against $M^*$ and compare the upper bound of $\underline{E}(f)$
against $M_*$, we rarely determine whether $f$ is interval dominant or
not before we can identify the (usually unique) $\Gamma$-maximin
gamble. In contrast, once we know the $\Gamma$-maximin gamble, then we
observed that we can very quickly determine whether $f$ is interval
dominant or not with much fewer iterations.
This leads us to try one more final idea: use \cref{alg:maximin3}
until a $\Gamma$-maximin gamble is identified but still allow
$M^*-M_*\ge\epsilon$, and then proceed as with \cref{alg:ID}. This is
shown in \cref{alg:ID3}. 
This version of the algorithm is simpler and 
easier to implement. However, we found that the performance of \cref{alg:ID,alg:ID3} is still nearly identical. Therefore, \cref{alg:ID3} still does not beat \cref{alg:ID2}, as we shall see in \cref{sec:benchmark}.

\begin{algorithm}
	\caption{Interval dominance}\label{alg:ID3}
	\begin{algorithmic}[1]
		\Require  a set of $k$ gambles $\mathcal{K} = \{f_1,\dots, f_k \}$;  a small number $\epsilon$; initial feasible states $x_i^P$ and $x_i^D$ for \linprogref{P1} and \linprogref{D1} corresponding to $f_i$ respectively; initial feasible states $y_i^P$ and $y_i^D$ for \linprogref{P2} and \linprogref{D2} corresponding to $f_i$ respectively;
		\Ensure a set of $\Gamma$-maximin gambles; an index set of $\opt_{\sqsupset}(\mathcal{K})$
		\Procedure{IntervalDominant3}{$\mathcal{K},\epsilon$}
		\State{\textbf{Stage I:} Apply \textsc{$\Gamma$-maximin3} to obtain $R$, $M_*$ and $M^*$} 
		
		\Comment{$R$ is an index set of the $\Gamma$-maximin gambles, $M_*$ is a lower bound for the $\Gamma$-maximin value and $M^*$ is an upper bound for the $\Gamma$-maximin value}
		~	
		\State{\noindent \textbf{Stage II: Interval dominance}}
		\State $I\gets R$ \Comment{index set of known interval dominant gambles}  
		\State $N \gets \emptyset$ \Comment{index set of known non interval dominant gambles}  
		\State{$J \gets R^{\mathsf{c}}$} \Comment{index set of potentially interval dominant gambles}
		\Repeat
                \State{$\forall j\in J$:}
		\State{\ \ \ $(y_j^P,y_j^D) \gets \psi(f_j,y_j^P,y_j^D)$} \Comment{next iteration to compute $\overline{E}(f_j)$} 
		\State{\ \ \ $(\ell_j,u_j) \gets (\overline{e}_*(f_j,y_j^P,y_j^D),\overline{e}^*(f_j,y_j^P,y_j^D))$}  \Comment{bounds on $\overline{E}(f_j)$}
		\State{\ \ \ $\gamma_j = u_j-\ell_j$} \Comment{duality gap}
		\State{\textbf{if } $M^{*} - M_{*} \ge \epsilon$ \textbf{ then } } 
		\State{\ \ \ $\forall r \in R$:}
                \State{\ \ \ \ \ \ $(x_r^P,x_r^D) \gets \phi(f_r,x_r^P,x_r^D)$}\Comment{continue compute $\underline{E}(f_r)$}
		\State{\ \ \ \ \ \ $(\ell_r, u_r)\gets (\underline{e}_*(f_r,x_r^P,x_r^D), \underline{e}^*(f_r,x_r^P,x_r^D))$}  
		\State{\ \ \ $M_{*} = \max_{r\in R}\ell_r$}
		\State{\ \ \ $M^{*} = \max_{r\in R}u_r$}
		\State{\textbf{end if}}
		\State{$I \gets I	\cup \{j \in J\colon\ell_j \geq M^*\}$} \Comment{$f_j$ is interval dominant}
		\State{$N \gets N \cup \{j \in J \colon u_j < M_*\}$}
		\Comment{$f_j$ is not interval dominant}
		\State{$J = (I \cup N)^{\mathsf{c}}$}
		\Until{$(J = \emptyset)$ \textbf{ or } $\big((\forall j \in J: \gamma_j < \epsilon)$\textbf{ and }$M^*-M_*<\epsilon\big)$} 
		\State $I \gets I	\cup J$ 
		\State \Return $(R,I)$
		\EndProcedure
	\end{algorithmic}
\end{algorithm}

We note that many of the computations can be parallelized, which could
be exploited to reduce the running time of the algorithm.
In particular, the parts under the
universal quantifiers (i.e. $\forall i\in R$ and $\forall j\in J$) are
independent, and can therefore in principle run
simultaneously. However, after each parallel step, the outcomes must
be combined to continue with the next step, to determine the updated
$R$ and $J$, as well as the updated bounds $M_*$ and $M^*$ on the
$\Gamma$-maximin value, so some communication between the parallel
processes is needed between iterations.
From a runtime perspective,
\cref{alg:ID} might still be an attractive option to consider even though
it takes more computational effort.
For this paper however, we have focused our implementations and benchmarking
on the serial running time to reflect the overall computational effort,
and leave a study of parallel algorithms,
where interesting trade-offs are to be had compared to serial algorithms,
to future work.

To summarise,  \cref{alg:maximin1,alg:maximin2,alg:maximin3}  are for finding $\Gamma$-maximin gambles, \cref{alg:maximax1,alg:maximax2,alg:maximax3} are for finding $\Gamma$-maximax gambles and \cref{alg:ID1,alg:ID2,alg:ID,alg:ID3}  are for finding interval dominant gambles. We will benchmark these algorithms in the next section. 

\section{Benchmarking results}\label{sec:benchmark}
\subsection{$\Gamma$-maximin and $\Gamma$-maximax}
In this section, we will compare our algorithms for $\Gamma$-maximin and $\Gamma$-maximax for the case that $|\mathcal{K}| = 2^i$ for $i \in \{1,2,\dots,10 \}$ and $|\Omega| = 2^j$ for $j \in \{1,2,\dots,10 \}$. This covers a good range of the experiment sizes that we expect to see in real-world applications of decision making (for example \citep{2012:TROFFAES,2018:Troffaes,2019:DENOEUX}), where $2^{10}$ gambles (or decisions) and $2^{10}$ outcomes should be plenty big enough for real-world decision making.

For each case, we first generate a lower prevision $\underline{P}$, where $|\dom \underline{P}| = 2^4$, on a finite domain which avoids sure loss as follows. To do so, we use \citep[algorithm 2]{2018:Nakharutai:Troffaes:Caiado} with $2^4$ coherent previsions to generate a lower prevision $\underline{E}$ on the set of all gambles, that avoids sure loss.  Next, we use \citep[stages 1 and 2 in algorithm 4]{2018:Nakharutai:Troffaes:Caiado} to restrict $\underline{E}$ to a lower prevision $\underline{P}$ that avoids sure loss, with a given finite size of domain. We also consider $|\dom \underline{P}| = 2^n$ for $n \in \{6,8,10\}$.
Next, to generate gambles $f_i \in \mathcal{K}$, for each $\omega$ and $i$, we sample $f_i(\omega)$ uniformly from interval $[0, 1]$.

For each generated set of gambles $\mathcal{K}$, we perform \cref{alg:maximin1,alg:maximin2,alg:maximin3} to find $\Gamma$-maximin gambles and \cref{alg:maximax1,alg:maximax2,alg:maximax3} to find $\Gamma$-maximax gambles. We apply \citep[\S 4.2]{2018:Nakharutai:Troffaes:Caiado} to find feasible starting points and implement early stopping criteria in the primal-dual method corresponding to different cases.
Specifically, the improved primal-dual method solves  
linear programs \linprogref{P1} and \linprogref{D1} inside \cref{alg:maximin2,alg:maximin3}, and solves  
linear programs \linprogref{P2} and \linprogref{D2} inside \cref{alg:maximax2,alg:maximax3}.
Based on the implementation used in \citep{2018:Nakharutai:Troffaes:Caiado}, we wrote our own implementation of the improved primal-dual method in MATLAB (R2019b) \citep{MATLAB:2019} rather than using  MATLAB's existing library. This is because the existing library does not allow us to inspect the algorithms at different iterations easily. Since we want to compare different algorithms that work in different ways, we need to ensure that the only difference is in our part of the algorithms. For \cref{alg:maximin1,alg:maximax1} where we solve linear programs inside these algorithms by the standard primal-dual method, to have a same advantage, we also wrote our own implementation of the standard primal-dual method instead of using the primal-dual method in MATLAB's toolbox.

\begin{figure}
	\centering
	\setlength{\tabcolsep}{2pt}
	\newcolumntype{C}{>{\centering\arraybackslash} m{0.48\linewidth} }
	
	\begin{tabular}{m{0.5em}CC}
		&
		Algorithms for $\Gamma$-maximin
		&
		Algorithms for $\Gamma$-maximax
		\\
		\rotatebox[origin=l]{90}{$|\mathcal{K}| = 2^4$}
		&
		\includegraphics[width=\hsize, trim={2cm 0 2cm 0},clip]{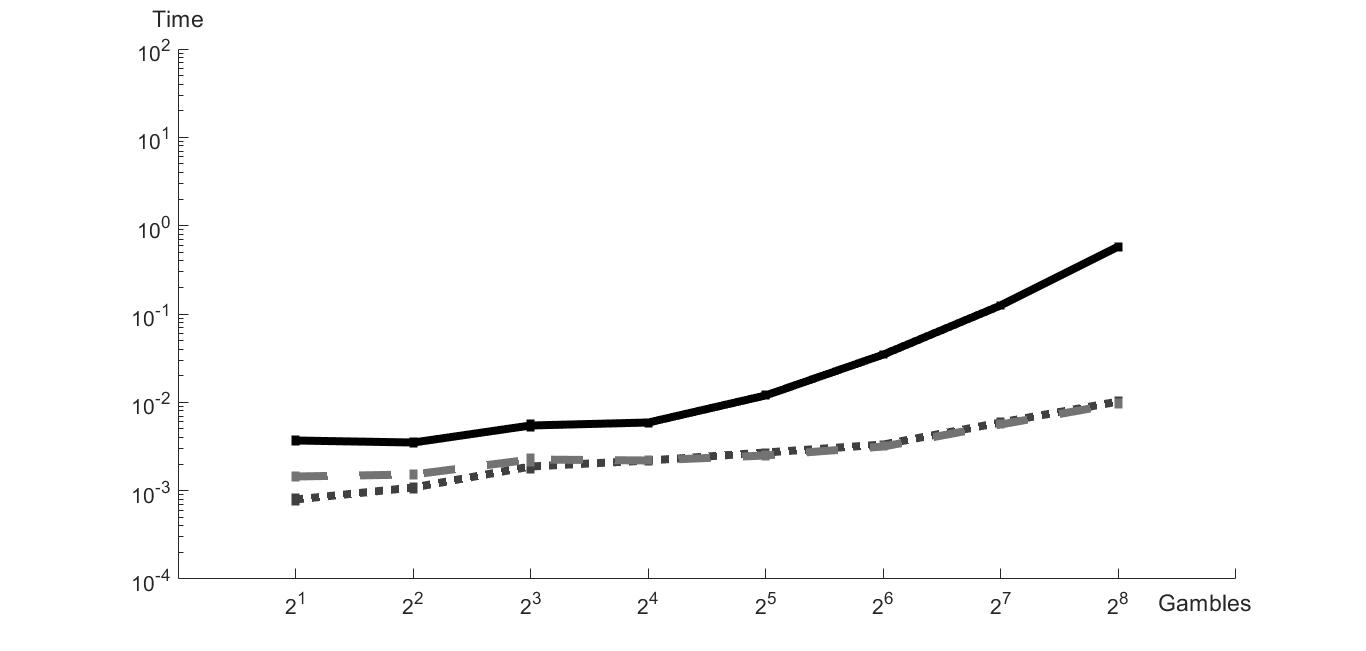}
		&
		\includegraphics[width=\hsize, trim={2cm 0 2cm 0},clip]{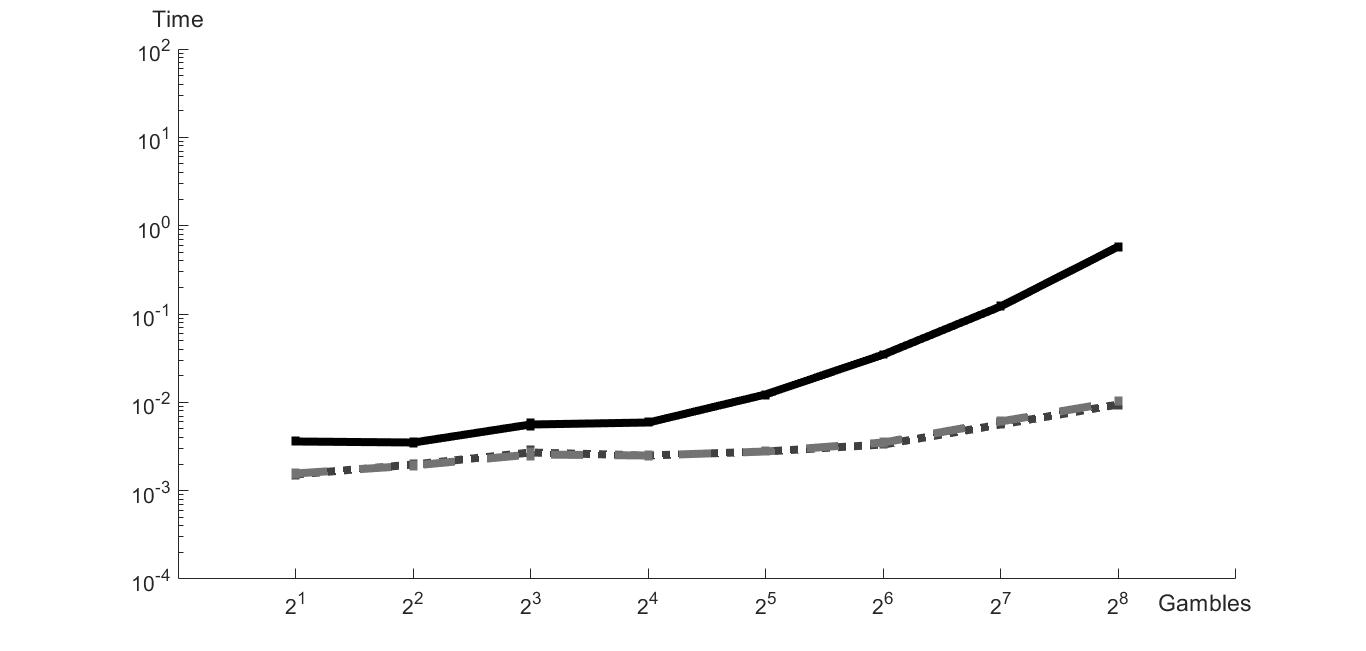}
		\\
		\rotatebox[origin=l]{90}{$|\mathcal{K}| = 2^6$}
		&
		\includegraphics[width=\hsize, trim={2cm 0 2cm 0},clip]{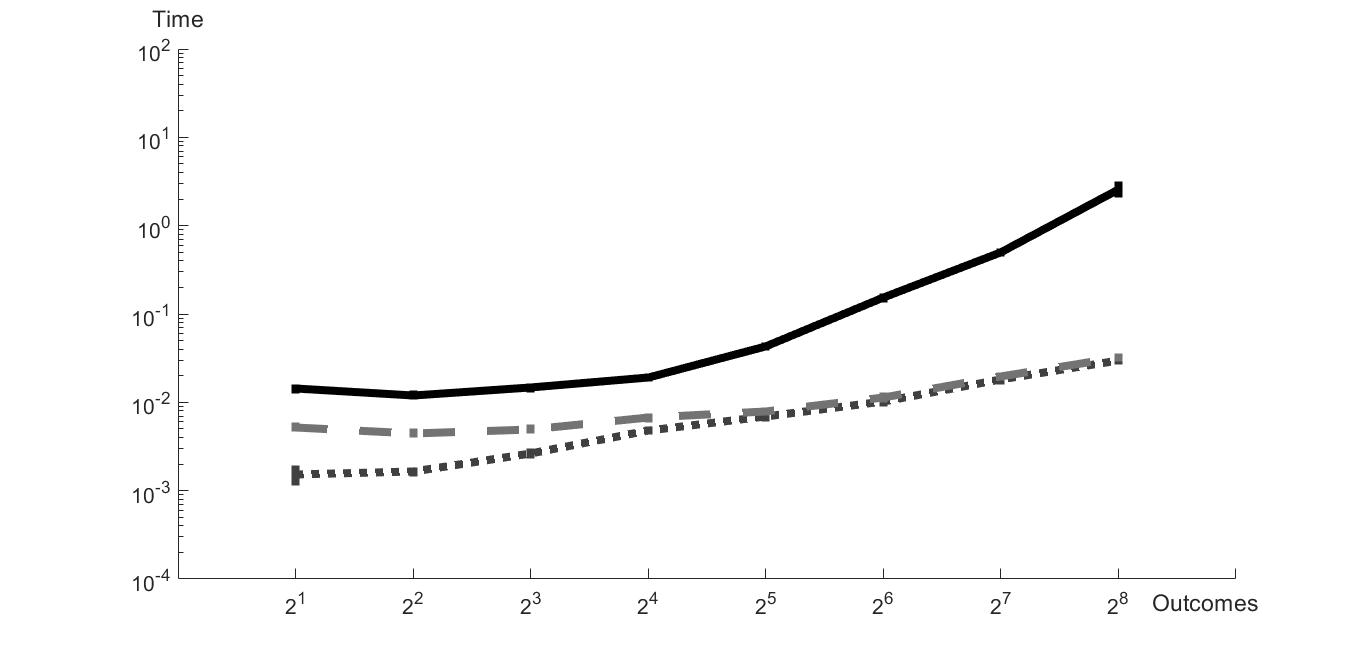}
		&
		\includegraphics[width=\hsize, trim={2cm 0 2cm 0},clip]{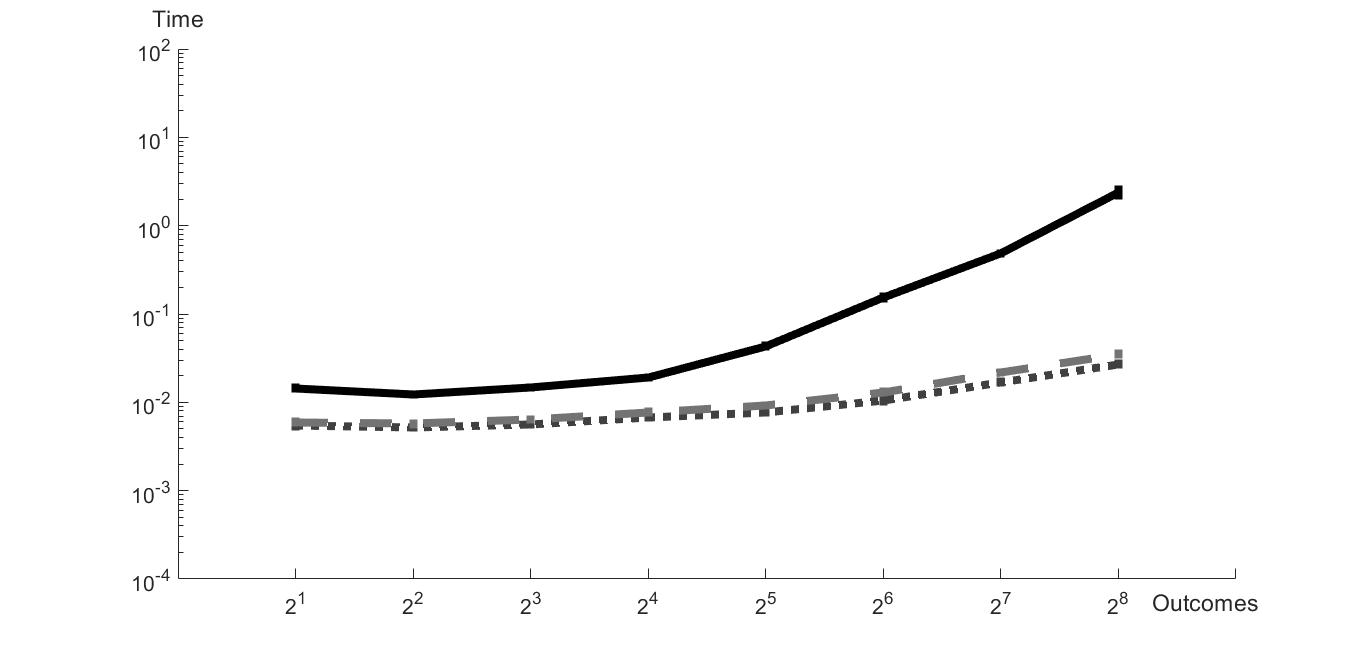}
		\\
		\rotatebox[origin=l]{90}{$|\mathcal{K}| = 2^8$}
		&
		\includegraphics[width=\hsize, trim={2cm 0 2cm 0},clip]{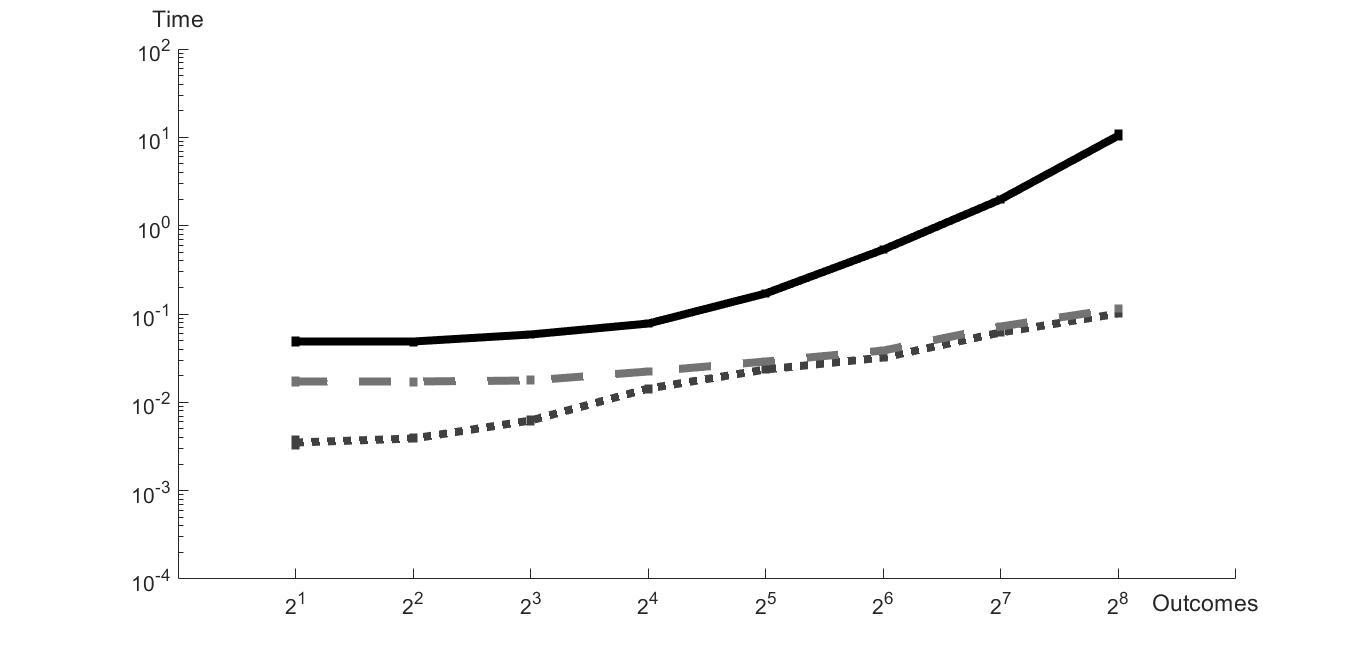}
		&
		\includegraphics[width=\hsize, trim={2cm 0 2cm 0},clip]{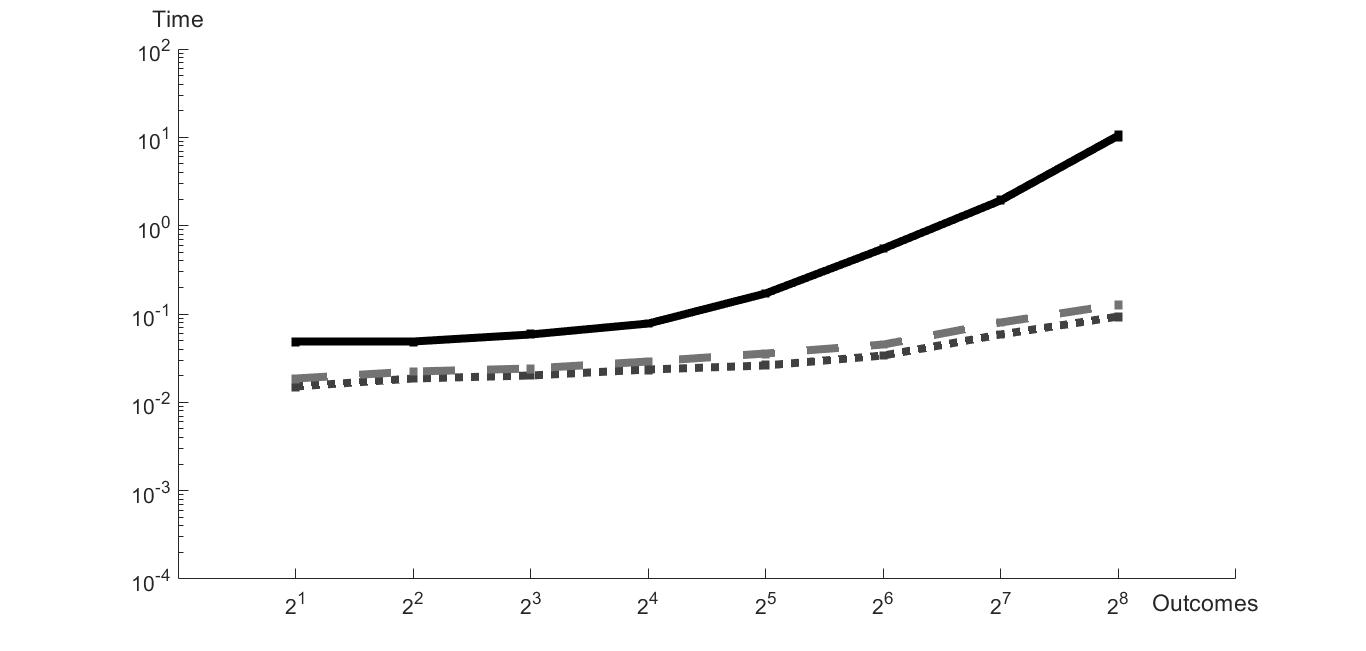}
		\\
		\rotatebox[origin=l]{90}{$|\mathcal{K}| = 2^{10}$}
		&
		\includegraphics[width=\hsize, trim={2cm 0 2cm 0},clip]{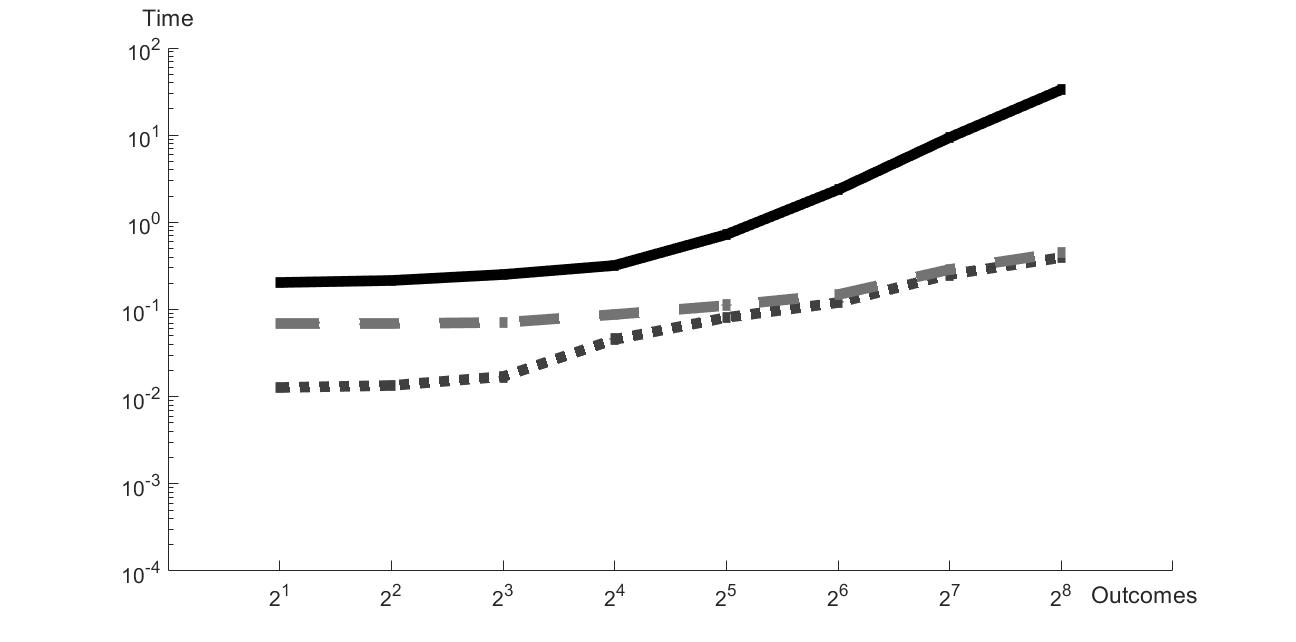}
		&
		\includegraphics[width=\hsize, trim={2cm 0 2cm 0},clip]{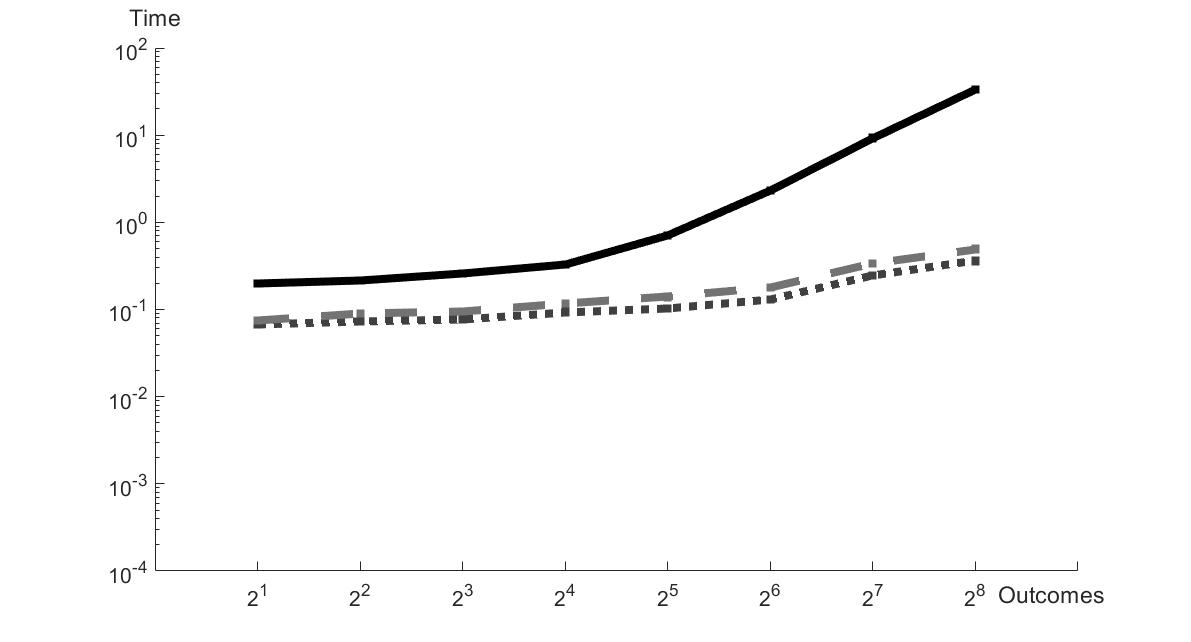}
		\\
		&
		\includegraphics[width=0.3\linewidth]{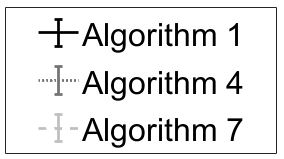}
		&
		\includegraphics[width=0.3\linewidth]{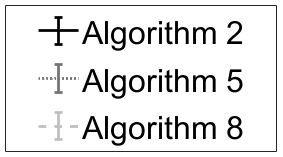}		
	\end{tabular}

	\caption{Comparison plots of the average total running time for different algorithms for $\Gamma$-maximin (the left column) and $\Gamma$-maximax (the right column).   Each row shows a different number of gambles with varying numbers of outcomes and fixed $|\dom \underline{P}| = 2^4$. The labels indicate different algorithms.}
	\label{fig:plot1}
\end{figure}

\begin{figure}
	\centering
	\setlength{\tabcolsep}{2pt}
	\newcolumntype{C}{>{\centering\arraybackslash} m{0.48\linewidth} }
	
	\begin{tabular}{m{0.5em}CC}
		&
		Algorithms for $\Gamma$-maximin
		&
		Algorithms for $\Gamma$-maximax
		\\
		\rotatebox[origin=l]{90}{$|\Omega| = 2^4$}
		&
		\includegraphics[width=\hsize, trim={2cm 0 2cm 0},clip]{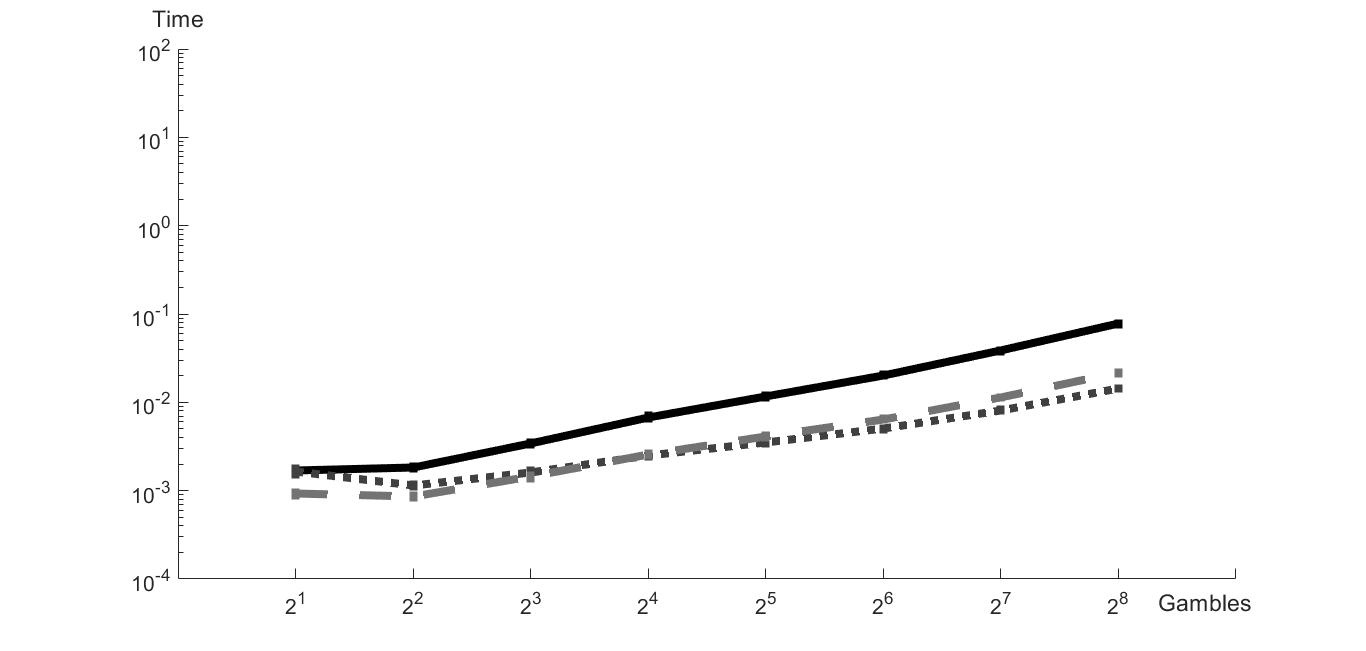}
		&
		\includegraphics[width=\hsize, trim={2cm 0 2cm 0},clip]{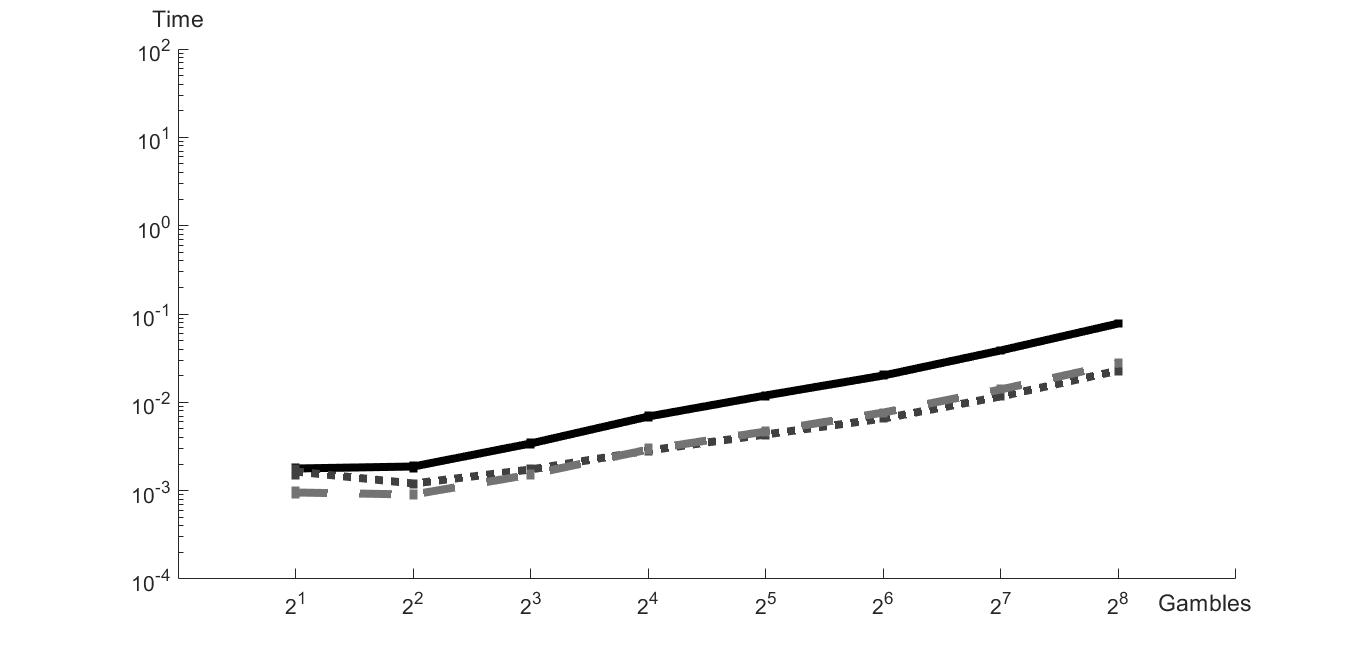}
		\\
		\rotatebox[origin=l]{90}{$|\Omega| = 2^6$}
		&
		\includegraphics[width=\hsize, trim={2cm 0 2cm 0},clip]{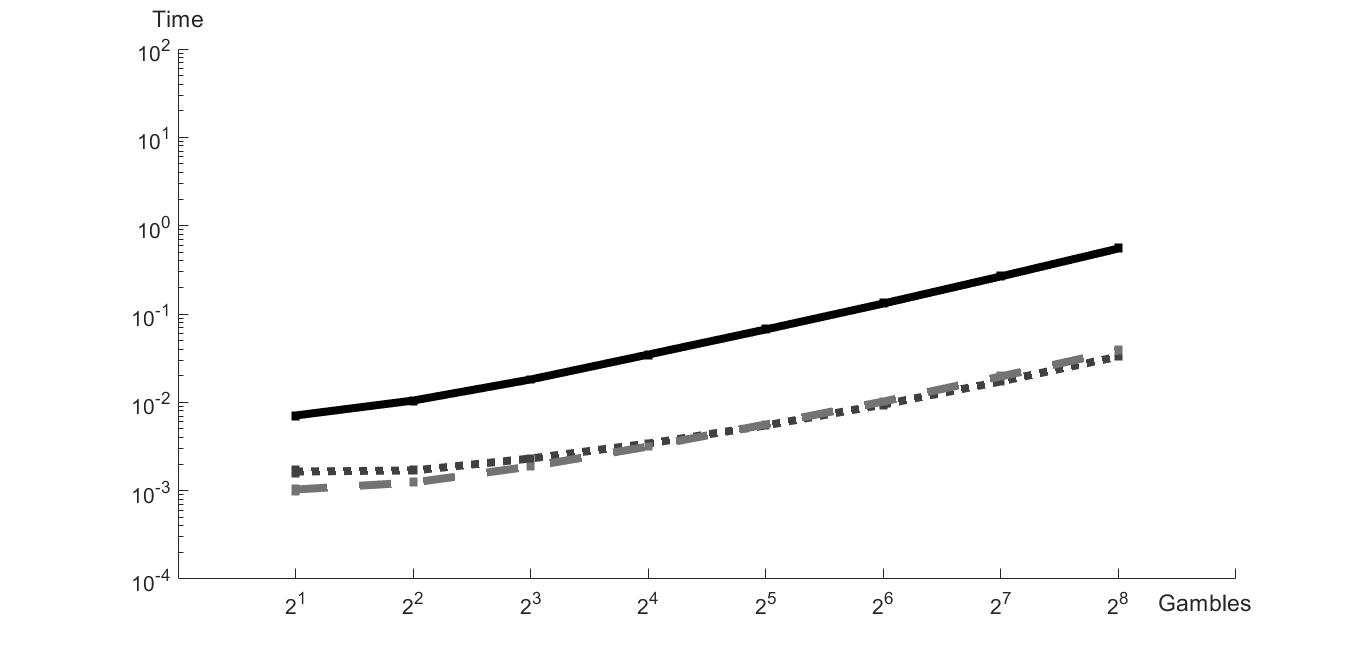}
		&
		\includegraphics[width=\hsize, trim={2cm 0 2cm 0},clip]{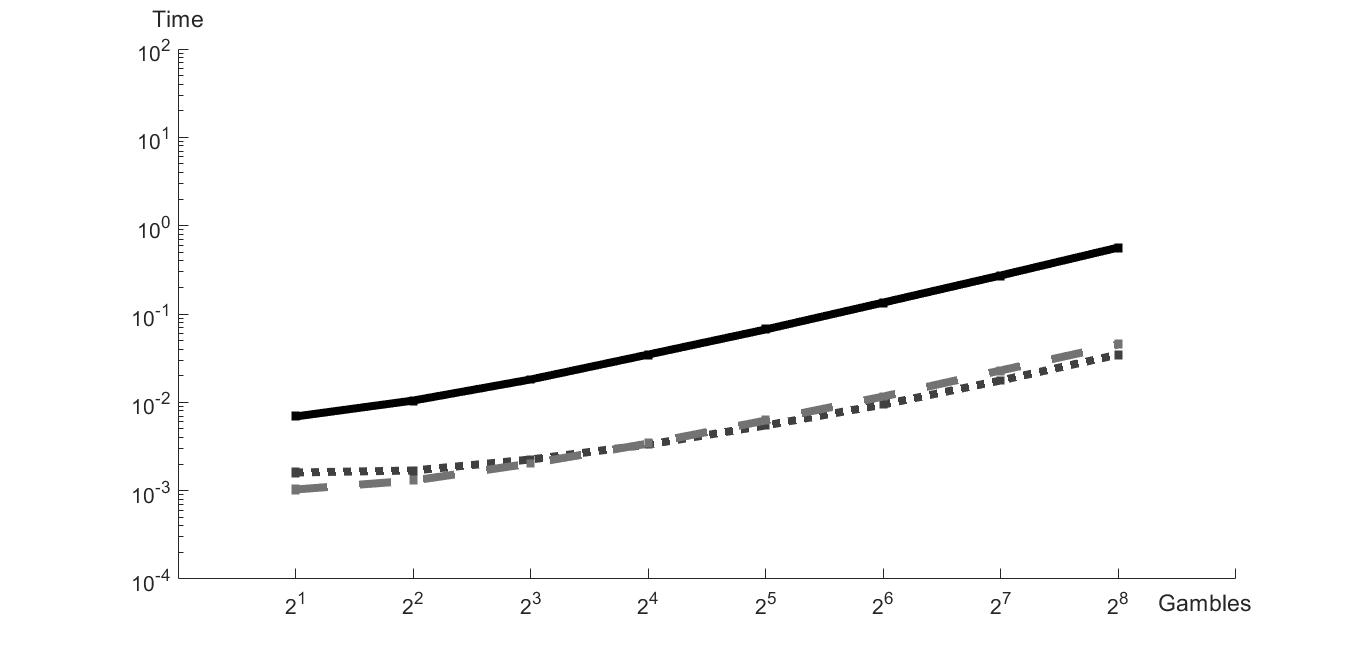}
		\\
		\rotatebox[origin=l]{90}{$|\Omega| = 2^8$}
		&
		\includegraphics[width=\hsize, trim={2cm 0 2cm 0},clip]{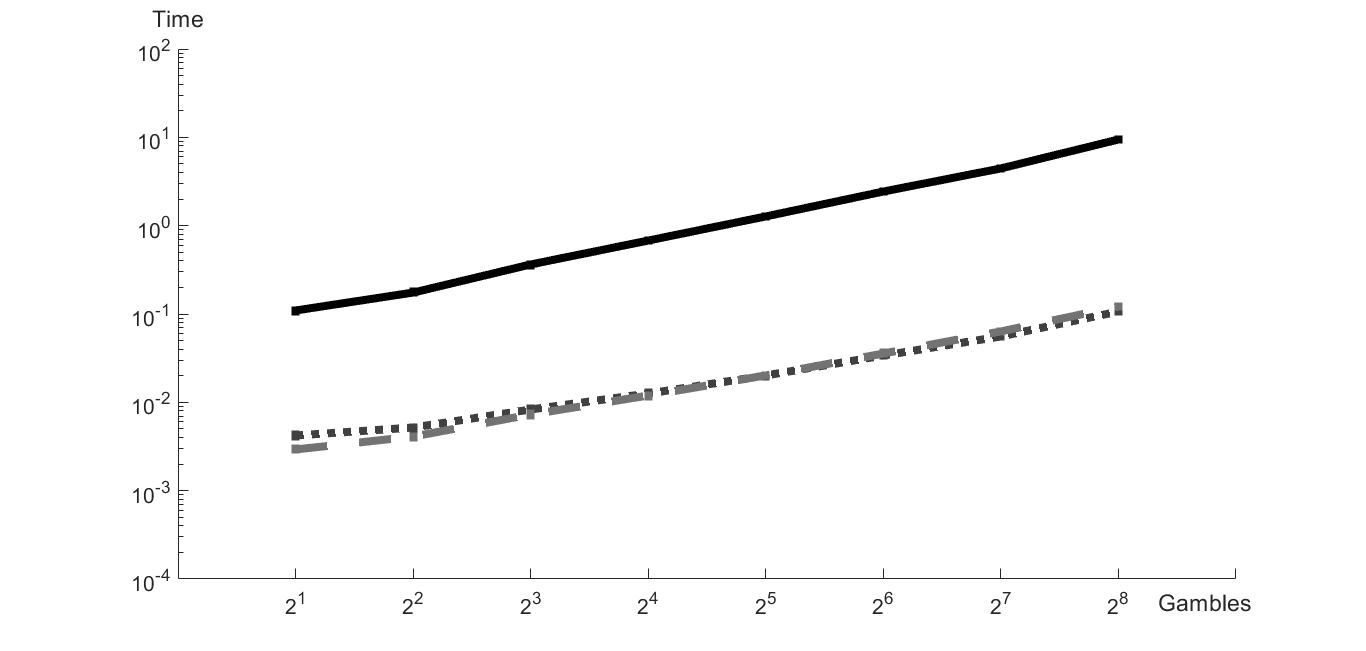}
		&
		\includegraphics[width=\hsize, trim={2cm 0 2cm 0},clip]{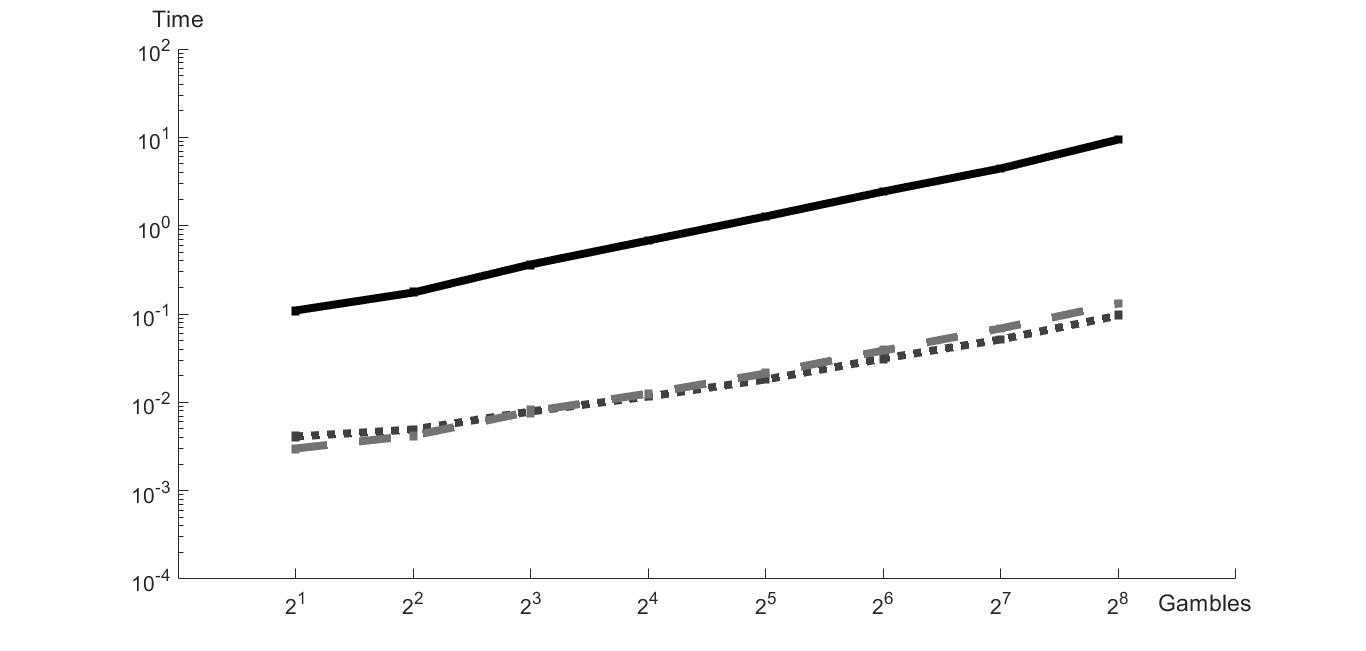}
		\\
		\rotatebox[origin=l]{90}{$|\Omega| = 2^{10}$}
		&
		\includegraphics[width=\hsize, trim={2cm 0 2cm 0},clip]{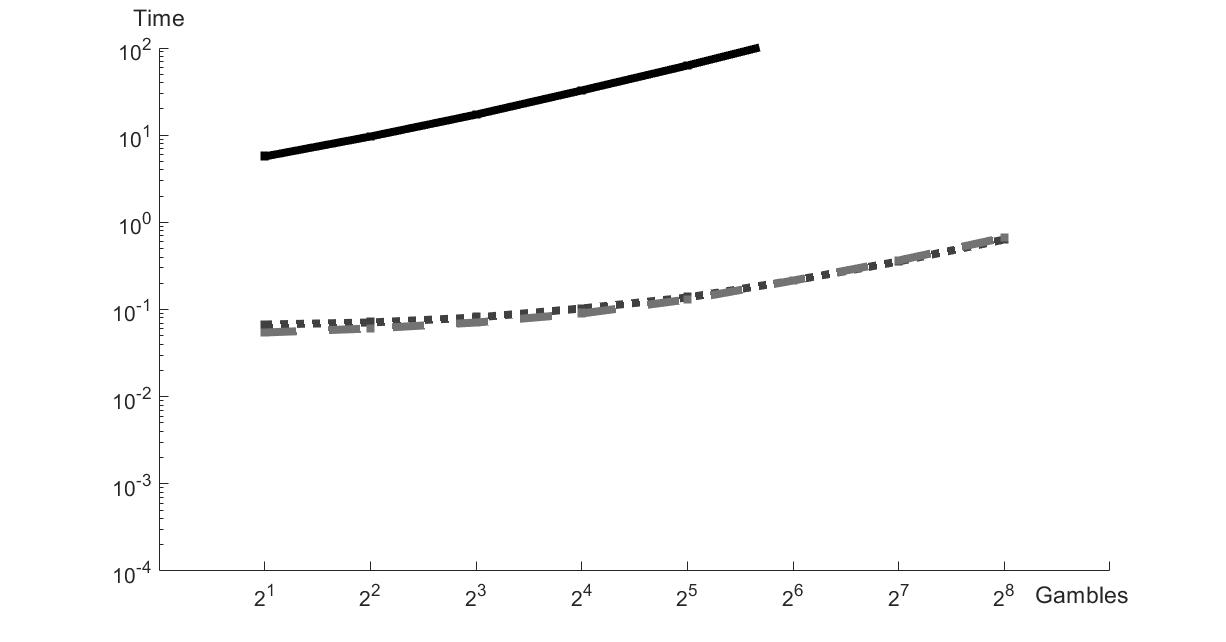}
		&
		\includegraphics[width=\hsize, trim={2cm 0 2cm 0},clip]{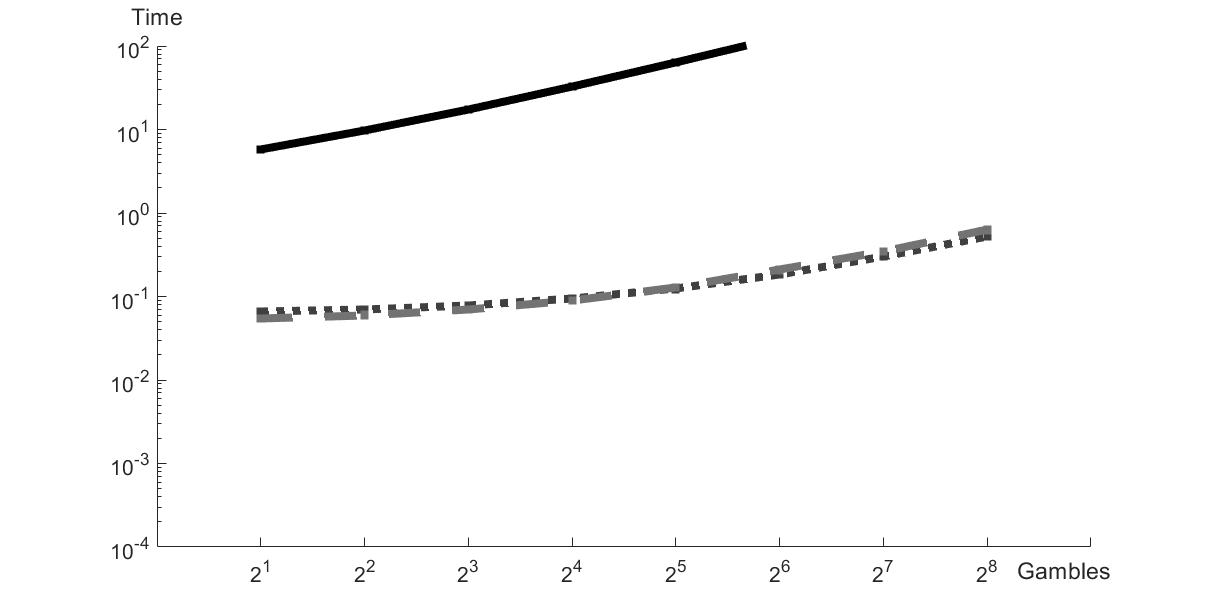}
		\\
		&
		\includegraphics[width=0.3\linewidth]{label_min}
		&
		\includegraphics[width=0.3\linewidth]{label_max}		
	\end{tabular}

	\caption{Comparison plots of the average total running time of different algorithms for $\Gamma$-maximin (the left column) and $\Gamma$-maximax (the right column).   Each row shows a different number of outcomes with varying numbers of gambles and fixed $|\dom \underline{P}| = 2^4$.  The labels indicate different algorithms.}
	\label{fig:plot2}
\end{figure}

For each algorithm, we record the total running time spent. In particular, for \cref{alg:maximin2,alg:maximax2} this includes the time to sort gambles with respect to their expectations as in \cref{eq:E_p(i-f_j)}. To sort gambles, we use quicksort which is a built-in function available in MATLAB \citep{MATLAB:2019}. Even though quicksort is an obvious choice for sorting, we state it for completeness, as there are other available functions for sorting.
This process is repeated 1000 times. A summary of the results is presented  in \cref{fig:plot1,fig:plot2}.

\begin{figure}
	\centering
	\setlength{\tabcolsep}{2pt}
	\newcolumntype{C}{>{\centering\arraybackslash} m{0.48\linewidth} }
	
	\begin{tabular}{m{0.5em}CC}
		&
		Algorithms for $\Gamma$-maximin
		&
		Algorithms for $\Gamma$-maximax
		\\
		\rotatebox[origin=l]{90}{\text{\tiny$|\dom \underline{P}| = 2^4$}}
		&
		\includegraphics[width=\hsize, trim={2cm 0 2cm 0},clip]{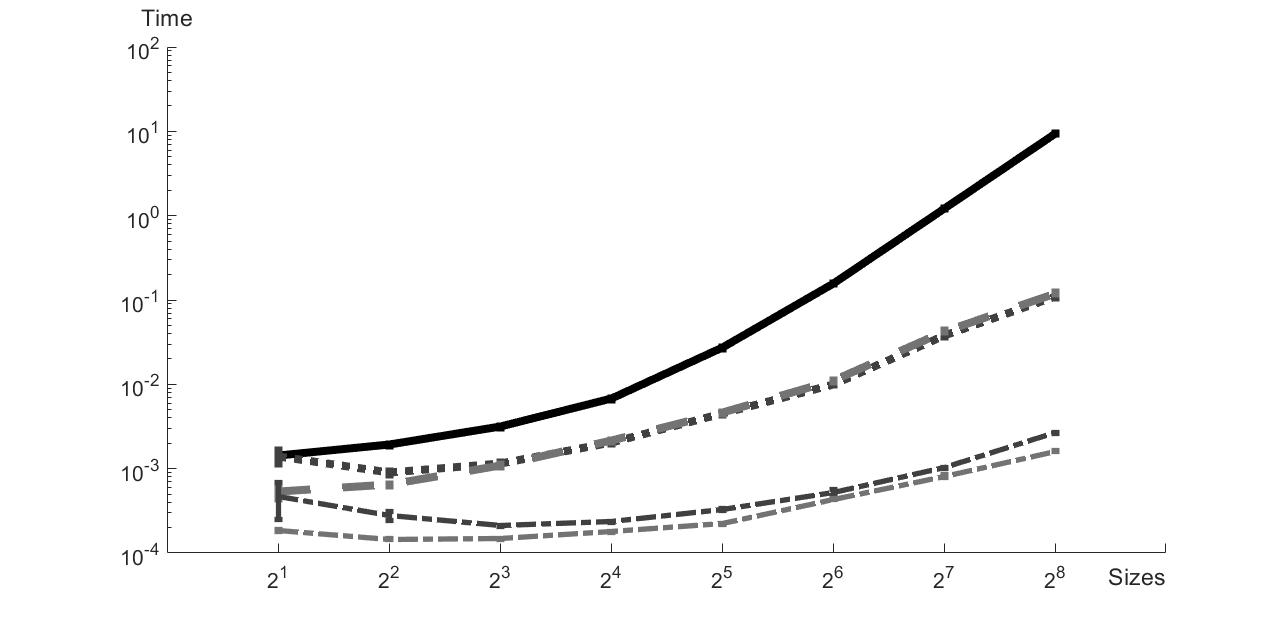}
		&
		\includegraphics[width=\hsize, trim={2cm 0 2cm 0},clip]{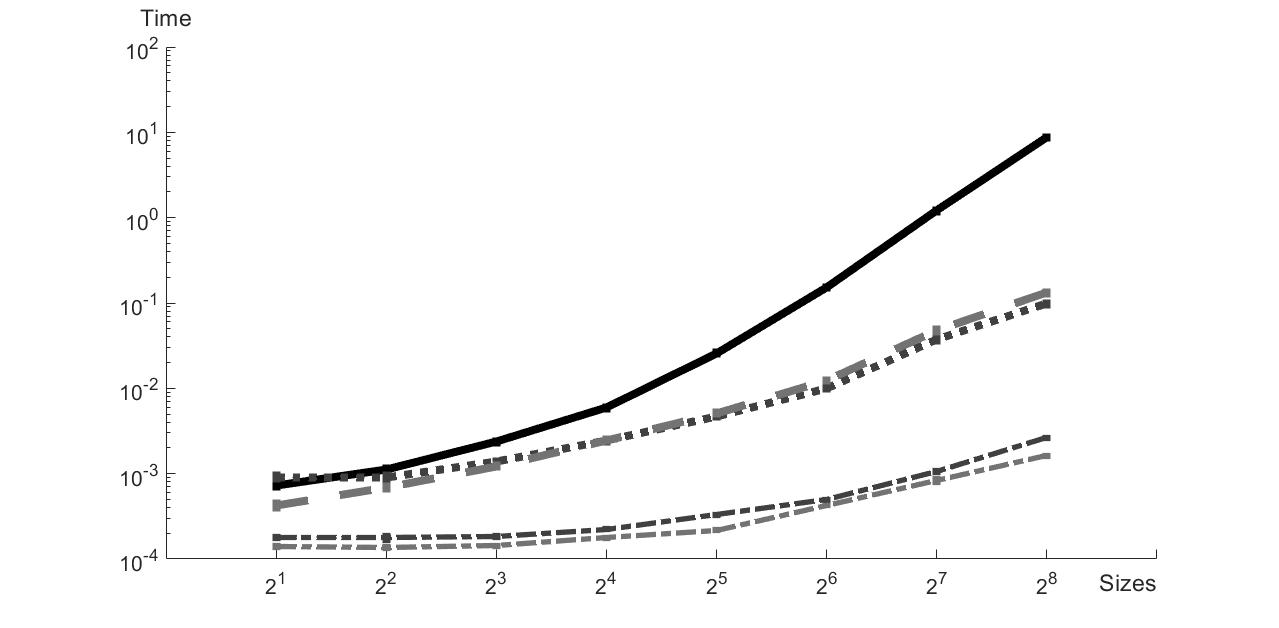}
		\\
		\rotatebox[origin=l]{90}{\text{\tiny$|\dom \underline{P}|= 2^6$}}
		&
		\includegraphics[width=\hsize, trim={2cm 0 2cm 0},clip]{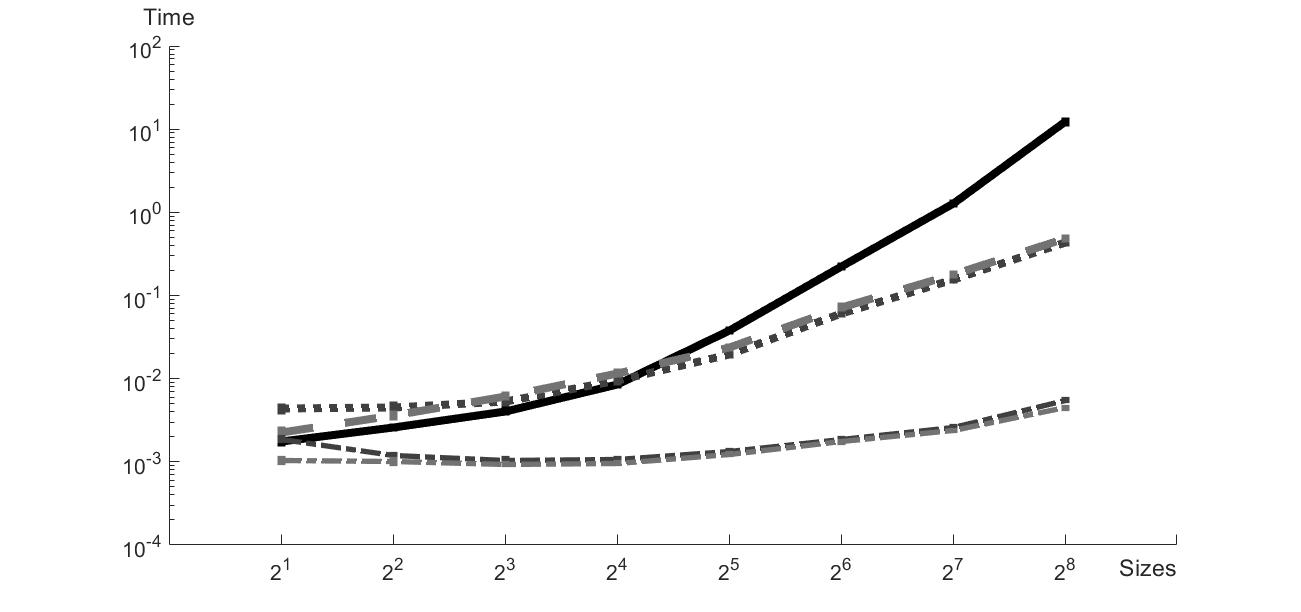}
		&
		\includegraphics[width=\hsize, trim={2cm 0 2cm 0},clip]{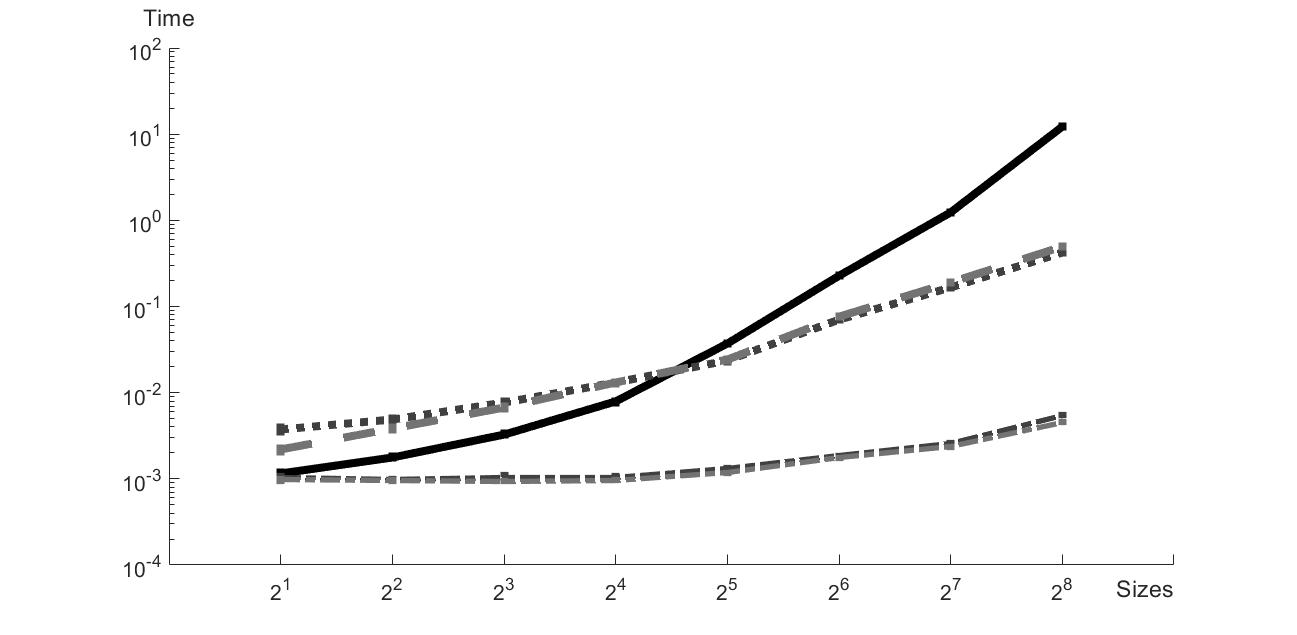}
		\\
		\rotatebox[origin=l]{90}{\text{\tiny$|\dom \underline{P}|= 2^8$}}
		&
		\includegraphics[width=\hsize, trim={2cm 0 2cm 0},clip]{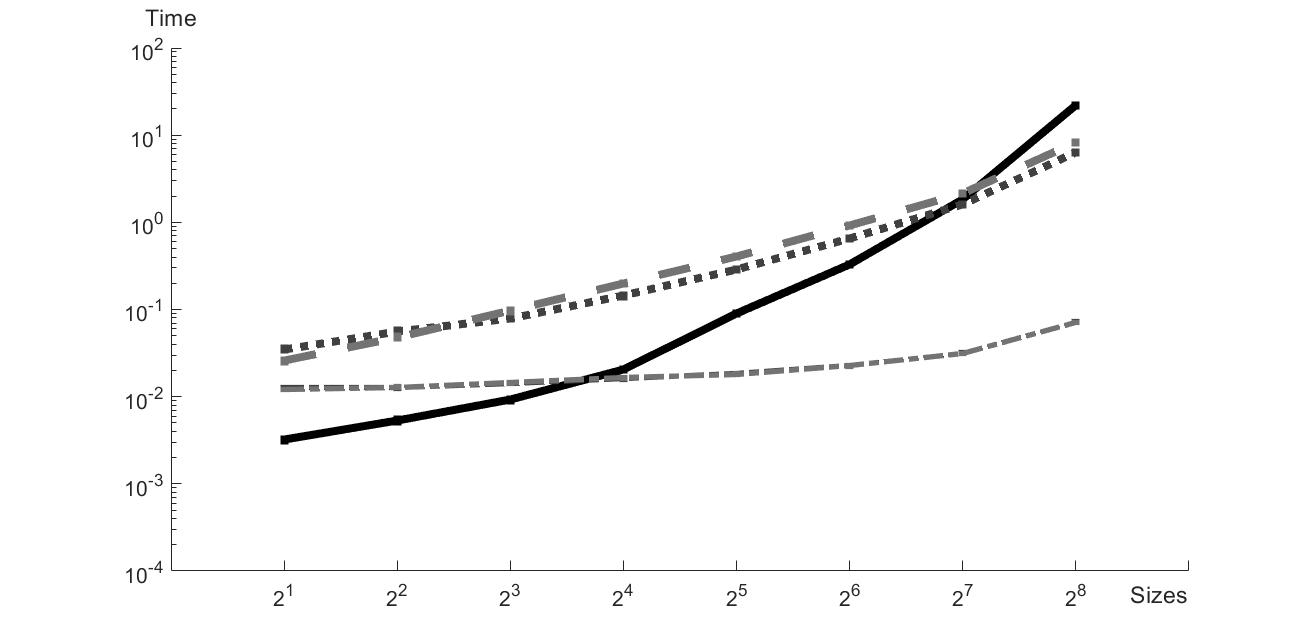}
		&
		\includegraphics[width=\hsize, trim={2cm 0 2cm 0},clip]{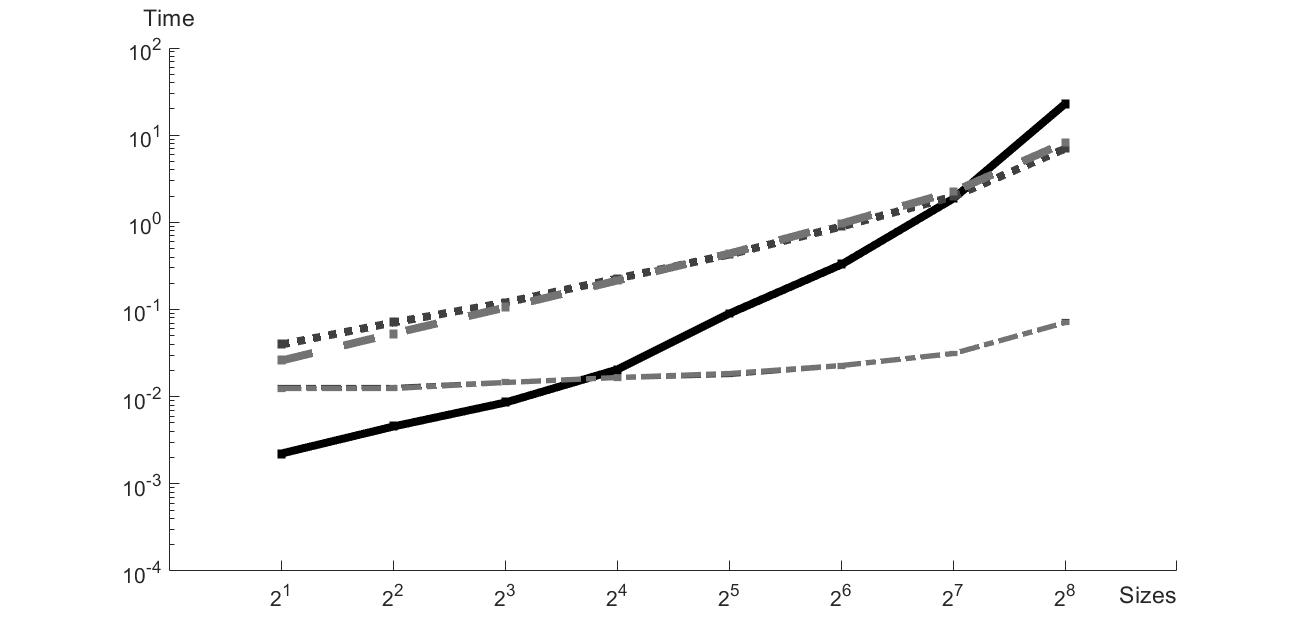}
		\\
				\rotatebox[origin=l]{90}{\text{\tiny$|\dom \underline{P}|= 2^{10}$}}
		&
		\includegraphics[width=\hsize, trim={2cm 0 2cm 0},clip]{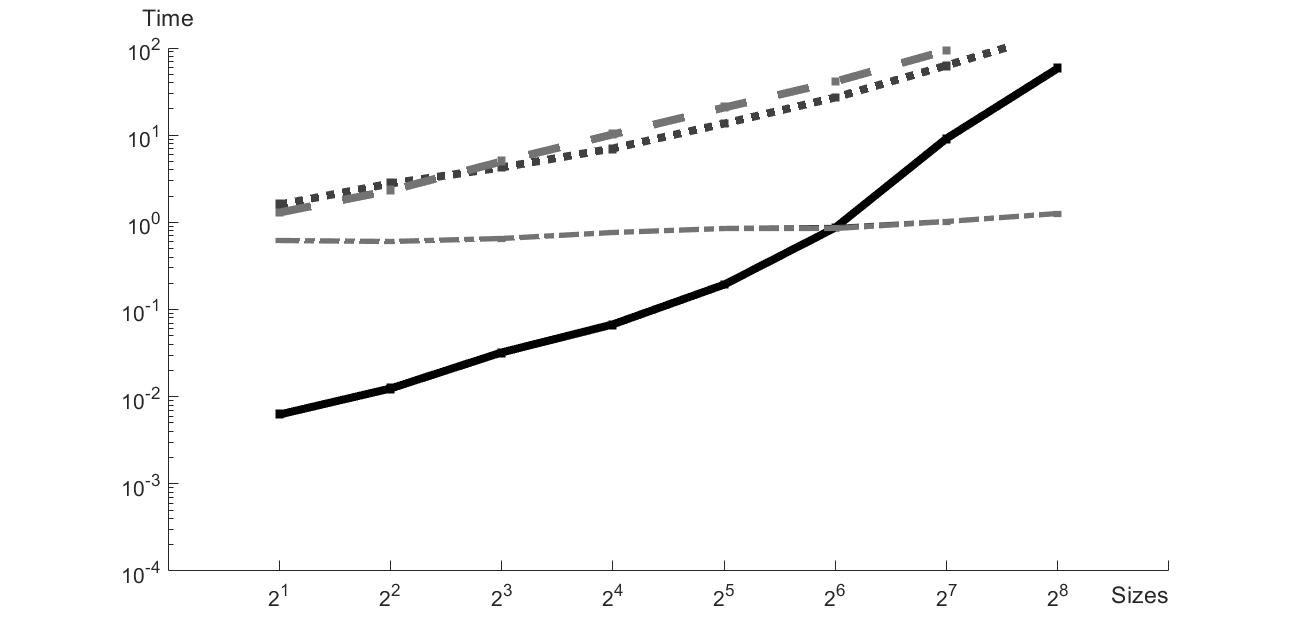}
		&
		\includegraphics[width=\hsize, trim={2cm 0 2cm 0},clip]{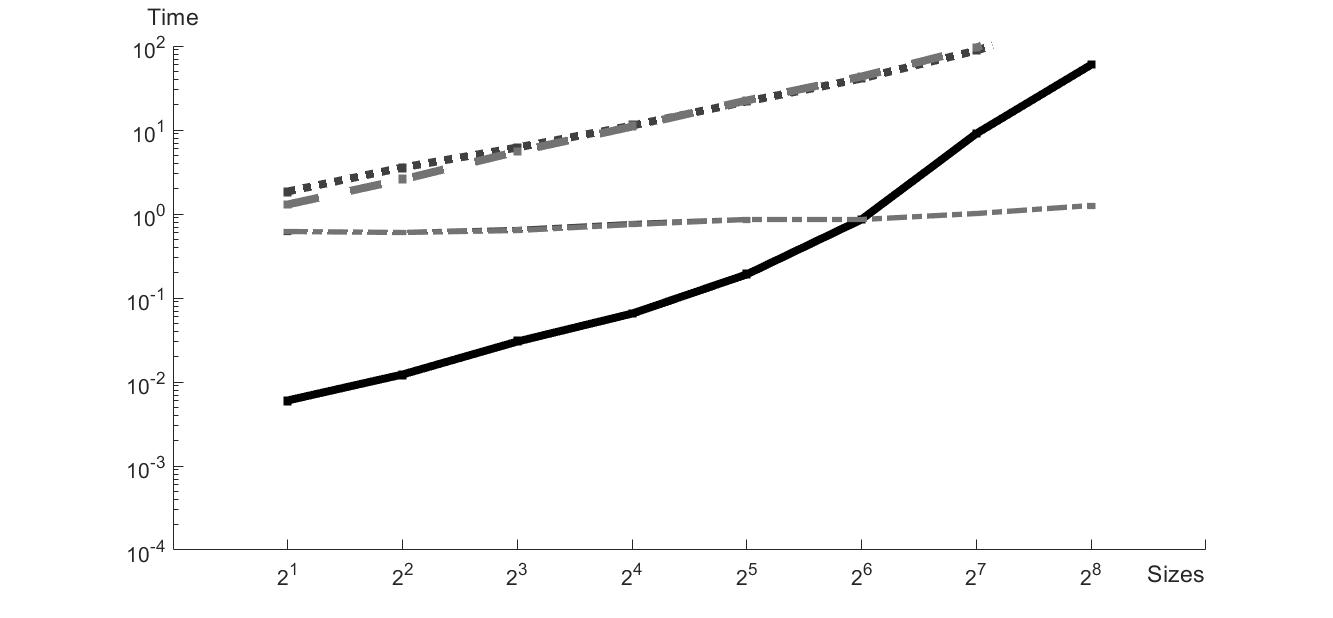}
		\\
		&
		\includegraphics[width=0.3\linewidth]{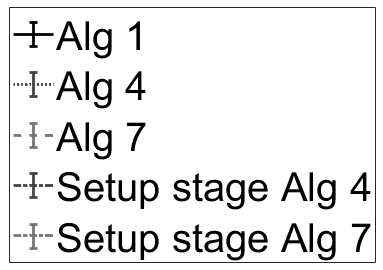}
		&
		\includegraphics[width=0.3\linewidth]{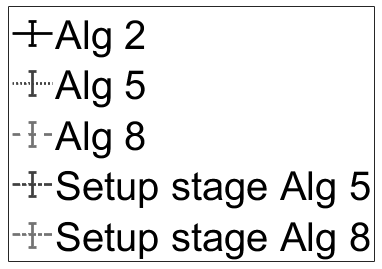}		
	\end{tabular}

	\caption{Comparison plots of the average running time for different improved algorithms for $\Gamma$-maximin (the left column) and $\Gamma$-maximax (the right column). Each row shows a different number of $|\dom \underline{P}|$ with varying numbers of gambles and outcomes, i.e. $|\Omega| = |\mathcal{K}| = 2^i$, for $i \in \{1,2,\dots, 8 \}$. The labels indicate different algorithms.}
	\label{fig:plot_domP}
\end{figure}

\Cref{fig:plot1,fig:plot2} present the average total running time of different algorithms for $\Gamma$-maximin in the left column and $\Gamma$-maximax in the right column. In \cref{fig:plot1}, each row presents a different number of gambles and the horizontal axis presents the number of outcomes. In \cref{fig:plot2}, each row presents a different number of outcomes and the horizontal axis presents the number of gambles.  In both figures, the vertical axis shows the running time which is averaged over 1000 random generated decision problems. 
The error bars on the figures show approximate 95\% confidence intervals on the average running time. They are barely visible except in few cases, e.g. in \cref{fig:plot_domP}, because of the large sample size.

We also consider the effect of the size of $\dom \underline{P}$. We ensure that the generated lower prevision has no special properties, for example, 2-monotonicity \citep[Chapter~6]{2014:troffaes:decooman::lower:previsions}) by considering $|\dom \underline{P}| = 2^k$ for $k \in \{4, 6, 8 ,10 \}$ and find that different size of $k$ has an effect on the results (see \cref{fig:plot_domP}).  To investigate this, we also plot the computational time spent in the setup stage of the proposed algorithms (\cref{alg:maximin2,alg:maximin3,alg:maximax2,alg:maximax3}), i.e. finding feasible starting points and sorting gambles.

\begin{figure}
	\centering
	\setlength{\tabcolsep}{2pt}
	\newcolumntype{C}{>{\centering\arraybackslash} m{0.48\linewidth} }
	
	\begin{tabular}{m{0.5em}CC}
		&
		Algorithms for $\Gamma$-maximin
		&
		Algorithms for $\Gamma$-maximax
		\\
		\rotatebox[origin=l]{90}{$|\mathcal{K}| = 2^5$}
		&
		\includegraphics[width=\hsize, trim={1cm 0 2cm 0},clip]{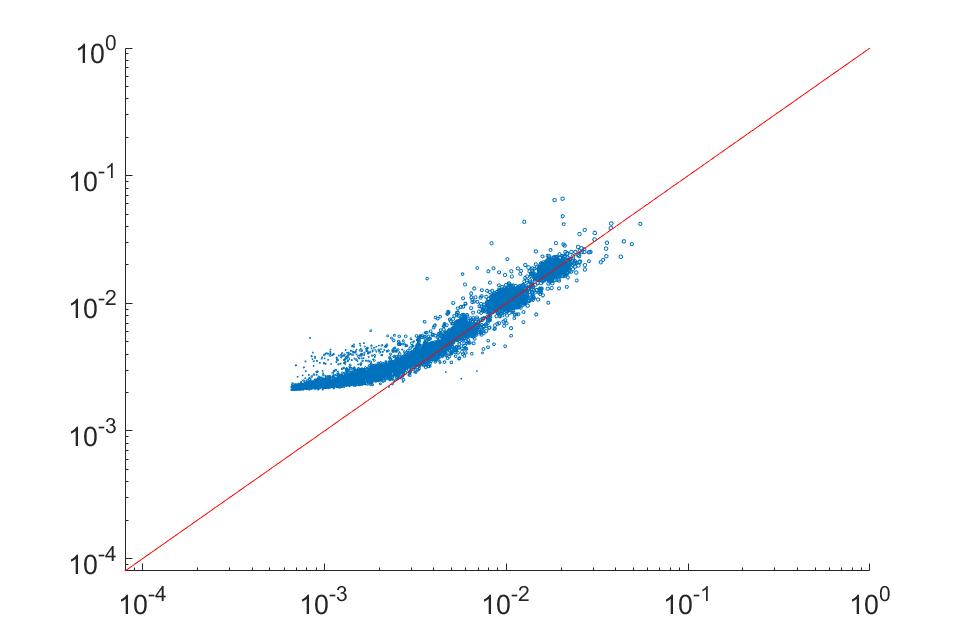}
		&
		\includegraphics[width=\hsize, trim={2cm 0 2cm 0},clip]{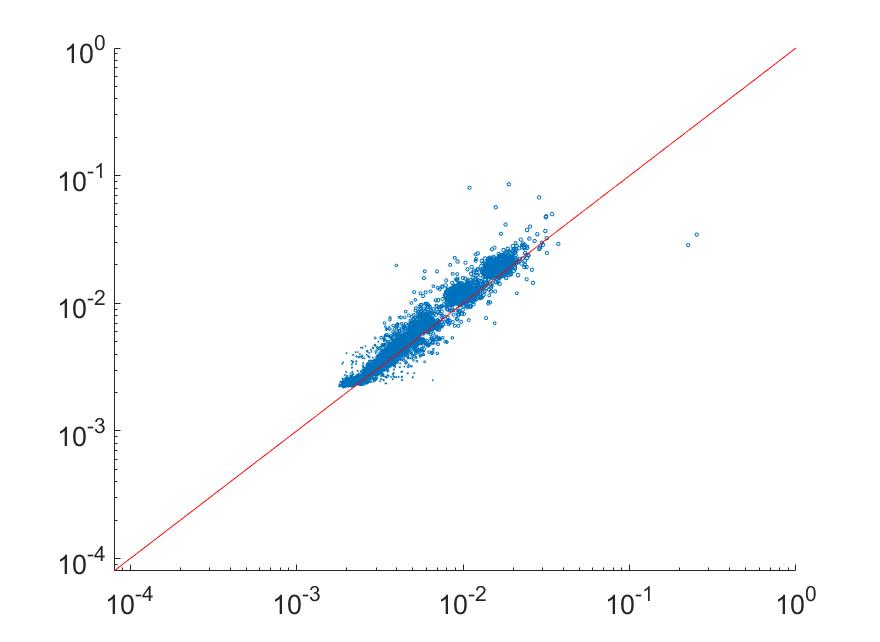}
		\\
		\rotatebox[origin=l]{90}{$|\mathcal{K}| = 2^7$}
		&
		\includegraphics[width=\hsize, trim={1cm 0 2cm 0.5cm},clip]{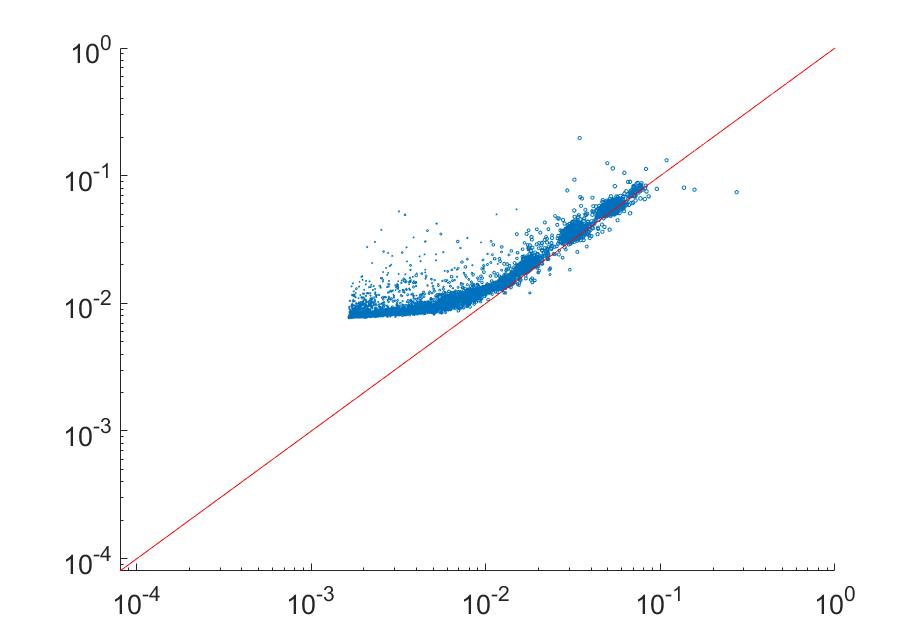}
		&
		\includegraphics[width=\hsize, trim={2cm 0 2cm 0},clip]{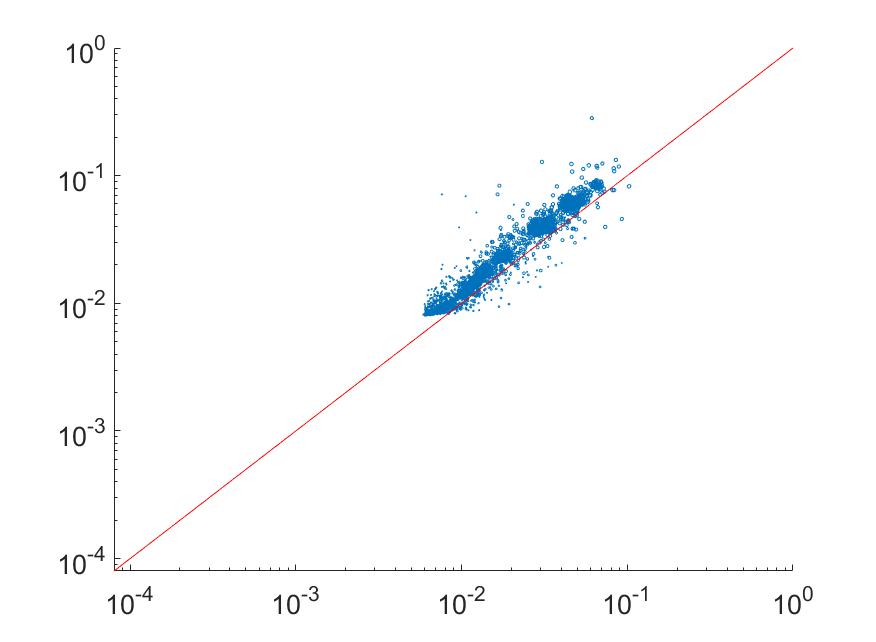}
	\end{tabular}

	\caption{Comparison scatter plots of the running time for different pairs of algorithms for $\Gamma$-maximin (the left column) and $\Gamma$-maximax (the right column).  Each represent a different number of gambles with varying numbers of outcomes and fixed $|\dom \underline{P}| = 2^4$. In the left column, the horizontal axis shows running time for \cref{alg:maximin2} and the vertical axis shows running time for  \cref{alg:maximin3}. In the right column, the horizontal axis shows running time for \cref{alg:maximax2} and the vertical axis shows running time for \cref{alg:maximax3}. The red line is the identity line.}
	\label{fig:scatter_fixGam}
\end{figure}

\begin{figure}
	\centering
	\setlength{\tabcolsep}{2pt}
	\newcolumntype{C}{>{\centering\arraybackslash} m{0.48\linewidth} }
	
	\begin{tabular}{m{0.5em}CC}
		&
		Algorithms for $\Gamma$-maximin
		&
		Algorithms for $\Gamma$-maximax
		\\
		\rotatebox[origin=l]{90}{$|\Omega| = 2^5$}
		&
		\includegraphics[width=\hsize, trim={1cm 0 2cm 0},clip]{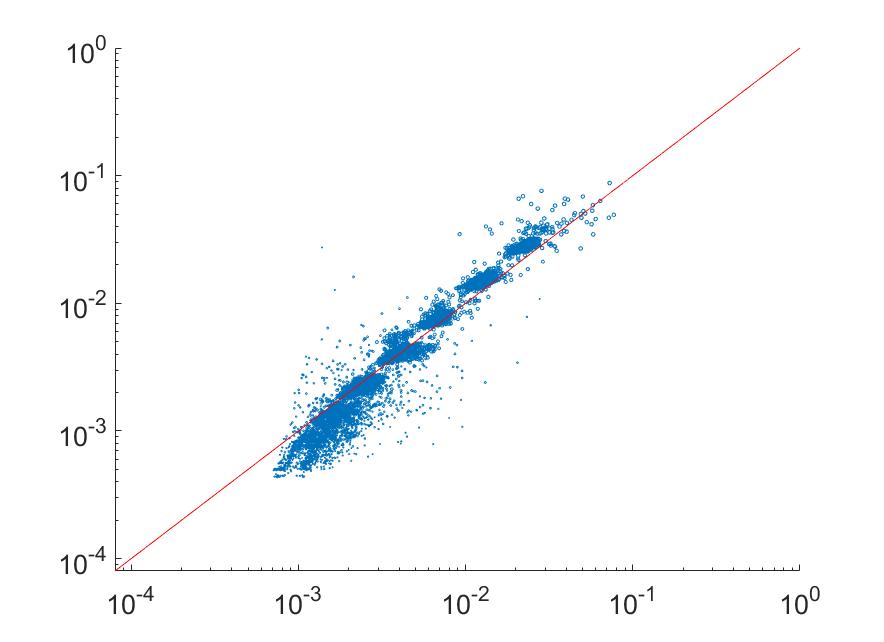}
		&
		\includegraphics[width=\hsize, trim={2cm 0 2cm 0},clip]{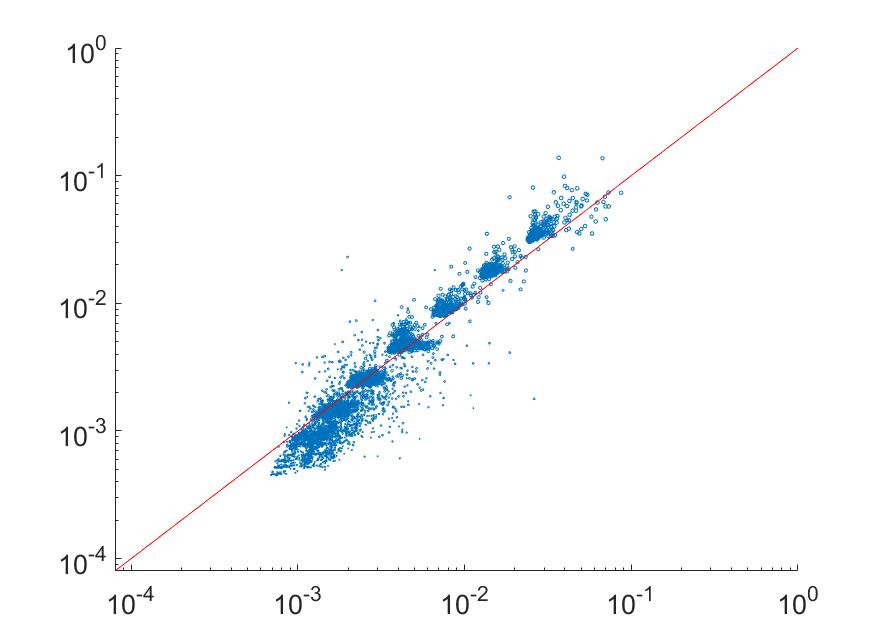}
		\\
		\rotatebox[origin=l]{90}{$|\Omega| = 2^7$}
		&
		\includegraphics[width=\hsize, trim={1cm 0 2cm 0},clip]{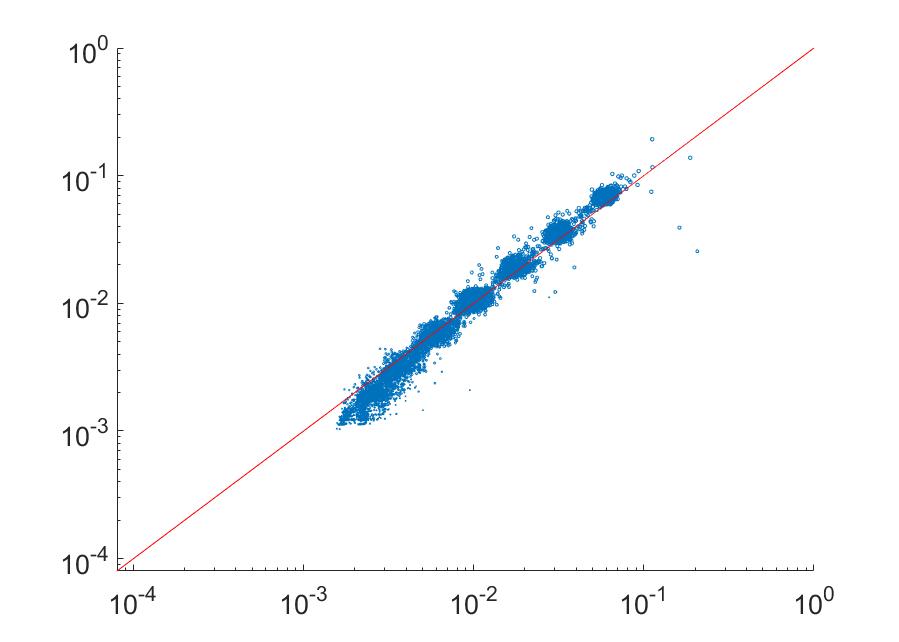}
		&
		\includegraphics[width=\hsize, trim={2cm 0 2cm 0},clip]{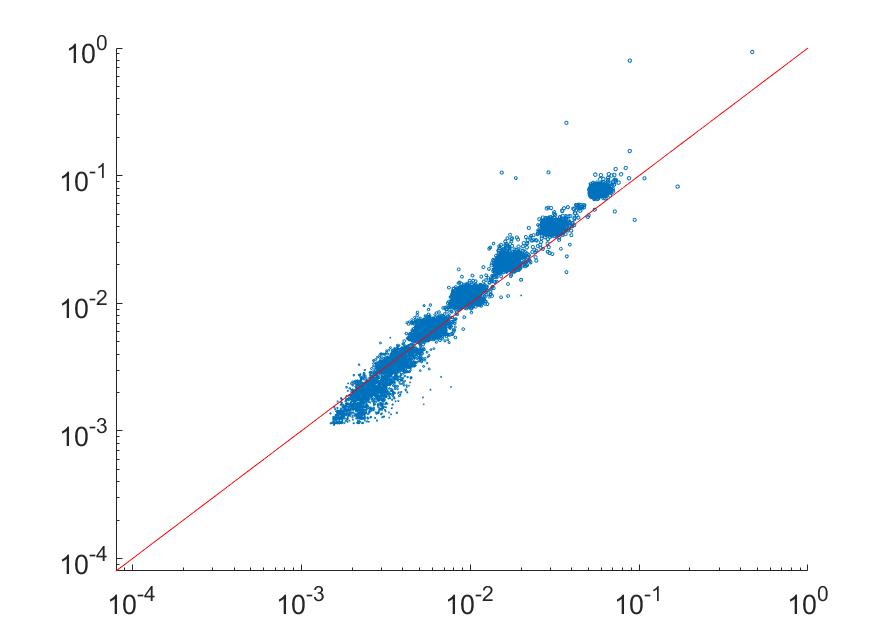}	
	\end{tabular}

	\caption{Comparison scatter plots of the running time for different pairs of algorithms for $\Gamma$-maximin (the left column) and $\Gamma$-maximax (the right column).  Each presents a different number of outcomes with varying numbers of gambles and fixed $|\dom \underline{P}| = 2^4$. In the left column, the horizontal axis shows running time for \cref{alg:maximin2} and the vertical axis shows running time for  \cref{alg:maximin3}. In the right column, the horizontal axis shows running time for \cref{alg:maximax2} and the vertical axis shows running  time for \cref{alg:maximax3}. The red line is the identity line.}
\label{fig:scatter_fixOut}
\end{figure}

We find that in \cref{fig:plot1,fig:plot2,fig:plot_domP}, the running times for finding $\Gamma$-maximin by \cref{alg:maximin2,alg:maximin3} are nearly identical. Similarly, for \cref{alg:maximax2,alg:maximax3}, we see that the running time for finding $\Gamma$-maximax gambles is very similar.
In addition to presenting average running times and error bars, we also provide some scatter plots of running time for finding $\Gamma$-maximin and $\Gamma$-maximax gambles in the case that they are very closed, i.e. comparing \cref{alg:maximin2} against \cref{alg:maximin3} and comparing \cref{alg:maximax2} against \cref{alg:maximax3} in \cref{fig:scatter_fixGam,fig:scatter_fixOut}, respectively. 
This shows any simulations where one algorithm has performed better than another or perhaps whether the running times under the two algorithms are related to each other or not.

\Cref{fig:scatter_fixGam,fig:scatter_fixOut} show scatter plots of running time for different pairs of algorithms for finding $\Gamma$-maximin gambles in the left column and $\Gamma$-maximax gambles in the right column. In \cref{fig:scatter_fixGam}, each row shows a different number of gambles and in \cref{fig:scatter_fixOut}, each row presents a different number of outcomes. In the left column where we compare running time taken for finding $\Gamma$-maximin gamble. The horizontal axis presents the running time of \cref{alg:maximin2} and the vertical axis presents the running time of \cref{alg:maximin3}. In the right column, we compare running time for finding $\Gamma$-maximax gambles, where the horizontal axis presents the running time of \cref{alg:maximax2} and the vertical axis presents the running time of \cref{alg:maximax3}.

\subsection{Interval dominance}

To benchmark the algorithms for interval dominance (\cref{alg:ID1,alg:ID2,alg:ID,alg:ID3}), we will generate random sets of gambles that have a precise number of $\Gamma$-maximin and interval dominant gambles by using \citep[algorithm 3]{2020:IUKM:Nakharutai}. We set the number of outcomes to be $|\Omega| = 2^4$ and $|\Omega| = 2^6$, the number of gambles in $\mathcal{K}$ to be $2^4, 2^6$ and $2^8$, and the size of $\dom \underline{P}$ to be $2^4$ and $2^6$.

\begin{figure}[p]
	\centering
	\begin{tikzpicture}
		\filldraw[draw=black, fill=gray!20] (2,0) -- (8,4.5) --(8,0);
		\draw (2,-1) --(2,5.5)node[anchor=east]{$m$};
		\draw (2,4.5) --(9.5,4.5);
		\draw [->](2,5.5) -- (2,6);
		\draw (6,3) --(8,3);
		\draw (4,1.5) -- (8,1.5);
		\draw (2,0)--(8,4.5);
		\draw (8,0)-- (8,1.5)node[circle,draw=black,fill=white!80!black,minimum size=20]{$g$} --(8,3)node[circle,draw=black,fill=white!80!black,minimum size=20]{$i$}--(8,4.5)node[circle,draw=black,fill=white!80!black,minimum size=20]{$j$} --(8,5.5);
		\draw (8.5,5.2)node[anchor=west]{$k$};
		\draw(8,5) -- (8,6);
		
		\draw (6,0)-- (6,1.5)node[circle,draw=black,fill=white!80!black,minimum size=20]{$f$} --(6,3)node[circle,draw=black,fill=white!80!black,minimum size=20]{$h$};
		
		\draw (4,0)-- (4,1.5)node[circle,draw=black,fill=white!80!black,minimum size=20]{$e$};
		
		\draw (1,0) -- (2,0)node[circle,draw=black,fill=white!80!black,minimum size=20]{$a$} --(4,0)node[circle,draw=black,fill=white!80!black,minimum size=20]{$b$}--(6,0)node[circle,draw=black,fill=white!80!black,minimum size=20]{$c$}--(8,0)node[circle,draw=black,fill=white!80!black,minimum size=20]{$d$} --(10,0)node[anchor=north]{$n$};
		\draw [->](10,0) -- (11,0);
	\end{tikzpicture}
	\caption{The area of $\ell \leq n \leq k$ and 10 options label the different $\ell$ and $n$ that we consider in the simulation (see \cref{table:10cases})}\label{fig:ell_n_k}
\end{figure}
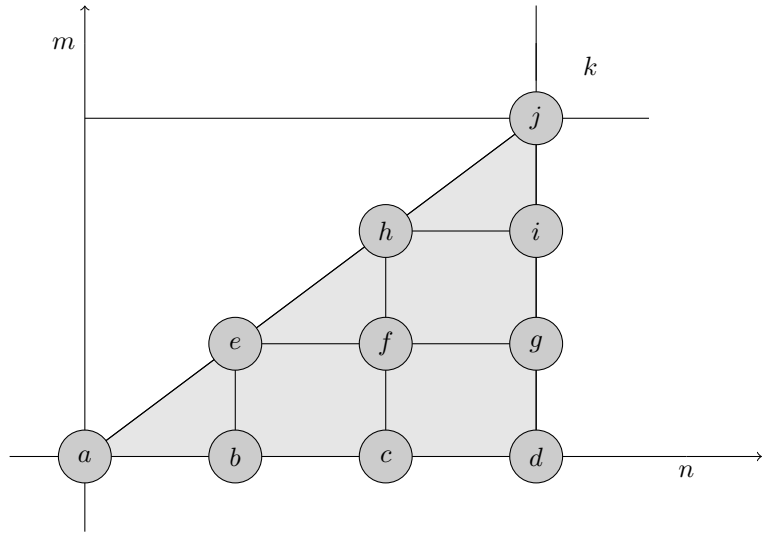

To do so, we first generate $\underline{P}$ with size of $\dom\underline{P}$  as we did before in benchmarking for $\Gamma$-maximin and $\Gamma$-maximax. Next, we apply \citep[algorithm 3]{2020:IUKM:Nakharutai} to generate a random set $\mathcal{K}$ for $|\mathcal{K}|=k\in\{2^4, 2^6, 2^8\}$, $|\opt_{\underline{E}}(\mathcal{K})|= \ell$ and $	|\opt_{\sqsupset}(\mathcal{K})|=n$ where $\ell \leq n\leq k$ by using the previous $\underline{P}$ to evaluate $\underline{E}$ and $\overline{E}$.
To cover a range of possible options of $\ell$, $n$, and $k$ that satisfy $\ell \leq n \leq k$, we follows different options with respect to the size of $\mathcal{K}$ as suggested in \citep[Table 1]{2019:Nakharutai:Troffase:Caiado:maximal}; see \cref{table:10cases} and \cref{fig:ell_n_k} that illustrate these cases. Specifically, options a to d represent the sets of gambles with $\ell=1$ while we increase $n$ from $1$ to $k$. Options d, g, i, and j represent the sets of gambles with $n=k$ while we increase $\ell$ from $1$ to $k$. Options a, e, h, and j represent the sets of gambles that $\ell=n$ while we increase them jointly from $1$ to $k$. Option f represents a set of gambles where $\ell<n<k$.

{
	\setlength{\arrayrulewidth}{0.5mm}
	\setlength{\tabcolsep}{12pt}
	\renewcommand{\arraystretch}{1.3}
	\begin{table}[p]
		\centering
		\begin{tabular}{|c|c|cV{1.5}c|cV{1.5}c|c|}
			\hline
			\multirow{2}{*}{Options} & \multicolumn{2}{lV{1.5}}{$|\mathcal{K}| = 2^4$} & \multicolumn{2}{lV{1.5}}{$|\mathcal{K}| = 2^6$} & \multicolumn{2}{l|}{$|\mathcal{K}| = 2^8$} \\ \cline{2-7} 
			&  $m$  &  $n$  &     $m$  &  $n$     &   $m$  &  $n$ \\ \hline
			\circled{a} &  1 & 1 &  1   & 1 &  1 & 1\\ \hline
			\circled{b} &  1 & 5  &  1 &  21 &   1 & 85\\ \hline
			\circled{c} &  1 & 11 & 1  &  42 &  1  &170 \\ \hline
			\circled{d} &  1 & 16 &  1 & 64  &  1  &256 \\ \hline
			\circled{e} &  5 & 5  &  21 &  21 &   85 &85 \\ \hline
			\circled{f} &  5 & 11 &  21 &  42 &  85  &170 \\ \hline
			\circled{g} &  5 & 16 &  21 & 64  &   85 & 256\\ \hline
			\circled{h} &  11 & 11 &  42 &  42 &  170  &170 \\ \hline
			\circled{i} &  11 & 16 &  42 &  64 &   170 & 256\\ \hline
			\circled{j}&  16 & 16  &  64 & 64  &  256  & 256\\ \hline
		\end{tabular}
		\caption{Table of points that indicate different sizes of set $\mathcal{K}$ with varying number of $\Gamma$-maximin gambles $\ell$ and number of interval dominant gambles $n$ in $\mathcal{K}$} \label{table:10cases}
	\end{table}
}

\begin{figure}
	\centering
	\setlength{\tabcolsep}{2pt}
	\newcolumntype{C}{>{\centering\arraybackslash} m{0.48\linewidth} }
	
	\begin{tabular}{m{0.5em}CC}
		&
		$|\Omega| = 2^4$
		&
		$|\Omega| = 2^6$
		\\
		\rotatebox[origin=l]{90}{$|\mathcal{K}| = 2^4$}
		&
		\includegraphics[width=\hsize, trim={2cm 0 2cm 0},clip]{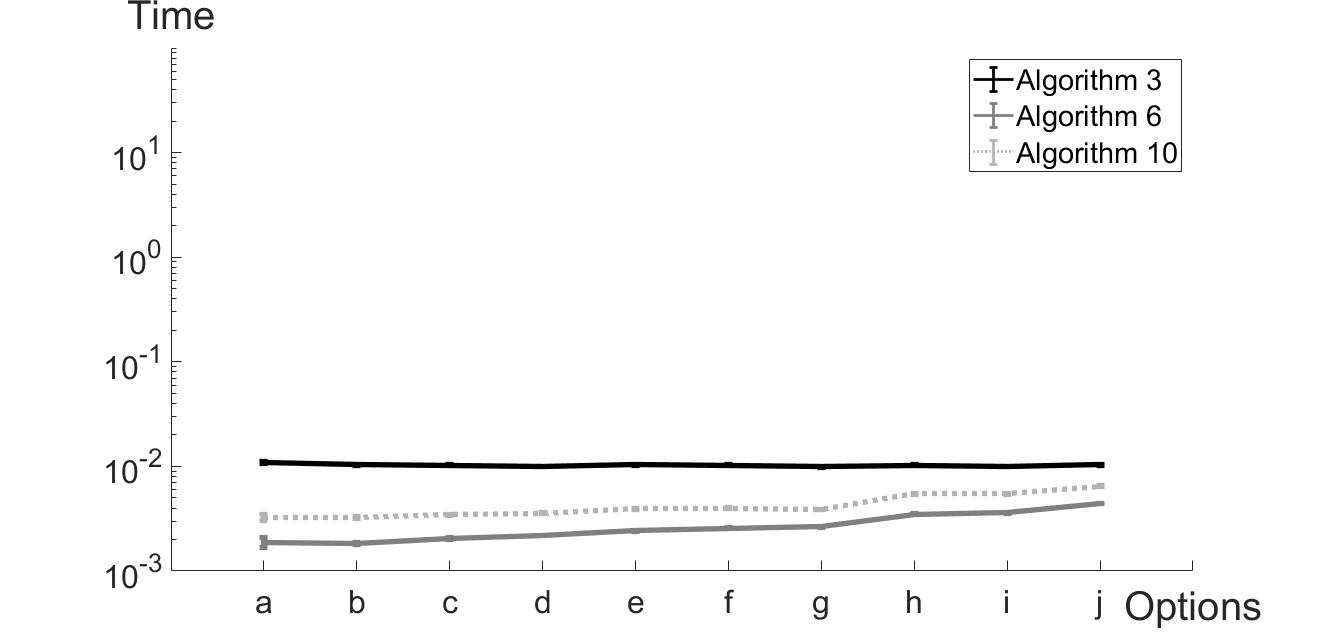}
		&
		\includegraphics[width=\hsize, trim={2cm 0 2cm 0},clip]{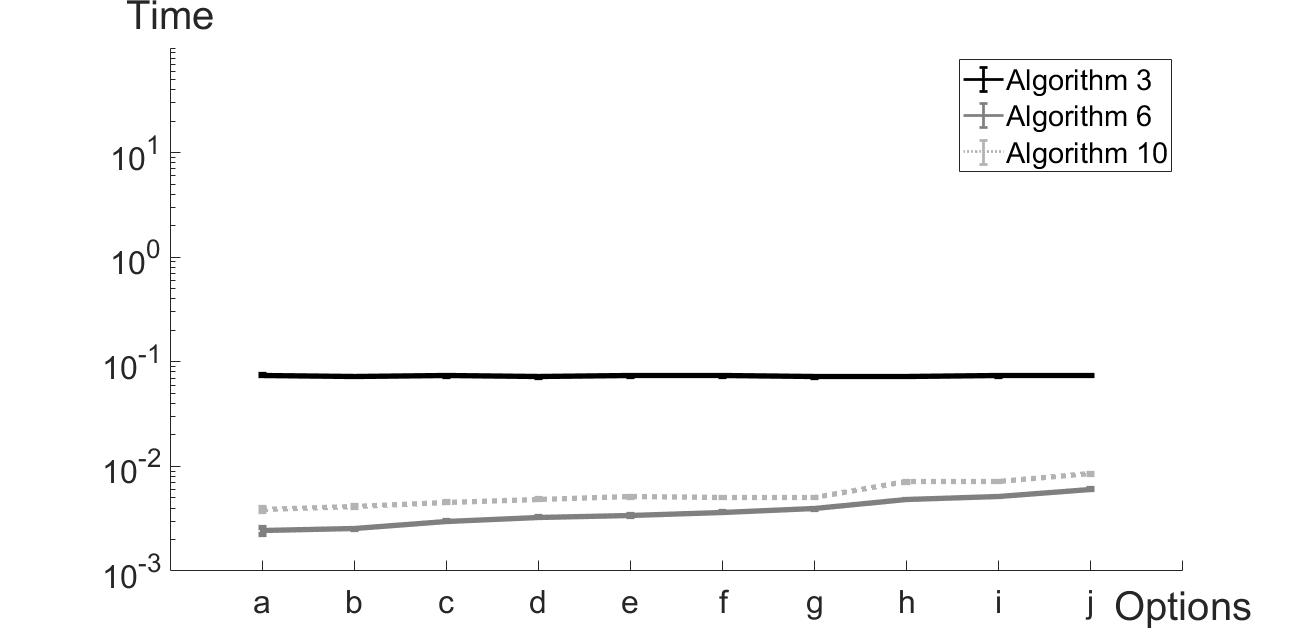}
		\\
		\rotatebox[origin=l]{90}{$|\mathcal{K}| = 2^6$}
		&
		\includegraphics[width=\hsize, trim={2cm 0 2cm 0},clip]{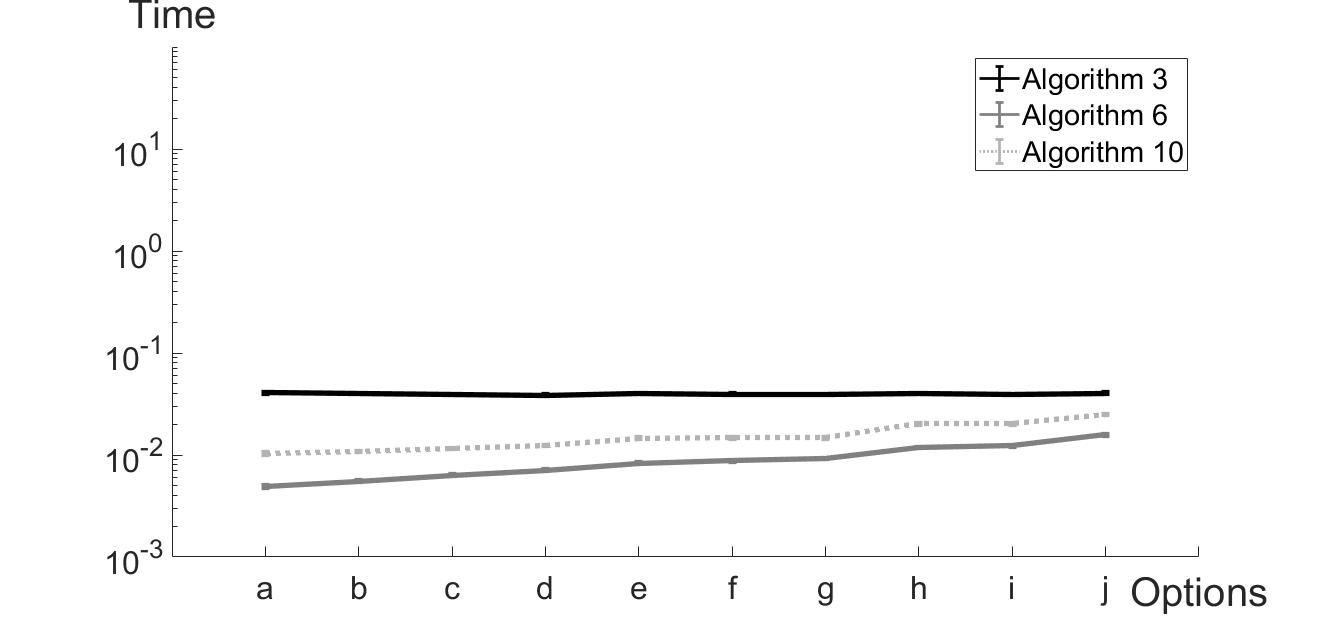}
		&
		\includegraphics[width=\hsize, trim={2cm 0 2cm 0},clip]{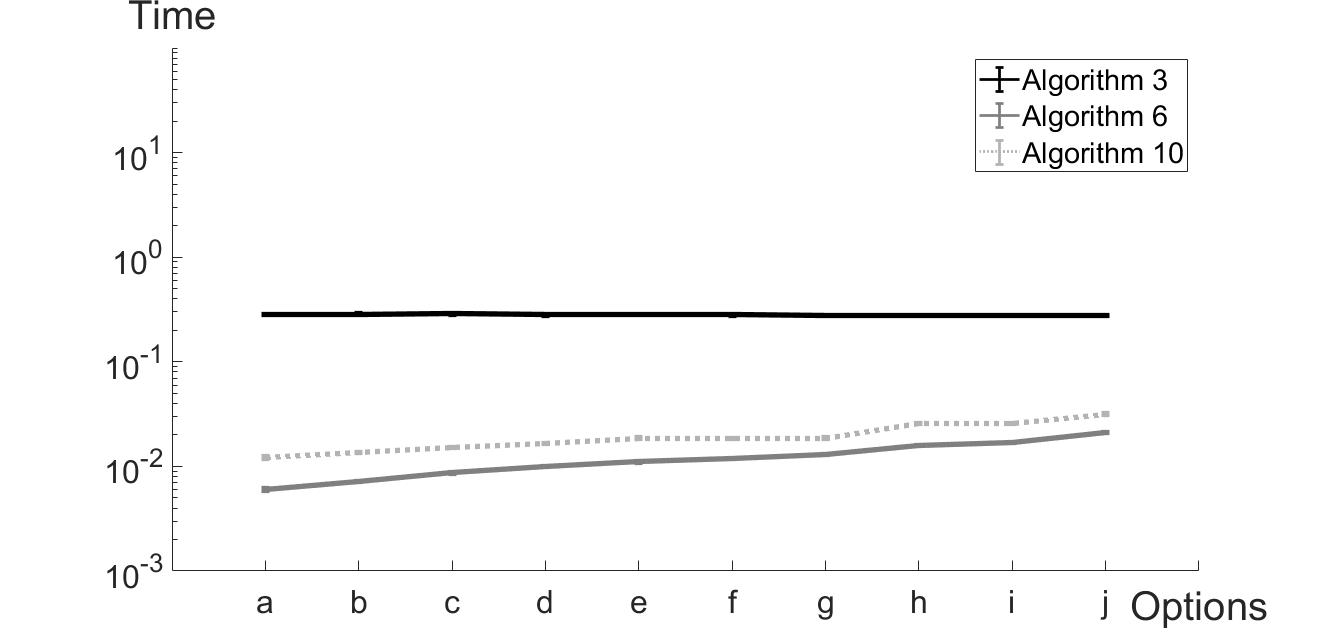}
		\\
		\rotatebox[origin=l]{90}{$|\mathcal{K}| = 2^8$}
		&
		\includegraphics[width=\hsize, trim={2cm 0 2cm 0},clip]{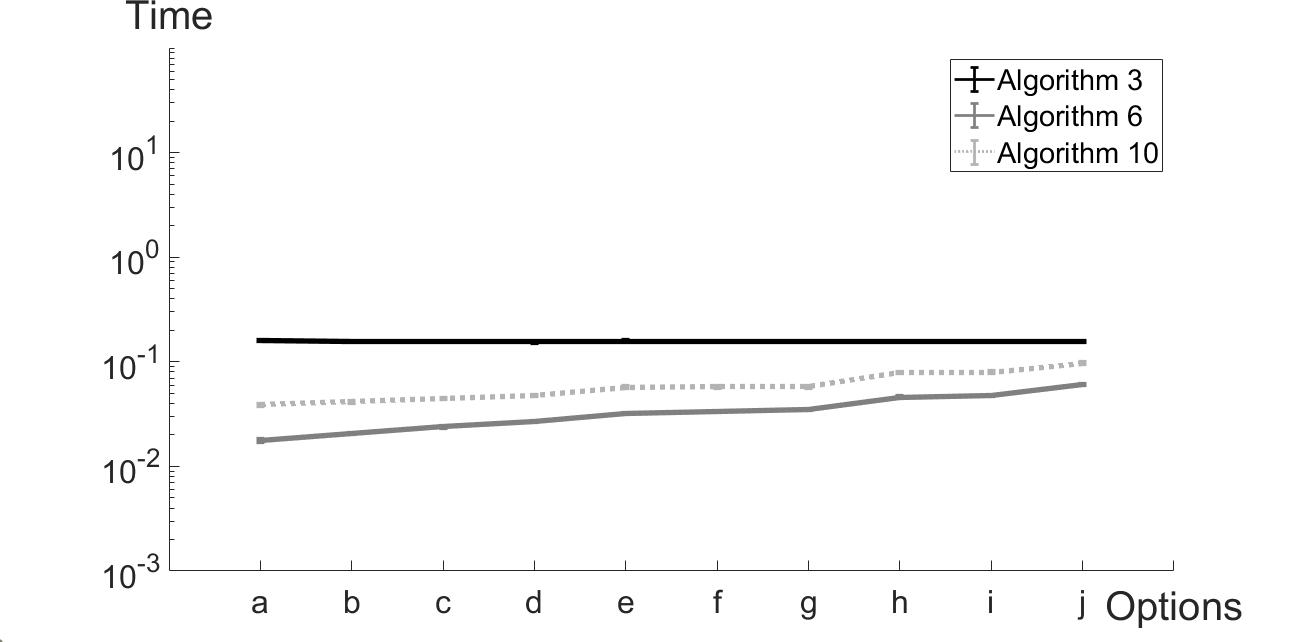}
		&
		\includegraphics[width=\hsize, trim={2cm 0 2cm 0},clip]{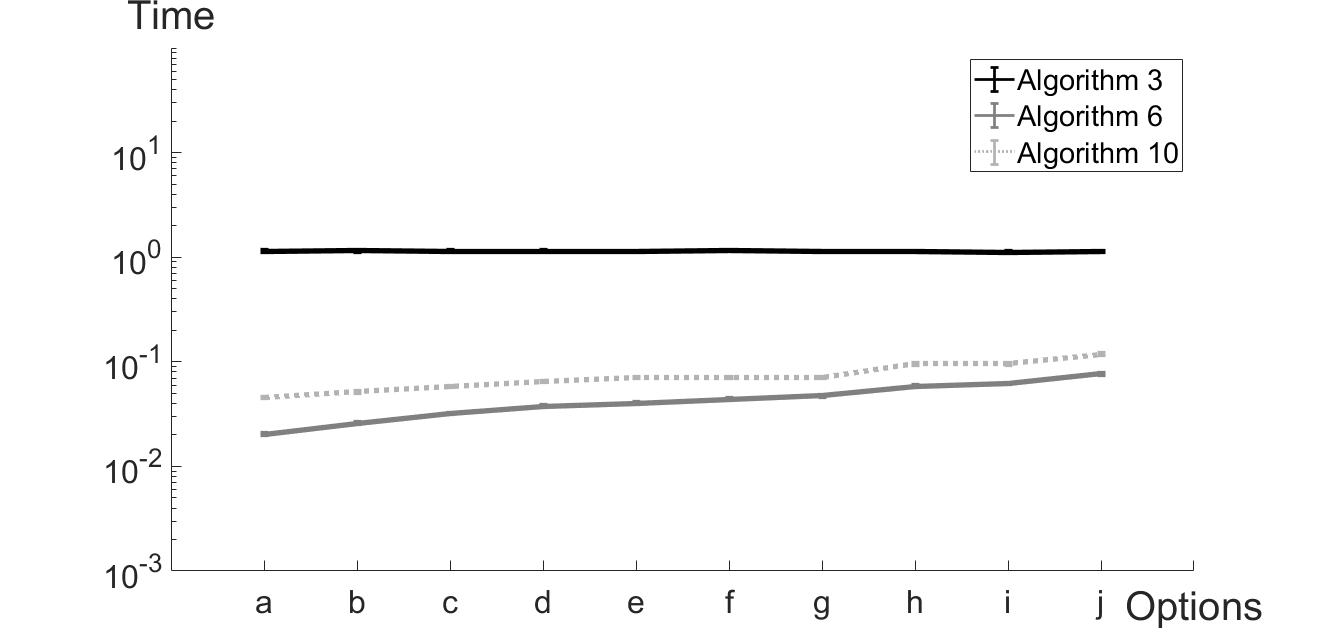}
	\end{tabular}
	\caption{Comparison plots of the average computational time for different improved algorithms for finding interval dominant gambles where we fix $|\dom \underline{P}| = 2^4$. The number of outcomes in left column is $2^4$ and $2^6$ in the right column.  Each row represents a different number of gambles with varying options of the numbers of $\Gamma$-maximin gambles and interval dominant gambles in the set (see \cref{table:10cases} for each option). The labels indicate different algorithms.}
	\label{fig:plot:ID:domP4}
\end{figure}

\begin{figure}
	\centering
	\setlength{\tabcolsep}{2pt}
	\newcolumntype{C}{>{\centering\arraybackslash} m{0.48\linewidth} }
	
	\begin{tabular}{m{0.5em}CC}
		&
		$|\Omega| = 2^4$
		&
		$|\Omega| = 2^6$
		\\
		\rotatebox[origin=l]{90}{$|\mathcal{K}| = 2^4$}
		&
		\includegraphics[width=\hsize, trim={2cm 0 2cm 0},clip]{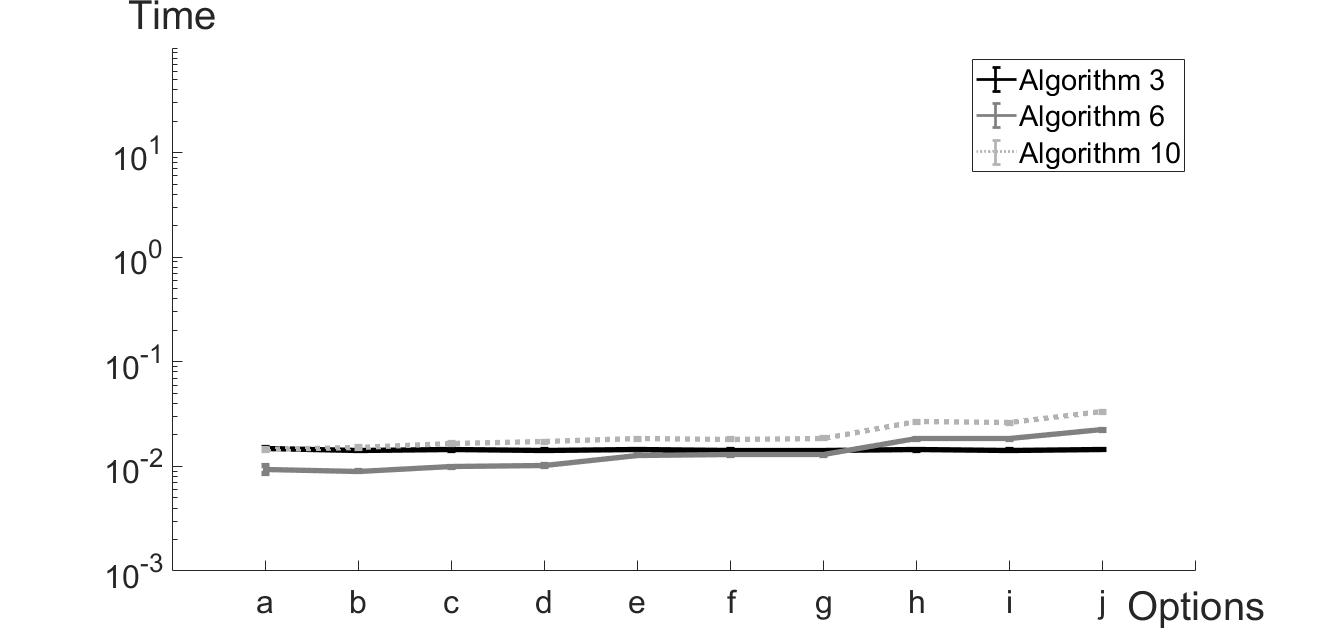}
		&
		\includegraphics[width=\hsize, trim={2cm 0 2cm 0},clip]{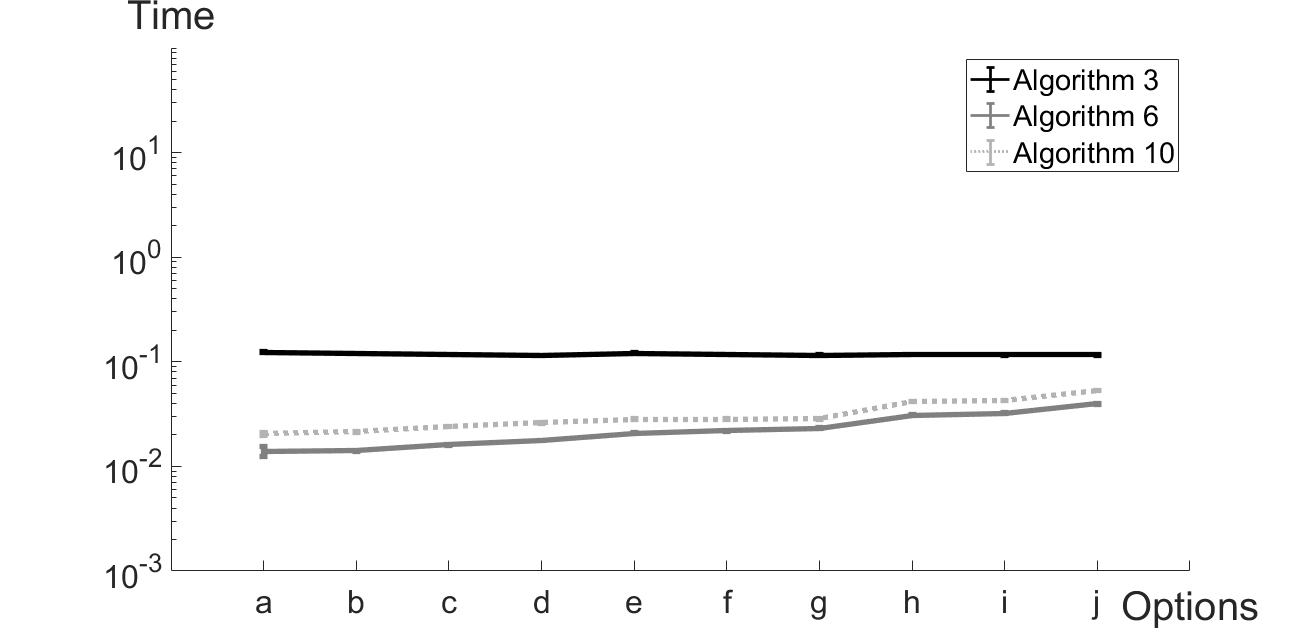}
		\\
		\rotatebox[origin=l]{90}{$|\mathcal{K}| = 2^6$}
		&
		\includegraphics[width=\hsize, trim={2cm 0 2cm 0},clip]{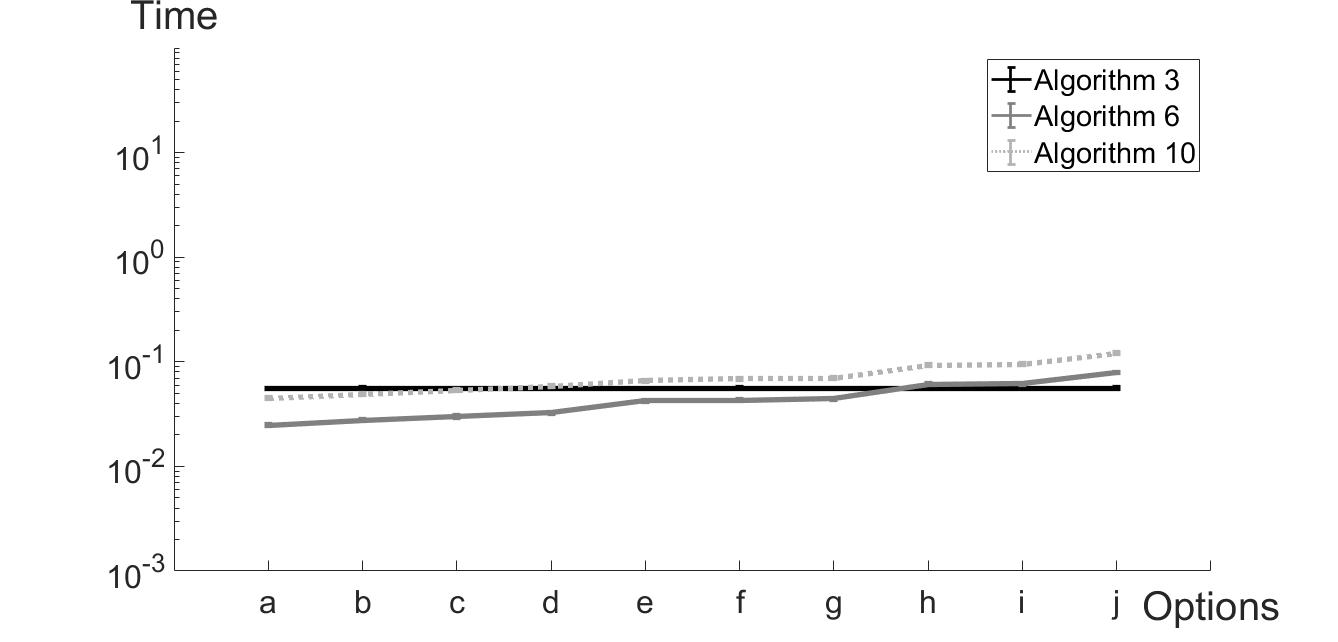}
		&
		\includegraphics[width=\hsize, trim={2cm 0 2cm 0},clip]{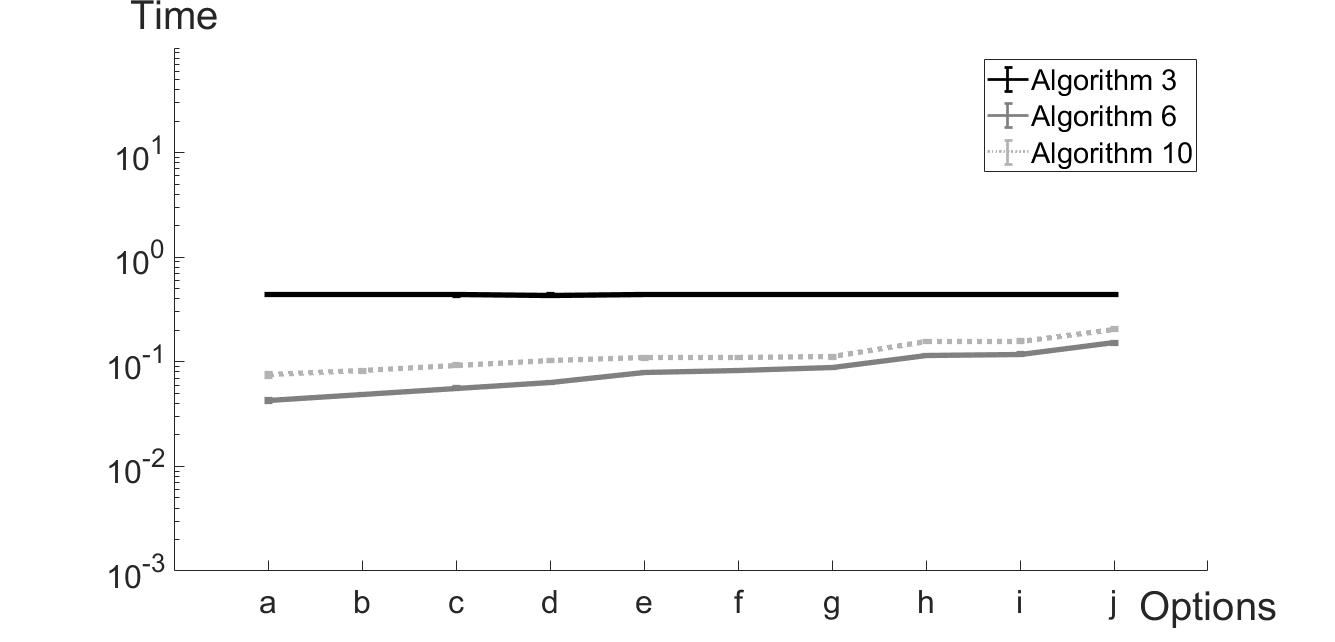}
		\\
		\rotatebox[origin=l]{90}{$|\mathcal{K}| = 2^8$}
		&
		\includegraphics[width=\hsize, trim={2cm 0 2cm 0},clip]{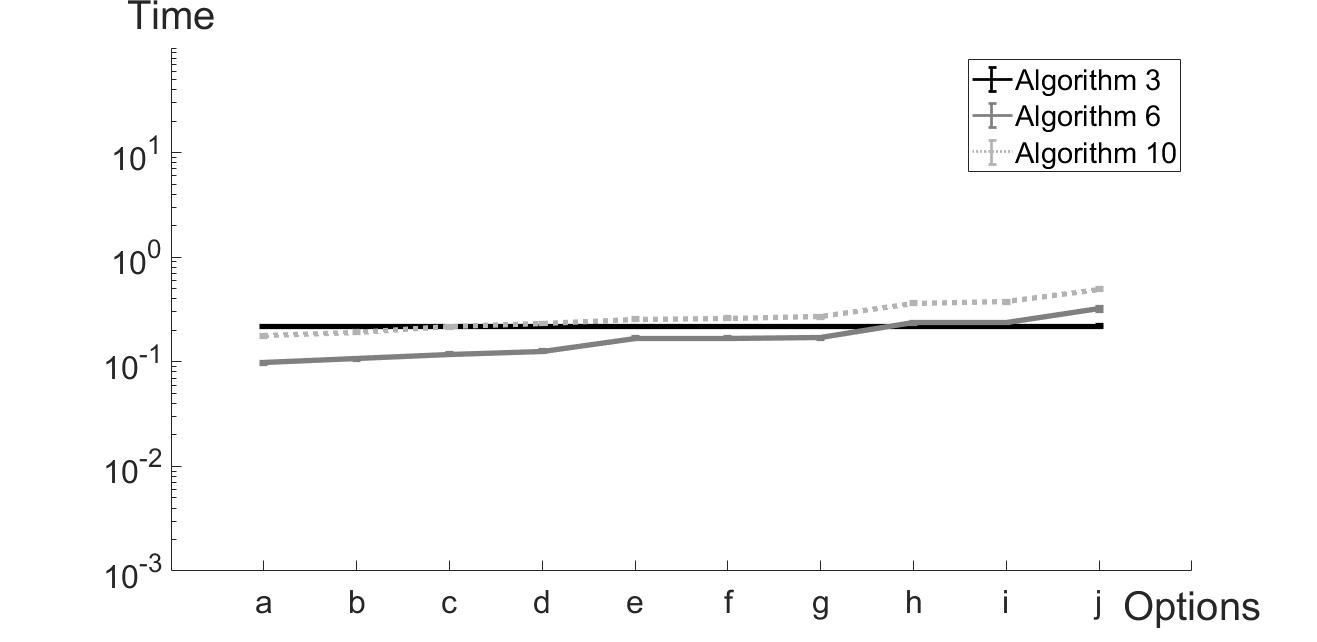}
		&
		\includegraphics[width=\hsize, trim={2cm 0 2cm 0},clip]{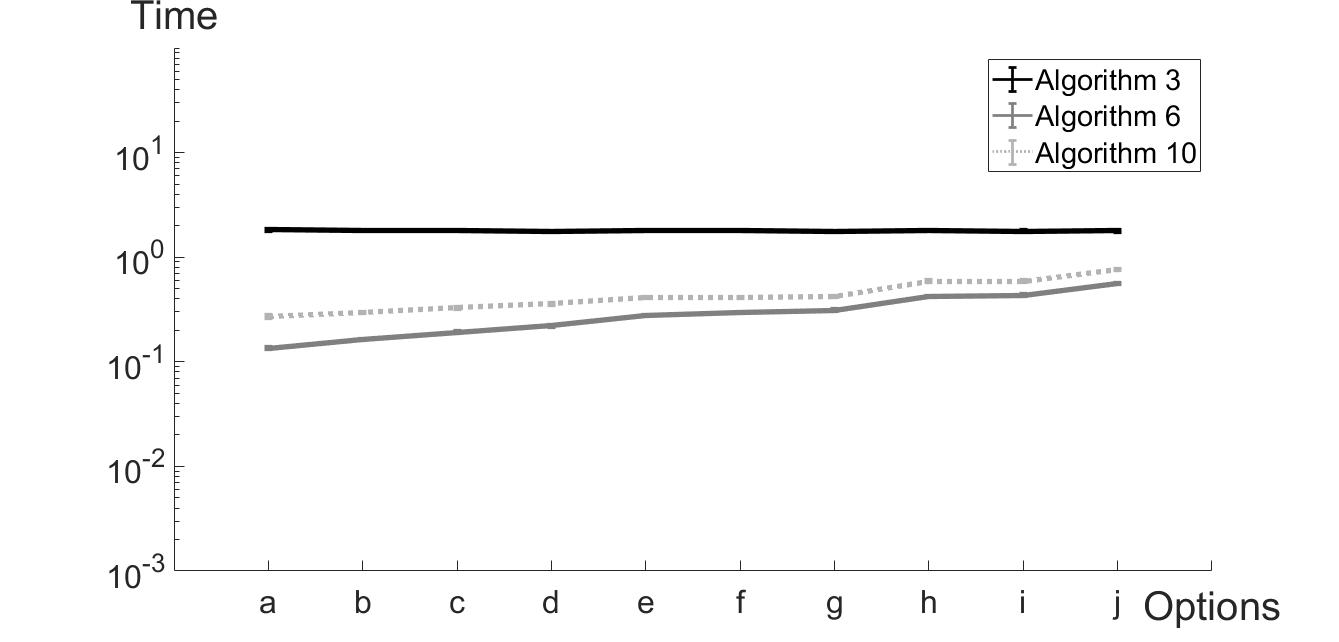}
	\end{tabular}
	\caption{Comparison plots of the average computational time for different improved algorithms for finding interval dominant gambles  where we fix $|\dom \underline{P}| = 2^6$. The number of outcomes in left column is $2^4$ and $2^6$ in the right column. Each row represents a different number of gambles with varying options of the numbers of $\Gamma$-maximin gambles and interval dominant gambles in the set (see \cref{table:10cases} for each option). The labels indicate different algorithms.}
	\label{fig:plot:ID:domP6}
\end{figure}

For each generated set of gambles $\mathcal{K}$, we apply \cref{alg:ID1,alg:ID2,alg:ID,alg:ID3} and record the total running time. For \cref{alg:ID1}, we solve linear programs by the standard primal-dual method while for \cref{alg:ID2,alg:ID,alg:ID3}, we solve the linear programs by the improved primal-dual method which has all improvements from \cref{subsec:improve_LP}. Again, to ensure a fair comparison, these algorithms are written in MATLAB, and use identical implementations for solving the linear programs. For sorting gambles we used the built-in quicksort function. For each case, we repeat the process 500 times and present summaries in \cref{fig:plot:ID:domP4,fig:plot:ID:domP6}.

\Cref{alg:ID} is not depicted in the plots,
as we have found it to perform nearly identically to \cref{alg:ID3}. This is because there is a trade-off between effort spent on bounding the $\Gamma$-maximin value, and the number of comparisons required against these bounds, as discussed in \cref{sec:furtherimprovements}.

\Cref{fig:plot:ID:domP4,fig:plot:ID:domP6} present the average total running time of three algorithms for interval dominance. In the left column, the number of outcomes is $2^4$ while the right column, the number of outcomes is $2^6$. Each row represents a different number of gambles in the set $\mathcal{K}$. The vertical axis shows the average computational time over 500 random generated sets and the horizontal axis shows the different ten options of $\ell$, $n$ and $k$ that we gave in \cref{table:10cases}. The error bars representing 95\% confidence intervals on the total running time are barely visible here.

\section{Discussion and conclusion}\label{sec:conclusion}
In this study, we proposed new algorithms for $\Gamma$-maximin (\cref{alg:maximin2,alg:maximin3}), for $\Gamma$-maximax (\cref{alg:maximax2,alg:maximax3}) and for interval dominance (\cref{alg:ID2,alg:ID,alg:ID3}) and compared their performance with the standard algorithms (\cref{alg:maximin1,alg:maximax1,alg:ID1} for $\Gamma$-maximin, for $\Gamma$-maximax and for interval dominance respectively). 

To find $\Gamma$-maximin and $\Gamma$-maximax gambles, the algorithms solve linear programs for each gamble in a set of gambles. For finding interval dominant gambles, the algorithms first need to solve a sequence of linear programs and then solve linear programs for each gamble in a set of gambles.
Applying improvements from \citep{2019:Nakharutai:Troffase:Caiado:maximal}, 
for each case, we gave early stopping criteria and a quick way to compute feasible starting points for linear programs. We found that the primal-dual method can benefit from these improvements, as suggested in \citep{2019:Nakharutai:Troffase:Caiado:maximal}.
Our proposed algorithms for $\Gamma$-maximin, $\Gamma$-maximax and interval dominance implemented these improvements.

For benchmarking these improvements, we used \citep[algorithms 2 and 4]{2018:Nakharutai:Troffaes:Caiado} to generate sets of gambles and applied \cref{alg:maximin1,alg:maximin2,alg:maximin3} for $\Gamma$-maximin and  \cref{alg:maximax1,alg:maximax2,alg:maximax3} for $\Gamma$-maximax on these generated sets. For benchmarking algorithms for interval dominance, we used \citep[algorithm 3]{2020:IUKM:Nakharutai} to generate sets of gambles with a precise number of $\Gamma$-maximin and interval dominant gambles  and applied \cref{alg:ID1,alg:ID2,alg:ID,alg:ID3}.
According to our simulation, the numbers of outcomes and gambles in the set have an impact on the performance of these different algorithms. Specifically,
the mean running time spent on the algorithm is larger if 
the numbers of outcomes or the number of gambles is increasing. The mean running time also depends on the size of the domain of the lower previsions.

Overall, both \cref{alg:maximin2,alg:maximin3} outperform \cref{alg:maximin1} for finding $\Gamma$-maximin gambles and both \cref{alg:maximax2,alg:maximax3} outperform \cref{alg:maximax1} for finding $\Gamma$-maximax gambles especially when either the number of gambles or the number of outcomes is large. There is no big difference in the time taken on \cref{alg:maximax2,alg:maximax3} for finding $\Gamma$-maximax gambles. For the time taken on \cref{alg:maximin2,alg:maximin3}, there is also no big difference except when the number of gambles is large and the number of outcome is small (see \cref{fig:plot1}). In this specific case, \cref{alg:maximin2} performs slightly better than \cref{alg:maximin3}.

Nevertheless, when we increased the size of the domain of lower previsions (see \cref{fig:plot_domP}), the result showed that the size of the domain has an impact on the setup stage of the proposed algorithms. We found that if the size of the domain is much larger than the sizes of gambles and outcomes, then \cref{alg:maximin1} performs better than  \cref{alg:maximin2,alg:maximin3} while \cref{alg:maximax1} performs better than \cref{alg:maximax2,alg:maximax3}. This is because \cref{alg:maximin2,alg:maximin3,alg:maximax2,alg:maximax3} spent longer time in the setup stage in order to get feasible starting points. However, in the case that the size of the domain is less or equal to the sizes of gambles and outcomes, our proposed \cref{alg:maximin2,alg:maximin3} outperform \cref{alg:maximin1} and our proposed \cref{alg:maximax2,alg:maximax3} outperform \cref{alg:maximax1} as there is little effect from the setup stage in these proposed algorithms.

We note that \cref{alg:maximin3} is the only algorithm that can return all $\Gamma$-maximin gambles while the other algorithms return just one.
Therefore, if we suspect that there might be more than one $\Gamma$-maximin gamble in the set, and we are interested in finding them all, then we should apply \cref{alg:maximin3}.

For the result of algorithms for interval dominance, we find that our improved algorithms \cref{alg:ID2,alg:ID,alg:ID3} outperform \cref{alg:ID1} in most cases except when the size of the domain of lower previsions is much larger than the sizes of gambles and outcomes where our improved algorithms need to take time for the setup stage. If either the number of $\Gamma$-maximin or interval dominant gambles is increasing, then the computational time on \cref{alg:ID2,alg:ID,alg:ID3} is slightly longer while \cref{alg:ID1} remains the same. \Cref{alg:ID2} slightly performs better than \cref{alg:ID,alg:ID3}.

Among the three proposed algorithms for interval dominance, we found that there was an interesting trade-off between being precise and the number of comparisons. In particular, if we aim to be precise at the beginning, that is, if we obtain a precise $\Gamma$-maximin value, then we require a smaller number of comparisons later on. On the other hand, if we allow being less precise at the beginning and start comparison immediately, then we could save computational effort as we may not need to identify the precise $\Gamma$-maximin value, at the expense of having to do more comparisons. It turns out that our new algorithm for finding the $\Gamma$-maximin value (\cref{alg:maximin2}) is fast enough to tip the balance in clear favour of \cref{alg:ID2}.

Based on our benchmarking results from all scenarios considered, both of our proposed algorithms are good choices for implementation when the size of the domain is less or equal to the sizes of gambles and outcomes, as they outperform the standard algorithms, and otherwise perform quite similarly for the most part. However, when the size of the domain is much larger than the sizes of gambles and outcomes, the proposed algorithms may not be recommended due to the effect from the setup stage.

\section*{Acknowledgements}
This work is supported by CMU Junior Research Fellowship Program. 

\bibliographystyle{plainnat}
\bibliography{references}

\end{document}